\documentclass[12pt]{amsart}

\usepackage{amscd}
\usepackage{amssymb}
\usepackage{amsmath}
\usepackage[dvips]{graphicx}
\usepackage{amsfonts}
\usepackage{amsopn}
\usepackage{amsthm}

\newcounter{mtheorem}

\newtheorem{theorem}{Theorem}[section]
\newtheorem{lemma}[theorem]{Lemma}
\newtheorem{prop}[theorem]{Proposition}
\newtheorem{corollary}[theorem]{Corollary}

\theoremstyle{definition}
\newtheorem{definition}[theorem]{Definition}
\newtheorem{example}[theorem]{Example}

\theoremstyle{remark}
\newtheorem{remark}[theorem]{Remark}
\numberwithin{equation}{section}
\newcommand{\C}{\mathbb{C}}

\newcommand{\R}{\mathbb{R}}

\newcommand{\dist}{\operatorname{dist}}

\newcommand{\ccinf}{C^\infty_c}

\newcommand{\Rd}{\mathbb{R}^d}
\newcommand{\RD}{\mathbb{R}^D}
\newcommand{\RDD}{\mathbb{R}^D \times \mathbb{R}^D}
\newcommand{\Rtwod}{\mathbb{R}^{2d}}
\newcommand{\RnD}{\mathbb{R}^{nD}}

\newcommand{\G}{\mathcal{G}}
\newcommand{\Gmu}{\mathcal{G}_\mu}
\newcommand{\calK}{\mathcal{K}}
\newcommand{\calL}{\mathcal{L}}

\newcommand{\M}{\mathcal{M}}
\newcommand{\N}{(\ccinf)^*}

\newcommand{\calS}{\mathbb{S}}

\newcommand{\xm}{\mathcal{X}(M)}

\newcommand{\cxrD}{\mathcal{X}_c(\R^D)}

\newcommand{\sympx}{\mbox{Symp}\,\mathcal{X}}
\newcommand{\hamfields}{\mbox{Ham}\,\mathcal{X}}
\newcommand{\chamfields}{\mbox{Ham}\,\mathcal{X}_c}
\newcommand{\lhamfields}{\overline{Ham\,\mathcal{X}_c}^\mu}

\newcommand{\fields}{\mathcal{X}}
\newcommand{\cfields}{\mathcal{X}_c}
\newcommand{\lfields}{L^2(\mu)}

\newcommand{\lgradfields}{\overline{\gradfields}^\mu}

\newcommand{\ltgradfields}{\overline{\gradfields}^{\sigma_t}}

\newcommand{\ltfields}{L^2(\sigma_t)}

\newcommand{\gradfields}{\nabla\,C^\infty_c}

\newcommand{\divmu}{div_\mu}

\newcommand{\divsigma}{div_\sigma}

\newcommand{\kerdiv}{\mbox{Ker}(\divmu)}

\newcommand{\kerdivt}{\mbox{Ker}(div_{\sigma_t})}

\newcommand{ \kerdivvarphi}{ \mbox{Ker} (div_{\varphi_\# \mu}) }

\newcommand{\cm}{C^\infty(M)}

\newcommand{\ccrD}{C^\infty_c(\R^D)}

\newcommand{\diffm}{\mbox{Diff}(M)}
\newcommand{\cdiffm}{\mbox{Diff}_c(M)}
\newcommand{\diffmum}{\mbox{Diff}_\mu(M)}
\newcommand{\cdiffrD}{\mbox{Diff}_c(\R^D)}
\newcommand{\cdiffmurD}{\mbox{Diff}_{c,\mu}(\R^D)}
\newcommand{\cdiffrtwod}{\mbox{Diff}_c(\R^{2d})}

\newcommand{\sympm}{\mbox{Symp}(M)}
\newcommand{\hamm}{\mbox{Ham}(M)}
\newcommand{\symprtwod}{\mbox{Symp}_c(\Rtwod)}
\newcommand{\hamrtwod}{\mbox{Ham}_c(\Rtwod)}

\newcommand{\calLone}{\mathcal{L}^1}
\newcommand{\calLD}{\mathcal{L}^D}

\newcommand{\calO}{\mathcal{O}}
\newcommand{\calOmu}{\mathcal{O}_\mu}

\newcommand{\tTM}{\mathcal{T} \M}
\newcommand{\ctTM}{\mathcal{T}^*\M}

\newcommand{\TmuM}{T_\mu \M}
\newcommand{\intR}{\int_{\R}}
\newcommand{\intD}{\int_{\RD}}
\newcommand{\intDD}{\int_{\RDD}}

\newcommand{\id}{{ \it Id}}

\newcommand{\rb}{{\bf r}}

\newcommand{\x}{{\bf x}}
\newcommand{\y}{{\bf y}}

\def\endproof{\ \hfill QED. \bigskip}

\def\proof#1{\noindent {\bf Proof{#1}:}}

\title[Differential forms on Wasserstein space]{Differential forms on Wasserstein space and infinite-dimensional Hamiltonian systems}
\author[W. Gangbo]{Wilfrid~Gangbo}
\address{Georgia Institute of Technology, Atlanta GA, USA} \email{gangbo@math.gatech.edu}
\author[H. K. Kim]{Hwa~Kil~Kim}
\address{Georgia Institute of Tecnology, Atlanta GA, USA} \email{hwakil@math.gatech.edu}
\author[T.~Pacini]{Tommaso~Pacini}
\address{Mathematical Institute, Oxford, UK} \email{pacini@maths.ox.ac.uk}

\subjclass[2000]{Primary 37Kxx, 49-XX; Secondary 35Qxx, 53Dxx}

\date{\today}

\begin{document}
\begin{abstract}
Let $\M$ denote the space of probability measures on $\RD$ endowed with the Wasserstein metric. A differential calculus for a certain class of absolutely continuous curves in $\M$ was introduced in \cite{ags:book}. In this paper we develop a calculus for the corresponding class of differential forms on $\M$. In particular we prove an analogue of Green's theorem for 1-forms and show that the corresponding first cohomology group, in the sense of de Rham, vanishes. For $D=2d$ we then define a symplectic distribution on $\M$ in terms of this calculus, thus obtaining a rigorous framework for the notion of Hamiltonian systems as introduced in \cite{ambrosiogangbo:hamiltonian}. Throughout the paper we emphasize the geometric viewpoint and the role played by certain diffeomorphism groups of $\RD$.
\end{abstract}
\maketitle


\section{Introduction}\label{s:intro}

Historically speaking, the main goal of Symplectic Geometry has been to provide the mathematical formalism and the tools to define and study the most fundamental class of equations within classical Mechanics, \textit{Hamiltonian ODEs}. Lie groups and group actions provide a key ingredient, in particular to describe the symmetries of the equations and to find the corresponding preserved quantities. 

As the range of physical examples of interest expanded to encompass continuous media, fields, \textit{etc.}, there arose the question of reaching an analogous theory for PDEs. It has long been understood that many PDEs should admit a reformulation as \textit{infinite-dimensional Hamiltonian systems}. A deep early example of this is the work of Born-Infeld \cite{born}, \cite{borninfeld} and Pauli \cite{pauli}, who started from a Hamiltonian formulation of Maxwell's equations to develop a quantum field theory in which the commutator of operators is analogous to the Poisson brackets used in the classical theory. Further examples include the wave and Klein-Gordon equations (cf. \textit{e.g.} \cite{chernoffmarsden}, \cite{marsdenratiu:book}), the relativistic and non-relativistic Maxwell-Vlasov equations \cite{bbetal}, \cite{chendraetal}, \cite{marsdenweinstein:maxwellvlasov}, and the Euler equations for incompressible fluids \cite{arnold}. 

In each case it is necessary to define an appropriate phase space, build a symplectic or Poisson structure on it, find an appropriate energy functional, then show that the PDE coincides with the corresponding Hamiltonian flow. For various reasons, however, the results are often more formal than rigorous. In particular, existence and uniqueness theorems for PDEs require a good notion of weak solutions which need to be incorporated into the configuration and phase spaces; the geometric structure of these spaces needs to be carefully worked out; the functionals need the appropriate degree of regularity, \textit{etc}. The necessary techniques can become quite complicated and \textit{ad hoc}.

The purpose of this paper is to provide the basis for a new framework for defining and studying Hamiltonian PDEs. The configuration space we rely on is the \textit{Wasserstein space} $\M$ of non-negative Borel measures on $\RD$ with total mass $1$ and finite second moment. Over the past decade it has become clear that $\M$ provides a very useful space of weak solutions for those PDEs in which total mass is preserved. One of its main virtues is that it provides a unified theory for studying these equations. In particular, the foundation of the theory of Wasserstein spaces comes from Optimal Transport and Calculus of Variations, and these provide a toolbox which can be expected to be uniformly useful throughout the theory. Working in $\M$ also allows for extremely singular initial data, providing a bridge between PDEs and ODEs when the initial data is a Dirac measure. 

The main geometric structure on $\M$ is that of a metric space. The geometric and analytic features of this structure have been intensively studied, cf. \textit{e.g.} \cite{ags:book}, \cite{carlengangbo}, \cite{cmv}, \cite{mccann}, \cite{otto}. In particular the work \cite{ags:book} has developed a theory of \textit{gradient flows} on metric spaces. In this work the technical basis for the notion of weak solutions to a flow on $\M$ is provided by the theory of \textit{2-absolutely continuous curves}. In particular, \cite{ags:book} develops a differential calculus for this class of curves including a notion of ``tangent space" for each $\mu\in\M$. Applied to $\M$, this allows for a rigorous reformulation of many standard PDEs as gradient flows on $\M$. Overall, this viewpoint has led to important new insights and results, cf. \textit{e.g.} \cite{agueh}, \cite{ags:book}, \cite{cmv}, \cite{jko}, \cite{otto}. Topics such as geodesics, curvature and connections on $\M$ have also received much attention, cf. \cite{ambrosiogigli}, \cite{lott:somecalcs}, \cite{lottvillani}, \cite{sturm:I}, \cite{sturm:II}. 

In the case $D=2d$, recent work \cite{ambrosiogangbo:hamiltonian} indicates that other classes of PDEs can be viewed as \textit{Hamiltonian flows} on $\M$. Developing this idea requires however a rigorous symplectic formalism for $\M$, adapted to the viewpoint of \cite{ags:book}. Our paper achieves two main goals. The first is to develop a general theory of differential forms on $\M$. We present this in Sections \ref{s:calculusonM} and \ref{s:pseudocalculus}. This calculus should be thought of as dual to the calculus of absolutely continuous curves. Our main result here, Theorem \ref{th:greencurveannulus}, is an analogue of Green's theorem for 1-forms and leads to a proof that, in a specific sense, every closed 1-form on $\M$ is exact. The second goal is to show that there exists a natural symplectic and Hamiltonian formalism for $\M$ which is compatible with this calculus of curves and forms. The appropriate notions are defined and studied in Sections \ref{s:symplectic} and \ref{s:foliationispoisson}.

Given any mathematical construction, it is a fair question if it can be considered ``the most natural" of its kind. It is well known for example that cotangent bundles admit a ``canonical" symplectic structure. It is an important fact, discussed in Section \ref{s:foliationispoisson}, that on a non-technical level our symplectic formalism turns out to be formally equivalent to the Poisson structure considered in \cite{marsdenweinstein:maxwellvlasov}, cf. also \cite{khesinlee} and \cite{lott:somecalcs}. From the geometric point of view it is clear that the structure in \cite{marsdenweinstein:maxwellvlasov} is indeed an extremely natural choice. The choice of $\M$ as a configuration space is also both natural and classical. The difference between our paper and the previous literature appears precisely on the technical level, starting with the choice of geometric structure on $\M$. Specifically, whereas previous work tends to rely on various adaptations of differential geometric techniques, we choose the methods of Optimal Transport. The technical effort involved is justified by the final result: while previous studies are generally forced to restrict to smooth measures and functionals, our methods allow us to present a uniform theory which includes all singular measures and assumes very little regularity on the functionals. Sections \ref{ss:regularforms} through \ref{ss:greenstheorem} are an example of the technicalities this entails. Section \ref{subse:smoothGreen} provides instead an example of the simplifications which occur when one assumes a higher degree of regularity. 

By analogy with the case of gradient flows we expect that our framework and results will provide new impulse and direction to the development of the theory of Hamiltonian PDEs. In particular, previous work and other work in progress inspired by these results lead to existence results for singular initial data \cite{ambrosiogangbo:hamiltonian}, existence results for Hamiltonians satisfying weak regularity conditions \cite{hwakil:thesis}, and to the development of a \textit{weak KAM} theory for the nonlinear Vlasov equation \cite{gangbotudorascu}. 
It is to be expected that in the process of these developments our regularity assumptions will be even further relaxed so as to broaden the range of applications. We likewise expect that the geometric ideas underlying Symplectic Geometry and Geometric Mechanics will continue to play an important role in the development of the Wasserstein theory of Hamiltonian systems on $\M$. For example, in a very rough sense the relationship between our methods and those implicit in \cite{marsdenweinstein:maxwellvlasov} can be thought of as analogous to the relationship between \cite{ebinmarsden} and \cite{arnold}. A connection between the choice of using Lie groups (as in \cite{ebinmarsden} and \cite{arnold}) or the space of measures as configuration spaces is provided by the process of \textit{symplectic reduction}, cf. \cite{5authors}, \cite{marsdenratiu:book}. 

It is an interesting question to what extent our results can be generalized to spaces of probability measures on other manifolds $M$.  Regarding this issue, the situation is as follows. Many of the analytic foundations of our paper are provided by the work \cite{ags:book}, which is based on the choice $M:=\R^D$. In theory many results of \cite{ags:book} should be extendible to general Riemannian manifolds, but at present such an extension does not exist. Assuming that this extension will be obtained, we have written our paper in such a way as to make it clear how one might then try to extend our own results. This partly explains our emphasis on the geometric ideas and intuition underpinning our analytic definitions and results: exactly the same ideas would continue to hold for general manifolds $M$. Section \ref{ss:section5discussion} discusses how our results on cohomology depend on the choice $M:=\R^D$. The situation regarding the symplectic structure is similar: one should expect most results to continue to hold for general symplectic manifolds $M$. 

The above considerations make it worthwhile to stress the geometric viewpoint throughout this paper, with particular attention to the role played by certain group actions. It is important to emphasize, however, that we never try to use any form of infinite-dimensional geometry to prove our results. The reason behind this is that the various existing rigorous formulations of infinite-dimensional manifolds and Lie groups do not seem to be easily adaptable to our needs, cf. Section \ref{ss:section3discussion} for details. The typical approach adopted throughout our paper is thus as follows: (i) use geometric intuition to guide us towards specific choices of rigorous definitions, within the framework of \cite{ags:book}; (ii) prove theorems using the methods of \cite{ags:book} and Monge-Kantorovich theory; (iii) provide informal discussions of the geometric consequences of our results.

In recent years Wasserstein spaces have also been very useful in the field of Geometric Inequalities, cf. \textit{e.g.} \cite{agk}, \cite{cnv}, \cite{cgh}, \cite{maggivillani}. Most recently, the theory of Wasserstein spaces has started producing results in Metric and Riemannian Geometry, cf. \textit{e.g.} \cite{lottvillani}, \cite{mccanntopping}, \cite{sturm:I}, \cite{sturm:II}. Thus there exist at least three distinct communities which may be interested in these spaces: people working in Analysis/PDEs/Calculus of Variations, people in Geometrical Mechanics, people in Geometry. Concerning the exposition of our results, we have tried to take this into account in various ways: (i) by incorporating into the presentation an abundance of background material; (ii) by emphasizing the general geometric setting behind many of our constructions; (iii) by sometimes avoiding maximum generality in the results themselves. As much as possible we have also tried to keep the background material and the purely formal arguments separate from the main body of the article via a careful subdivision into sections and an appendix. 

We now briefly summarize the contents of each section. Section \ref{s:manifold} contains a brief introduction to the topological and differentiable structure (in the weak sense of \cite{ags:book}) of $\M$. Likewise, Appendix \ref{s:geometryreview} reviews various notions from Differential Geometry including Lie derivatives, differential forms, Lie groups and group actions. The material in both is completely standard, but may still be useful to some readers. Section \ref{s:Mrevisited} provides a bridge between these two parts by revisiting the differentiable structure of $\M$ in terms of group actions. Although this point of view is maybe implicit in \cite{ags:book}, it seems worthwhile to emphasize it. On a purely formal level, it leads to the conclusion that $\M$ should roughly be thought of as a \textit{stratified} rather than a smooth manifold, see Section \ref{ss:tangentspacesbis}. It also relates the sets $\RD\subset\M\subset\N$. The first inclusion, based on Dirac measures, shows that the theory on $\M$ specializes by restriction to the standard theory on $\RD$: this should be thought of as a fundamental test in this field, to be satisfied by any new theory on $\M$. The second inclusion provides background for relating the constructions of Section \ref{ss:symplecticfoliation} to the work \cite{marsdenweinstein:maxwellvlasov}. Overall, Section \ref{s:Mrevisited} is perhaps more intuitive than rigorous; however it does seem to offer a useful point of view on $\M$, providing intuition for the developments in the following sections. Section \ref{s:calculusonM} defines the basic objects of study for a calculus on $\M$, namely differential forms, push-forward operations and an exterior differential operator. It also introduces the more general notion of \textit{pseudo forms}. Pseudo forms are closely related to the group action: this is discussed in Section \ref{ss:section4discussion}. Pseudo forms reappear in Section \ref{s:pseudocalculus} as the main object of study, mainly because it seems both more natural and easier to control  their regularity. The main result of this section is an analogue of Green's theorem for certain \textit{annuli} in $\M$, Theorem \ref{th:greencurveannulus}. Stating and proving this result requires a good understanding of the differentiability and integrability properties of pseudo 1-forms. We achieve this in Sections \ref{ss:regularforms} and \ref{ss:morediff}. Our main application of Theorem \ref{th:greencurveannulus} is Corollary \ref{co:closedexact}, which shows that the 1-form defined by any closed regular pseudo 1-form on $\M$ is exact. Section \ref{ss:section5discussion} discusses the cohomological consequences of this result. In Section \ref{s:symplectic} we move on towards Symplectic Geometry, specializing to the case $D=2d$. The main material is in Section \ref{ss:symplecticfoliation}: for each $\mu\in\M$ we introduce a particular subspace of the tangent space $T_\mu\M$ and show that it carries a natural symplectic structure. We also study the geometric properties of this symplectic distribution and define the notion of Hamiltonian systems on $\M$, thus providing a firm basis to the notion already introduced in \cite{ambrosiogangbo:hamiltonian}. Formally speaking, this distribution of subspaces is integrable and the above defines a \textit{Poisson structure} on $\M$. The existence of a Poisson structure on $\N$ had already been noticed in \cite{marsdenweinstein:maxwellvlasov}: their construction is a formal infinite-dimensional analogue of Lie's construction of a canonical Poisson structure on the dual of any finite-dimensional Lie algebra. We review this construction in Section \ref{s:foliationispoisson} and show that the corresponding 2-form restricts to ours on $\M$. In this sense our construction is formally equivalent to the Kirillov-Kostant-Souriau construction of a symplectic structure on the coadjoint orbits of the dual Lie algebra.


\section{The topology on $\M$ and a differential calculus of curves} \label{s:manifold}
Let $\M$ denote the space of \textit{Borel probability measures on $\RD$ with bounded second moment}, \textit{i.e.}
$$\M:=\{\mbox{Borel measures on }\RD: \mu\geq 0, \int_{\RD}d\mu=1, \int_{\RD}|x|^2\,d\mu<\infty\}.$$
The goal of this section is to show that $\M$ has a natural metric structure and to introduce a differential calculus due to \cite{ags:book} for a certain class of curves in $\M$. We refer to \cite{ags:book} and \cite{villani:book} for further details.


\subsection{The space of distributions}\label{ss:distributions}
Let $\ccinf$ denote the space of compactly-supported smooth functions on $\RD$. Recall that it admits the structure of a complete locally convex Hausdorff topological vector space, cf. \textit{e.g.} \cite{rudin} Section 6.2. Let $(\ccinf)^*$ denote the topological dual of $\ccinf$, \textit{i.e.} the vector space of continuous linear maps $\ccinf\rightarrow\R$. We endow $(\ccinf)^*$ with the \textit{weak-* topology}, defined as the coarsest topology such that, for each $f\in \ccinf$, the induced evaluation maps
$$(\ccinf)^*\rightarrow \R,\ \ \phi\mapsto \langle \phi,f \rangle$$
are continuous. In terms of sequences this implies that, $\forall f\in \ccinf$, 
$$\phi_n\rightarrow \phi\Leftrightarrow \,\langle\phi_n,f\rangle\rightarrow \langle \phi,f\rangle.$$ 
Then $\N$ is a locally convex Hausdorff topological vector space, cf. \cite{rudin} Section 6.16. As such it has a natural differentiable structure.

The following fact may provide a useful context for the material of Section \ref{ss:topology}.
We denote by $\mathcal{P}$ the set of all Borel probability measures on $\RD$. A function $f$ on $\RD$ is said to be of \textit{p-growth} (for some $p>0$) if there exist constants $A,B\geq 0$ such that $|f(x)|\leq A+B|x|^p$. Let $C_b(\RD)$ denote the set of continuous functions with 0-growth, \textit{i.e.} the space of bounded continuous functions. As above we endow $(C_b(\RD))^*$ with its natural weak-* topology, defined using test functions in $C_b(\RD)$: this is also known as the \textit{narrow} topology. Since $\mathcal{P}$ is contained in both $(C_b(\RD))^*$ and $(\ccinf)^*$, it inherits two natural topologies. It is well known, cf. \cite{ags:book} Remarks 5.1.1 and 5.1.6, that the corresponding two notions of convergence of sequences coincide, but that the stronger topology induced from $(C_b(\RD))^*$ is more interesting in that it is metrizable.


\subsection{The topology on $\M$} \label{ss:topology}

Let $C_2(\RD)$ denote the set of continuous functions with 2-growth, as in Section \ref{ss:distributions}. We endow $(C_2(\RD))^*$ with its natural weak-* topology, defined using test functions in $C_2(\RD)$. As in Section \ref{ss:distributions}, $\M$ is contained in both $(C_2(\RD))^*$ and $(\ccinf)^*$. We will endow $\M$ with the topology induced from $(C_2(\RD))^*$. Notice that $\M$ is a convex affine subset of $(C_2(\RD))^*$. In particular it is contractible, so for $k\geq 1$ all its homology groups $H_k$ and cohomology groups $H^k$, defined topologically, vanish. As in Section \ref{ss:distributions}, it turns out that this topology is metrizable. A compatible metric can be defined as follows.

\begin{definition}\label{def:wasserstein}
Let $\mu,\,\nu\in\M.$ Consider
\begin{equation} \label{eq:wasserstein}
W_2(\mu,\nu):= \left(\inf_{\gamma\in\Gamma(\mu,\nu)}\int_{\RD\times\RD} |x-y|^2 d\gamma(x,y)\right)^{1/2}.
\end{equation}
Here, $\Gamma(\mu, \nu)$ denotes the set of Borel measures $\gamma$ on $\RD \times \RD$ which have $\mu$ and $\nu$ as \textit{marginals}, \textit{i.e.} satisfying $\pi^1_{\#}(\gamma)=\mu$ and $\pi^2_{\#}(\gamma)=\nu$ where $\pi^1$ and $\pi^2$ denote the standard projections $\RD\times\RD\rightarrow \RD$. 

Equation \ref{eq:wasserstein} defines a distance on $\M$. It is known that the infimum in the right hand side of Equation \ref{eq:wasserstein} is always achieved. We will denote by $\Gamma_o(\mu,\nu)$ the set of $\gamma$ which minimize this expression.
\end{definition}

It can be shown that $(\M, W_2)$ is a separable complete metric space, cf. \textit{e.g.} \cite{ags:book} Proposition 7.1.5. It is an important result from Monge-Kantorovich theory that
\begin{equation} \label{eq:dualwasserstein}
W^2_2(\mu, \nu)= \sup_{u, v \in C(\RD)} \Bigl\{ \intD u\,d\mu + \intD v\,d\nu : \; u(x) +v(y) \leq |x-y|^2  \ \ \forall x,y \in \RD \Bigr\}.
\end{equation}
Recall that $\mu$ is \textit{absolutely continuous} with respect to Lebesgue measure $\calLD$, written $\mu << \calLD$, if it is of the form $\mu=\rho\,\calLD$ for some function $\rho\in L^1(\RD)$. In this case for any $\nu\in\M$ there exists a unique map $T: \RD \rightarrow \RD$ such that $T_\# \mu=\nu$ and
\begin{equation} \label{eq:wasserstein2}
W^2_2(\mu,\nu)= \int_{\RD} |x-T(x)|^2 d\mu(x),
\end{equation}
cf. \textit{e.g.} \cite{ags:book} or \cite{gangbomccann:optimal}. One refers to $T$ as the \textit{optimal map} that pushes $\mu$ forward to $\nu.$

\begin{example} \label{e:discrete1} Given $x\in \RD$, let $\delta_x$ denote the corresponding \textit{Dirac measure} on $\RD$. Consider the set of such measures: this is a closed subset of $\M$ isometric  to $\RD$. More generally, let $a_i$ $(i=1,\dots,n)$ be a fixed collection of distinct positive numbers such that $\sum a_i=1$. Then the set of measures of the form $\sum a_i\delta_{x_i}$ constitutes a closed subset of $\M$, homeomorphic to $\R^{nD}$.

If $a_i\equiv 1/n$ then the set of measures of the form $\mu=\sum (1/n)\,\delta_{x_i}$ can be identified with $\R^{nD}$ quotiented by the set of permutations of $n$ letters. This space is not a manifold in the usual sense; in the simplest case $D=1$ and $n=2$, it is homeomorphic to a closed half plane, which is a manifold with boundary.
\end{example}

\begin{example} The set of all absolutely continuous measures is dense in $\M$. The set of all discrete measures, as in Example \ref{e:discrete1}, is also dense in $\M.$ Since these two sets are disjoint, neither is open nor closed in $\M$.
\end{example}


\subsection{Tangent spaces and the divergence operator} \label{ss:tangentspaces}

Let $\cfields$ denote the space of compactly-supported smooth vector fields on $\RD$. Set $\nabla C_c^\infty:=\{\nabla f: f\in C_c^\infty\}\subset\cfields$. For $\mu \in \M$ let $\lfields$ denote the set of Borel maps $X: \RD \rightarrow \RD$ such that $||X||^2_\mu:= \intD |X|^2 d\mu$ is finite. Recall that $\lfields$ is a Hilbert space with the inner product
\begin{equation}\label{eq:innerproduct}
\hat{G}_\mu(X,Y):=\intD \langle X , Y\rangle \,d\mu.
\end{equation}

\begin{remark} If $\mu=\rho \calLD$ for some $\rho:\Rd\rightarrow (0,\infty)$ such that $\int\rho dx=1$ then the natural map $\cfields\rightarrow\lfields$ is injective. But in general it is not: for example if $\mu$ is the Dirac mass at $x$ then two vector fields $X$, $Y$ will be identified as soon as $X(x)=Y(x)$.  However, the image of this map is always dense in $\lfields$.

\end{remark}

In \cite{ags:book} Section 8.4, a ``tangent space" is defined for each $\mu\in\M$ as follows.

\begin{definition} \label{def:tangentspaces}
Given $\mu \in \M$, let $T_\mu \M$ denote the closure of $\nabla C_c^\infty$ in $\lfields$. We call it the \textit{tangent space} of $\M$ at $\mu$. The \textit{tangent bundle} $T\M$ is defined as the disjoint union of all $T_\mu \M$.
\end{definition}

\begin{definition}\label{def:divergence}
Given $\mu \in \M$ we define the \textit{divergence operator}
\begin{equation*}\label{eq:divergence}
div_\mu:\cfields\rightarrow (\ccinf)^*,\ \  \langle div_\mu(X),f \rangle:=-\int_{\RD} df(X)\,d\mu.
\end{equation*}

Notice that the divergence operator is linear and that $\langle div_\mu(X),f\rangle \leq ||\nabla f||_\mu ||X||_\mu$. This proves that the operator $div_\mu$ extends to $\lfields$ by continuity; we will continue to use the same notation for the extended operator, so that $\kerdiv$ is now a closed subspace of $\lfields$.
\end{definition}
It follows from \cite{ags:book} Lemma 8.4.2 that, given any $\mu \in \mathcal{M}$, there is an orthogonal decomposition
\begin{equation}\label{eq10.11.1}
\lfields=\lgradfields\oplus \kerdiv.
\end{equation}

We will denote by $\pi_\mu:\lfields\rightarrow\lgradfields$ the corresponding projection. Notice that each tangent space has a natural Hilbert space structure $G_\mu$, obtained by restriction of $\hat{G}_\mu$ to $\lgradfields$.

\begin{remark}\label{rem:altertangentspace}
Decomposition \ref{eq10.11.1} shows that $T_\mu \M$ can also be identified with the quotient space $\lfields/\kerdiv$: the map $\pi_\mu$ provides a Hilbert space isomorphism between these two spaces.
\end{remark}

\begin{example}\label{ex:discreteextension} 
Suppose that $x_1, \cdots, x_n$ are points in $\RD$ and $\mu= \sum_{i=1}^n (1/n)\,\delta_{x_i}$. Fix $\xi \in L^2(\mu)$. Set $4r:= \min_{x_i \not = x_j} |x_i -x_j|$ and define
\begin{equation}\label{d:varphi}
\varphi(x)=
\left\{\begin{array}{rl}
\langle x, \xi(x_i)\rangle & \hbox{if} \quad x \in B_{2 r}(x_i) \quad i=1, \cdots, n\\
0 & \hbox{if} \quad x \not \in \cup_{i=1}^n B_{2 r}(x_i).
\end{array}\right.
\end{equation}
Let $\eta \in C_c^\infty$ be a symmetric function such that $\intD \eta dx=1$, $\eta \geq 0$ and $\eta$ is supported in the closure of $B_r(0).$ Then $\bar \varphi:=\eta \ast \varphi \in C_c^\infty$ and $\nabla \bar \varphi$ coincides with $\xi$ on $\cup_{i=1}^n B_{ r}(x_i).$ Consequently,  $\lfields=T_\mu\mathcal{M}$ and $\kerdiv =\{0 \}.$ In particular if the points $x_i$ are distinct then $\lfields$ can be identified with $\RnD$. If on the other hand all the points coincide, \textit{i.e.} $x_i\equiv x$, then $\mu=\delta_x$ and $\lfields\simeq\RD$.

Consider for example the simplest case $D=1$, $n=2$. As seen in Example \ref{e:discrete1} the corresponding space of Dirac measures is homeomorphic to a closed half plane. We now see that at any interior point, corresponding to $x_1\neq x_2$, the tangent space is $\R^2$. At any boundary point, corresponding to $x_1=x_2$, the tangent space is $\R$. One should compare this with the usual differential-geometric definition of tangent planes on a manifold with boundary, cf. \textit{e.g.} \cite{guilleminpollack:book}: in that case the tangent plane at a boundary point would be $\R^2$. We will come back to this in Section \ref{ss:tangentspacesbis}.  
\end{example}

\begin{remark} 
Decomposition \ref{eq10.11.1} extends  the standard orthogonal Hodge decomposition of a smooth $L^2$ vector field $X$ on $\RD$:
$$X=\nabla u+X',$$
where $u$ is defined as the unique smooth solution in $W^{1,2}$ of $\Delta u=div(X)$ and $X':=X-\nabla u$.

In particular, Decomposition \ref{eq10.11.1} shows that $\lgradfields\cap \kerdiv=\{0\}$. The analogous statement with respect to the measure $\mathcal{L}^D$ is that the only harmonic function on $\RD$ in $W^{1,2}$ is the function $u\equiv 0$.
\end{remark}


\subsection{Analytic justification for the tangent spaces}\label{ss:antangentspaces}
Following \cite{ags:book} we now provide an analytic justification for the above definition of tangent spaces for $\M$. A more geometric justification, using group actions, will be given in Section \ref{ss:tangentspacesbis}.

Suppose we are given a curve $\sigma: (a,b)\rightarrow \M$ and a Borel vector field $X:(a,b) \times \RD \rightarrow \RD$ such that $X_t \in   \ltfields.$ Here, we have written $\sigma_t$ in place of $\sigma(t)$ and $X_t$ in place of $X(t).$   We will write
\begin{equation}\label{eq:cont:eq}
\frac{\partial\, \sigma}{\partial t}+\divsigma(X)=0
\end{equation}
if the following condition holds: for all $\phi\in C^\infty_c((a,b)\times\RD)$,
\begin{equation}\label{eq:int:cont:eq}
\int_a^b\int_{\RD}\Bigl( \frac{\partial\,\phi}{\partial\,t}+d\phi(X_t) \Bigr) \,d\sigma_t\,dt=0,
\end{equation}
\textit{i.e.} if Equation \ref{eq:cont:eq} holds in the sense of distributions. Given $\sigma_t$, notice that if Equation \ref{eq:cont:eq} holds for $X$ then it holds for $X+W$, for any Borel map $W: (a,b) \times \RD \rightarrow \RD$ such that $W_t\in\kerdivt$.

The following definition and remark can be found in \cite{ags:book} Chapter 1.

\begin{definition}\label{deac2curves}
Let $(\calS, \dist)$ be a complete metric space. A curve $t \in (a,b) \mapsto \sigma_t \in \calS$ is \textit{2-absolutely continuous} if there exists $\beta \in L^2(a,b)$ such that
$\dist(\sigma_t, \sigma_s) \leq \int_s^t \beta(\tau)d\tau$
for all $a <s<t<b.$ We then write $\sigma \in AC_2(a,b; \calS).$ For such curves the limit $|\sigma'|(t):=\lim_{s \rightarrow t} \dist(\sigma_{t}, \sigma_s)/|t-s| $ exists for $\calL^1$-almost every $t \in (a,b)$. We call this limit the \textit{metric derivative} of $\sigma$ at $t.$ It satisfies $|\sigma'| \leq \beta$ $\calL^1$-almost everywhere.
\end{definition}

\begin{remark}\label{reac2curves}
(i) If $\sigma \in AC_2(a,b; \calS)$ then $|\sigma'| \in L^2(a,b)$ and $\dist(\sigma_s, \sigma_t) \leq \int_s^t |\sigma'|(\tau)d\tau$ for $a<s<t<b.$ We can apply H\"older's inequality to conclude that $\dist^2(\sigma_s, \sigma_t) \leq c |t-s|$
where $c= \int_a^b |\sigma'|^2(\tau)d\tau.$ 

(ii) It follows from (i) that $\{\sigma_t | \; t \in [a,b]\}$ is a compact set, so it is bounded. For instance,  given $x \in \calS$, the triangle inequality proves that $\dist(\sigma_s, x) \leq \sqrt{c |s-a|} +\dist(\sigma_a, x).$
\end{remark}

We now recall \cite{ags:book} Theorem 8.3.1. It shows that the definition of tangent space given above is flexible enough to include the velocities of any ``good" curve in $\M$.

\begin{prop}\label{prac2curves} 
If $\sigma \in AC_2(a,b; \M)$ then there exists a Borel map $v:  (a,b)\times\RD\rightarrow\RD$ such that $\frac{\partial\,\sigma}{\partial t}+\divsigma(v)=0$ and $v_t \in L^2(\sigma_t)$ for $\calLone$-almost every $t \in (a,b).$ We call $v$ a \textit{velocity} for $\sigma.$ If $w$ is another velocity for $\sigma$ then the projections $\pi_{\sigma_t}(v_t)$, $\pi_{\sigma_t}(w_t)$ coincide for $\calLone$-almost every $t \in (a,b).$ One can choose $v$ such that $ v_t  \in  \ltgradfields$ and $||v_t||_{\sigma_t}=|\sigma'|(t)$ for $\calL^1$-almost every $t \in (a,b)$. In that case, for $\calL^1$-almost every $t \in (a,b)$, $v_t$ is uniquely determined. We denote this velocity $\dot{\sigma}$ and refer to it as the \textit{velocity of minimal norm}, since if $w_t$ is any other velocity associated to $\sigma$ then $||\dot \sigma_t||_{\sigma_t} \leq ||w_t||_{\sigma_t}$ for $\calL^1$-almost every $t \in (a,b).$ 
\end{prop}

The following remark can be found in \cite{ags:book} Lemma 1.1.4 in a more general context.

\begin{remark} [Lipschitz reparametrization]\label{reparametrization1}
Let $\sigma \in AC_2(a,b; \M)$ and $v$ be a velocity associated to $\sigma.$ Fix $\alpha>0$ and define $S(t)=\int_a^t \bigl(\alpha +||v_\tau||_{\sigma_\tau} \bigr)d\tau.$ Then $S:[a,b]\rightarrow [0,L]$ is absolutely continuous and increasing, with $L=S(b).$  The inverse of $S$ is a function whose Lipschitz constant is less than or equal to $1/\alpha.$ Define
$$\bar \sigma_s:=\sigma_{S^{-1}(s)},\qquad \bar v_s:=\dot S^{-1}(s) v_{S^{-1}(s)}.$$
One can check that $\bar \sigma \in AC_2(0,L; \M)$ and that $\bar v$ is a velocity associated to $\bar \sigma.$ Fix $t \in (a,b)$ and set $s:=S(t).$ Then  $v_t= \dot S(t) \bar v_{S(t)}$ and $||\bar v_{s}||_{\sigma_s}={||v_{t}||_{\sigma_t} \over \alpha + ||v_{t}||_{\sigma_t}}<1.$
\end{remark}


\section{The calculus of curves, revisited}\label{s:Mrevisited}
The goal of this section is to revisit the material of Section \ref{s:manifold} from a more geometric viewpoint. Many of the results presented here are purely formal, but they may provide some insight into the structure of $\M$. They also provide useful intuition into the more rigorous results contained in the sections which follow. We refer to Appendix \ref{s:geometryreview} for notation and terminology.

\subsection{Embedding the geometry of $\RD$ into $\M$}\label{ss:embedding}

We have already seen in Example \ref{e:discrete1} that Dirac measures provide a continuous embedding of $\RD$ into $\M$. Many aspects of the standard geometry of $\RD$ can be recovered inside $\M$, and various techniques which we will be using for $\M$ can be seen as an extension of standard techniques used for $\RD$. 

One example of this is provided by Example \ref{ex:discreteextension}, which shows that the standard notion of tangent space on $\RD$ coincides with the notion of tangent spaces on $\M$ introduced by \cite{ags:book}.

Another simple example concerns calculus on $\R^D$, as follows. Consider the space of \textit{volume forms} on $\RD$, \textit{i.e.} the smooth never-vanishing D-forms. Under appropriate normalization and decay conditions, these define a subset of $\M$. Given a vector field $X\in\cfields$ and a volume form $\alpha$, there is a standard geometric definition of $div_\alpha(X)$ in terms of Lie derivatives: namely, $\mathcal{L}_X\alpha$ is also a D-form so we can define $div_\alpha(X)$ to be the unique smooth function on $\RD$ such that
\begin{equation}
div_\alpha(X)\alpha=\mathcal{L}_X\alpha.
\end{equation}
In particular, it is clear from this definition and Lemma \ref{l:preservedform} that $X\in\mbox{Ker}(div_\alpha)$ iff the corresponding flow preserves the volume form.

Cartan's formula \ref{eq:cartan} together with Green's theorem for $\R^D$ shows that $div_\alpha$ is the negative formal adjoint of $d$ with respect to $\alpha$, \textit{i.e.}
$$\int_{\R^D} f\,div_\alpha(X)\alpha=-\int_{\R^D} df(X)\,\alpha,\ \ \forall f\in C_c^\infty.$$
In particular, $div_\alpha(X)\alpha$ satisfies Equation \ref{eq:divergence}. In this sense Definition \ref{def:divergence} extends the standard geometric definition of divergence to the whole of $\M$.


\subsection{The intrinsic geometry of $\M$}\label{ss:tangentspacesbis}
It is appealing to think that, in some weak sense, the results of Section \ref{ss:antangentspaces} can be viewed as a way of using the Wasserstein distance to describe an ``intrinsic" differentiable structure on $\M$. This structure can be alternatively viewed as follows. 

Let $\phi:\R^D\rightarrow\R^D$ be a Borel map and $\mu\in\M$. Recall that the \textit{push-forward} measure $\phi_\#\mu\in\M$ is defined by setting $\phi_\#\mu(A):=\mu(\phi^{-1}(A))$, for any open subset $A\subseteq\R^D$. Let $\cdiffrD$ denote the $Id$-component of the Lie group of diffeomorphisms of $\R^D$ with compact support, cf. Section \ref{ss:diffM}. Choose any $X\in\cfields$ and let $\phi_t$ denote the flow of $X$. Given any $\mu\in\M$, it is simple to verify that $\mu_t:=\phi_{t\#}\mu$ is a path in $\M$ with velocity $X$ in the sense of Proposition \ref{prac2curves}. Notice that in this case the velocity is defined for all $t$, rather than only for almost every $t$. In particular the minimal velocity of $\mu_t$ at $t=0$ is $\pi_{\mu}(X)\in T_{\mu}\M$. From the point of view of Section \ref{ss:groupactions}, this construction can be rephrased as follows. The map 
\begin{equation}\label{eq:cdiffonM}
\cdiffrD\times\M\rightarrow\M, \ \ (\phi,\mu)\mapsto\phi_{\#}\mu
\end{equation}
is continuous and defines a \textit{left action} of $\cdiffrD$ on $\M$. The map 
$$\M\rightarrow T\M,\ \ \mu\rightarrow \pi_{\mu}(X)\in T_{\mu}\M$$
then defines the \textit{fundamental vector field} associated to $X$ in the sense of Section \ref{ss:groupactions}. 

According to Section \ref{ss:groupactions}, the orbit and stabilizer of any fixed $\mu\in\M$ are:
$$\mathcal{O}_\mu:=\{\nu\in\M:\nu=\phi_\#\mu, \mbox{ for some }\phi\in \cdiffrD\},$$
$$\cdiffmurD:=\{\phi\in\cdiffrD:\phi_\#\mu=\mu\}.$$
Formally, $\cdiffmurD$ is a Lie subgroup of $\cdiffrD$ and $\kerdiv$ is its Lie algebra. The map
$$j:\cdiffrD/\cdiffmurD\rightarrow\mathcal{O}_\mu,\ \ [\phi]\mapsto \phi_\#\mu$$
defines a 1:1 relationship between the quotient space and the orbit of $\mu$. Lemma \ref{l:orbits} suggests that $\mathcal{O}_\mu$ is a smooth manifold inside the topological space $\M$ and that the isomorphism $\nabla j:\cfields/\kerdiv\rightarrow T_{\mu}\mathcal{O}_\mu$ coincides with the map determined by the construction of fundamental vector fields. Notice that, up to $L^2_\mu$-closure, the space $\cfields/\kerdiv$ is exactly the space introduced in Definition \ref{def:tangentspaces}. This indicates that the tangent spaces of Section \ref{ss:tangentspaces} should be thought of as ``tangent" not to the whole of $\M$, but only to the leaves of the foliation induced by the action of $\cdiffrD$. In other words $\M$ should be thought of as a \textit{stratified manifold}, \textit{i.e.} as a topological space with a foliation and a differentiable structure defined only on each leaf of the foliation. This point of view is purely formal but it corresponds exactly to the situation already described for Dirac measures, cf. Example \ref{ex:discreteextension}.

Recall from Proposition \ref{prac2curves} the relationship between the class of 2-absolutely continuous curves and these tangent spaces. This result can be viewed as the expression of a strong compatibility between two natural but \textit{a priori} distinct structures on $\M$: the Wasserstein topology and the group action.

\begin{remark} The claim that the Lie algebra of $\cdiffmurD$ is $\kerdiv$ can be supported in various ways. For example, assume $\phi_t$ is a curve of diffeomorphisms in $\cdiffmurD$ and that $X_t$ satisfies Equation \ref{eq:tflow}. The following calculation is the weak analogue of Lemma \ref{l:preservedform}. It shows that $X_t\in\kerdiv$:
\begin{align*}
\int df(X_t)\,d\mu &= \int df(X_t)\,d(\phi_{t\#}\mu)=\int df_{|\phi_t}(X_{t|\phi_t})\,d\mu=\int d/dt(f\circ\phi_t)\,d\mu\\
&= d/dt\int f\circ\phi_t\,d\mu=d/dt\int f\,d(\phi_{t\#}\mu)= d/dt\int f\,d\mu\\
&=0.
\end{align*}
It is also simple to check that $\kerdiv$ is a Lie subalgebra of $\cfields$, \textit{i.e.} if $X,Y\in\kerdiv$ then $[X,Y]\in\kerdiv$. To show this, let $f\in \ccinf$. Then:
\begin{align*}
\langle div_\mu([X,Y]),f \rangle &=-\int_{\RD} df([X,Y])\,d\mu=-\int_{\Rd} dg(X)\,d\mu+\int_{\Rd} dh(Y)\,d\mu\\
&=\langle div_\mu(X),g\rangle -\langle div_\mu(Y),h\rangle \\
&=0,
\end{align*}
where $g:=df(Y)$ and $h:=df(X)$.

Finally, assume $\mu$ is a smooth volume form on a compact manifold $M$. In this situation Hamilton \cite{hamilton:nashmoser} proved that $\diffmum$ is a Fr\'echet Lie subgroup of $\diffm$ and that the Lie algebra of $\diffmum$ is the space of vector fields $X\in\xm$ satisfying the condition $\mathcal{L}_X\mu=0$. As seen in Section \ref{ss:embedding} this space coincides with $\kerdiv$.
\end{remark}


\subsection{Embedding the geometry of $\M$ into $\N$}\label{ss:embeddingbis}
We can also view $\M$ as a subspace of $\N$. It is then interesting to compare the corresponding geometries, as follows.

Consider the natural left action of $\cdiffrD$ on $\RD$ given by $\phi\cdot x:=\phi(x)$. As in Section \ref{ss:groupactions}, this induces a left action on the spaces of forms $\Lambda^k$ and in particular on the space of functions $\ccinf=\Lambda^0$, as follows:
\begin{equation*}
\cdiffrD\times\ccinf\rightarrow \ccinf,\ \ \phi\cdot f:=(\phi^{-1})^*f=f\circ \phi^{-1}.
\end{equation*}
By duality there is an induced left action on the space of distributions given by
\begin{equation*}
\cdiffrD\times\N\rightarrow \N,\ \ \langle(\phi\cdot\mu),f\rangle:=\langle\mu,(\phi^{-1}\cdot f)\rangle=\langle\mu,(f\circ\phi)\rangle.
\end{equation*}
Notice that we have introduced inverses to ensure that these are left actions, cf. Remark \ref{r:rightactions}. It is clear that this extends the action already defined in Section \ref{ss:tangentspacesbis} on the subset $\M\subset\N$. In other words, the natural immersion $i:\M\rightarrow\N$ is \textit{equivariant} with respect to the action of $\cdiffrD$, \textit{i.e.} $i(\phi_\#\mu)=\phi\cdot i(\mu)$. 

As mentioned in Section \ref{ss:distributions}, $\N$ has a natural differentiable structure. In particular it has well-defined tangent spaces $T_{\mu}\N=\N$. For each $\mu\in\M$, using the notation of Section \ref{ss:tangentspacesbis}, composition gives an immersion
\begin{equation*}
i\circ j:\cdiffrD/\cdiffmurD\rightarrow\mathcal{O}_\mu\rightarrow\N.
\end{equation*}
This induces an injection between the corresponding tangent spaces
\begin{equation*}
\nabla(i\circ j):\cfields/\kerdiv\rightarrow T_\mu\N.
\end{equation*}
Notice that, using the equivariance of $i$,
\begin{align*}
\langle\nabla(i\circ j)(X),f\rangle &= \langle \nabla i (d/dt(\phi_{t\#}\mu)_{|t=0}),f\rangle=\langle d/dt(i(\phi_{t\#}\mu))_{|t=0},f\rangle=\langle d/dt(\phi_t\cdot\mu)_{|t=0},f\rangle\\
&= d/dt\,\langle \mu,f\circ\phi_t\rangle_{|t=0}=\langle\mu,d/dt (f\circ\phi_t)_{|t=0}\rangle=\langle\mu,df(X)\rangle\\
&= -\langle div_{\mu}(X),f\rangle.
\end{align*}
In other words, the negative divergence operator can be interpreted as the natural identification between $T_{\mu}\M$ and the appropriate subspace of $\N$. 

More generally, we can compare the calculus of curves in $\M$ with the calculus of the corresponding curves in $\N$. Given any sufficiently regular curve of distributions $t \rightarrow \mu_t\in\N$, we can define tangent vectors $\tau_t:=\lim_{h\rightarrow 0} \frac{\mu_{t+h}-\mu_t}{h}\in T_{\mu_t}\N$. Assume that $\mu_t$ is strongly continuous, in the sense that the evaluation map
$$(a,b)\times \ccinf\rightarrow\R,\ \ (t,f)\mapsto\langle\mu_t,f\rangle$$
is continuous. Notice that $\mu=\mu_t$ defines a distribution on the product space $(a,b)\times\RD$: $\forall f=f_t(x)\in C^\infty_c((a,b)\times\RD)$,
$$\langle \mu,f\rangle:=\int_a^b\langle\mu_t,f_t\rangle \,dt.$$ 
One can check that $\frac{d}{dt}\langle\mu_t,f_t\rangle=\langle \tau_t,f_t\rangle+\langle \mu_t,\frac{\partial{f_t}}{\partial t}\rangle$, so
\begin{equation}\label{eq:moreconteq}
\int_a^b\langle\mu_t,\frac{\partial{f_t}}{\partial t}\rangle+\langle\tau_t,f_t\rangle\,dt=0.
\end{equation}
Equation \ref{eq:moreconteq} shows that if $\mu_t\in\M$ and $\tau_t=-div_{\mu_t}(X_t)$ then $\mu_t$ satisfies Equation \ref{eq:int:cont:eq}. In other words, the defining equation for the calculus on $\M$, Equation \ref{eq:cont:eq}, is the natural weak analogue of the statement $\lim_{h\rightarrow 0} \frac{\mu_{t+h}-\mu_t}{h}=-div_{\mu_t}(X_t)$. 

Roughly speaking, the content of Proposition \ref{prac2curves} is that if $\mu_t\in\M$ is 2-absolutely continuous then, for almost every $t$, $\tau_t$ exists and can be written as $-div_{\mu_t}(X_t)$ for some $t$-dependent vector field $X_t$ on $\RD$.
\begin{remark}
One should think of Equation \ref{eq:cont:eq}, \textit{i.e.} $d/dt(\mu_t)=-div_{\mu_t}(X_t)$, as an ODE on the submanifold $\M\subset\N$ rather than on the abstract manifold $\M$, in the sense that the right hand side is an element of $T_{\mu_t}\N$ rather than an element of $T_{\mu_t}\M$.  Using $\nabla(i\circ j)^{-1}$ we can rewrite this equation as an ODE on the abstract manifold $\M$, \textit{i.e.} $d/dt(\mu_t)=\pi_{\mu_t}(X)$.
\end{remark}


\subsection{Further comments}\label{ss:section3discussion}
There exists an extensive literature concerning how to make infinite-dimensional geometry rigorous. The first step is to provide rigorous definitions of infinite-dimensional manifolds and of infinite-dimensional Lie groups. The works \cite{ebinmarsden}, \cite{krieglmichor},  \cite{milnor:liegroups} and \cite{omori} are examples of standard references in this field. In all these cases the starting point is a notion of manifold built by gluing together charts which are open subsets of locally convex vector spaces (plus some completeness condition). These references are often useful when one wants to make a ``formal" study of PDE rigorous; in particular, when the space of solutions is, in some sense, an infinite-dimensional Lie group (as in \cite{arnold}). 

In this paper, however, we do not rely on the above frameworks. The main reason is, quite simply, the fact that the ``differentiable structure" on $\M$ introduced by \cite{ags:book} is not based on the above notion of manifold: as discussed in Section \ref{ss:tangentspacesbis}, it uses a much weaker notion and none of the results presented in this paper require anything more than this. It should also be emphasized that the relationship, discussed in Chapters \ref{s:manifold} and \ref{s:Mrevisited}, between the Wasserstein metric on $\M$, group actions and the theory of \cite{ags:book} shows that the latter is extremely natural within this context.

Another reason is that we want to keep regularity assumptions to a minimum. In particular, we want to avoid making unnecessary restrictions on the smoothness of measures (required by \cite{ebinmarsden}) and of velocity fields.

It may also be worthwhile to mention that it is not clear if the above references lead to a general theory of ``infinite-dimensional homogeneous spaces", which is in some sense the geometry underlying this paper. Specifically, Section \ref{ss:tangentspacesbis} introduces the idea that $\M$ is a ``stratified manifold" and that each stratum is the orbit of a certain group action. These orbits are of the form $G/H$, where $G$ is the infinite-dimensional Lie group of diffeomorphisms (as in the above references) and $H$ is the subgroup of diffeomorphisms which preserve a given measure. Thus our space of solutions $\M$ can be viewed as a collection of homogeneous spaces $G/H$ of varying dimension: some finite-dimensional, others infinite-dimensional. However, except in the case of smooth measures discussed in \cite{ebinmarsden}, these $H$ are not known to be ``infinite-dimensional Lie groups" and, to our knowledge, the corresponding homogeneous spaces are not known to be ``infinite-dimensional manifolds" in the sense of the above references.


\section{Tangent and cotangent bundles} \label{s:calculusonM}
We now define some further elements of calculus on $\M$. As opposed to Section \ref{s:Mrevisited}, the definitions and statements made here are completely rigorous. We will often refer back to the ideas of Section \ref{s:Mrevisited} and to the appendix, however, to explain the geometric intuition underlying this theory.

\subsection{Push-forward operations on $\M$ and $T\M$}
The following results concern the push-forward operation on $\M$.

\begin{lemma} \label{l:dualwasserstein2}
If $\phi: \RD \rightarrow \RD$ is a Lipschitz map with Lipschitz constant $Lip\,\phi$ then $\phi_\# : \M \rightarrow \M$ is also a Lipschitz map with the same Lipschitz constant.
\end{lemma}
\proof{} Let $\mu, \nu \in \M$. Note that if $u(x) +v(y) \leq |x-y|^2$ for all $x, y \in \RD$ then
$$u\circ\phi(a)+v\circ\phi(b)\leq |\phi(a)-\phi(b)|^2\leq (Lip\,\phi)^2 |a-b|^2.$$
This, together with Equation \ref{eq:dualwasserstein}, yields
\begin{equation} \label{eq:dualwasserstein2}
\intD u\,d\phi_\# \mu + \intD v\,d\phi_\# \nu= \intD u \circ \phi\,d \mu + \intD v \circ \phi\,d \nu \leq (Lip\,\phi)^2 W^2_2(\mu, \nu).
\end{equation}
We maximize the expression on the left hand side of Equation \ref{eq:dualwasserstein2}  over the set of pairs $(u,v)$ such that $u(x)+v(y)\leq |x-y|^2$ for all $x,y\in\RD.$ Then we use again Equation \ref{eq:dualwasserstein} to conclude the proof.
\endproof

\begin{lemma} \label{le:pushforwardvelocity}
For any $\mu \in \M$ and $\phi \in \cdiffrD$, the map $\phi_*:\cfields\rightarrow\cfields$ has a unique continuous extension $\phi_*: L^2(\mu) \rightarrow L^2(\phi_\# \mu)$. Furthermore $\phi_* \bigl( \kerdiv \bigr) \subseteq \kerdivvarphi$. Thus $\phi_*$ induces a continuous map $\phi_*: T_\mu\mathcal{M} \rightarrow T_{\phi_\# \mu}\mathcal{M}$.
\end{lemma}

\proof{} Let $\mu \in \M$, $\phi \in \cdiffrD$, $f \in \ccrD$ and let $X \in  \kerdiv$. If $C_\phi$ is the $L^\infty$-norm of $\nabla \phi$ we have $||\phi_* X||_{\phi_\# \mu} \leq C_\phi ||X||_\mu$. Hence $\phi_*$ admits a unique continuous linear extension. Furthermore
\begin{equation*}
\intD df(\varphi_\ast X)\,d\varphi_\# \mu
=\intD df_{|\varphi}(\varphi_\ast X_{|\varphi})\,d \mu
=\intD df_{|\varphi}(\nabla \varphi\cdot X)\,d \mu
=\intD d(f \circ \varphi)(X)\,d \mu=0.
\end{equation*}
\endproof

\begin{lemma}\label{lepullforwardcurve}

Let $\sigma \in AC_2(a,b; \M)$ and let $v$ be a velocity for $\sigma.$ Let $\varphi\in\cdiffrD$. Then $ t \rightarrow \varphi_\# (\sigma_t) \in AC_2(a,b; \M)$ and $\varphi_* v$ is a velocity for $\varphi_\# \sigma.$
\end{lemma}
\proof{} If $a<s<t<b$ then, by Lemma \ref{l:dualwasserstein2},  $W_2(\varphi_\# \sigma_t, \varphi_\# \sigma_s) \leq (Lip\,\varphi ) \, W_2( \sigma_t, \sigma_s).$ Since $\sigma \in AC_2(a,b; \M)$ one  concludes that $\varphi_\# (\sigma) \in AC_2(a,b; \M).$ If $f \in C_c^\infty((a,b) \times \RD)$ we have
\begin{align*}
\int_a^b \intD \Bigl( {\partial f_t \over \partial t} + df_t(\phi_*v_t) \Bigr) d (\varphi_\# \sigma_t) dt
&=
\int_a^b \intD \Bigl( {\partial f_t \over \partial t} \circ \varphi +( df_t (\phi_*v_t)) \circ \varphi \Bigr) d  \sigma_t dt \\
&=
\int_a^b \intD \Bigl( {\partial (f_t \circ \varphi) \over \partial t}  + d(f_t \circ \varphi) (v_t) \Bigr) d  \sigma_t dt
\\
&=0.
\end{align*}
To obtain the last equality we have used that $(t,x) \rightarrow f(t, \varphi (x))$ is in $C_c^\infty((a,b) \times \RD)$.  

\endproof


\subsection{Differential forms on $\M$} \label{ss:forms}
Recall from Definition \ref{def:tangentspaces} that the tangent bundle $T\M$ of $\M$ is defined as the union of all spaces $T_\mu\M$, for $\mu \in \M$. We now define the \textit{pseudo tangent bundle} $\tTM$ to be the union of all spaces $L^2(\mu)$. Analogously, the union of the dual spaces $T^*_\mu\M$ defines the \textit{cotangent bundle} $T^*\M$; we define the \textit{pseudo cotangent bundle} $\ctTM$ to be the union of the dual spaces $L^2(\mu)^*$.

It is clear from the definitions that we can think of $T\M$ as a subbundle of $\tTM$. Decomposition \ref{eq10.11.1} allows us also to define an injection $T^*\M\rightarrow\ctTM$ by extending any covector $T_\mu\M\rightarrow\R$ to be zero on the complement of $T_\mu\M$ in $\lfields$. In this sense we can also think of $T^*\M$ as a subbundle of $\ctTM$. The projections $\pi_\mu$ from Section \ref{ss:tangentspaces} combine to define a surjection $\pi:\tTM\rightarrow T\M$. Likewise, restriction yields a surjection $\ctTM\rightarrow T^*\M$.

\begin{remark} \label{r:downwiththemetric}
The above constructions make heavy use of the Hilbert structure on $\lfields$. Following the point of view of Remark \ref{rem:altertangentspace} and Section \ref{ss:tangentspacesbis}, \textit{i.e.} emphasizing the differential, rather than the Riemannian, structure of $\M$ one could decide to define $T_\mu\M$ as $\lfields/\kerdiv$. Then the projections $\pi_\mu:\lfields\rightarrow T_\mu\M$ would still define by duality an injection $T^*\M\rightarrow \ctTM$: this would identify $T^*\M$ with the \textit{annihilator} of $\kerdiv$ in $\lfields$. However there would be no natural injection $T\M\rightarrow \tTM$ nor any natural surjection $\ctTM\rightarrow T^*\M$.
\end{remark}

\begin{definition}\label{d:kforms} A \textit{1-form} on $\M$ is a section of the cotangent bundle $T^*\M$, \textit{i.e.} a collection of maps $\mu\mapsto \Lambda_\mu \in T^*_\mu\M$. A \textit{pseudo 1-form} is a section of the pseudo cotangent bundle $\ctTM$. Analogously, a \textit{2-form} on $\M$ is a collection of alternating multilinear maps
$$\mu\mapsto \Lambda_\mu:T_\mu\M\times T_\mu\M\rightarrow\R.$$
A \textit{pseudo 2-form} is a collection of alternating multilinear maps
$$\mu\mapsto \bar\Lambda_\mu:\lfields \times \lfields \rightarrow\R.$$
It is natural (but in practice sometimes too strong) to further assume that each $\Lambda_\mu$ (or $\bar\Lambda_\mu$) satisfies a continuity assumption such as the following: there exists $c_\mu\in\R$ such that 
\begin{equation*}
|\Lambda_\mu(X_1,X_2)|\leq c_\mu\|X_1\|_\mu\cdot\|X_2\|_\mu.
\end{equation*}
For $k=1,2$ we let $\Lambda^k\M$ (respectively, $\bar\Lambda^k\M$) denote the space of k-forms (respectively, pseudo k-forms). We define a \textit{0-form} to be a function $F:\M\rightarrow \R$.
\end{definition}
Notice that, for $k=1$, continuity implies that any 1-form is uniquely defined by its values on any dense subset of $T_\mu\M$, \textit{e.g.} on the dense subset defined by smooth gradient vector fields. The analogue holds also for pseudo forms and for $k=2$, as long as the continuity condition holds. As above, extension defines a natural injection 
\begin{equation}\label{eq:bar}
\Lambda^k\M\rightarrow\bar\Lambda^k\M, \ \ \Lambda\mapsto\bar\Lambda,
\end{equation} 
\textit{i.e.} to every k-form one can associate a canonical pseudo k-form. Conversely, restriction defines a surjection $\bar\Lambda^k\M\rightarrow\Lambda^k\M$.

Since $\TmuM$ is a Hilbert space, by the Riesz representation theorem every 1-form $\Lambda_\mu$ on $\TmuM$ can be written  $ \Lambda_\mu(Y)= \intD \langle  A_\mu , Y \rangle d\mu$ for a unique $A_\mu \in \TmuM$ and all $Y\in\TmuM$. The analogous fact is true also for pseudo 1-forms.

\begin{example}\label{e:deflinearforms}
Any $f\in\ccinf$ defines a function on $\M$, \textit{i.e.} a 0-form, as follows:
$$F(\mu):=\int_{\R^D}f d\mu.$$
We will refer to these as the \textit{linear} functions on $\M$, in that the natural extension to the space $\N$ defines a function which is linear with respect to $\mu$.

Any $\bar A\in\cfields$ defines a pseudo 1-form on $\M$ as follows:
\begin{equation}\label{eq:linear1forms}
\bar\Lambda_\mu(X):=\int_{\R^D}\langle \bar A,X\rangle d\mu.
\end{equation}
We will refer to these as the \textit{linear} pseudo 1-forms. Notice that if $\bar A=\nabla f$ for some $f\in C_c^\infty$ then $\bar\Lambda$ is actually a 1-form.

Any bounded field $B=B(x)$ on $\RD$ of $D\times D$ matrices defines a \textit{linear} pseudo 2-form via
\begin{equation}\label{eq:linear2forms}
\hat{B}_\mu(X,Y):=\int_{\R^D}\langle BX,Y\rangle d\mu.
\end{equation}
\end{example}

\begin{remark}
When $k=1$ H\"{o}lder's inequality for the product of two functions shows that any vector field $\bar A_\mu\in L^p(\mu)$, for $p\in [2,\infty]$, defines a continuous map $L^2(\mu)\rightarrow\R$ as in Equation \ref{eq:linear1forms}. When $k=2$ H\"{o}lder's inequality for the product of three functions shows that any field of matrices $B_\mu\in L^\infty(\mu)$ defines a continuous map $L^2(\mu)\times L^2(\mu)\rightarrow \R$ as in Equation \ref{eq:linear2forms}. However, for $k \geq 3$ there do not exist analogous constructions of (non-trivial) continuous k-linear maps $L^2(\mu)\times\dots\times L^2(\mu)\rightarrow \R.$ It is for this reason that we restrict our attention to the case $k\leq 2$. In any case, this is sufficient for the applications of interest to us.
\end{remark}

As in Section \ref{ss:groupactions}, the action of $\cdiffrD$ on $\M$ can be lifted to forms and pseudo forms as follows.
\begin{definition} \label{d:pullback}
For $k=1,2$, let $\bar\Lambda$ be a pseudo k-form on $\M$. Then any $\phi\in\cdiffrD$ defines a \textit{pull-back} k-multilinear map $\phi^*\bar\Lambda$ on $\M$ as follows:
$$(\phi^*\bar\Lambda)_{\mu}(X_1,\dots,X_k):=\bar\Lambda_{\phi_\#\mu}(\phi_*X_1,\dots,\phi_*X_k).$$
It follows from Lemma \ref{le:pushforwardvelocity} that the push-forward operation preserves Decomposition \ref{eq10.11.1}. This implies that the pull-back preserves the space of k-forms, \textit{i.e.} the pull-back of a k-form is a k-form.
\end{definition}

\begin{definition} \label{def:Fdifferentiable}
Let $F:\M\rightarrow\R$ be a function on $\M$. We say that
$\xi \in L^2(\mu)$ belongs to the \textit{subdifferential} $\partial_{-} F(\mu)$ if
$$ F(\nu) \geq F(\mu) + \sup_{\gamma \in \Gamma_o(\mu,\nu)}\iint_{\RD\times \RD} \langle\xi(x),y-x\rangle\,d\gamma(x,y) + o(W_2(\mu,\nu)),$$
as $\nu \rightarrow \mu$. If $-\xi \in \partial_{-} (-F)(\mu)$ we say that $\xi$ belongs to the \textit{superdifferential} $\partial^{+} F(\mu)$. 

If $\xi \in \partial_{-} F(\mu) \cap \partial^{+} F(\mu)$ then, for any $\gamma\in \Gamma_o(\mu,\nu)$, 
\begin{equation}\label{eq:differentiable0} 
F(\nu) = F(\mu) + \iint_{\RD\times \RD} \langle\xi(x),y-x\rangle\,d\gamma(x,y) + o(W_2(\mu,\nu)).
\end{equation} 
If such $\xi$ exists we say that $F$ is \textit{differentiable} at $\mu$ and we define the \textit{gradient vector} $\nabla_\mu F:= \pi_\mu(\xi)$. Using barycentric projections (cf. \cite{ags:book} Definition 5.4.2) one can show that, for $\gamma\in \Gamma_o(\mu,\nu)$, 
$$\iint_{\RD\times \RD} \langle\xi(x),y-x\rangle\,d\gamma(x,y)=\iint_{\RD\times \RD} \langle\pi_\mu(\xi)(x),y-x\rangle\,d\gamma(x,y).$$ 
Thus $\pi_\mu(\xi)\in\partial_{-} F(\mu) \cap \partial^{+} F(\mu)\cap T_{\mu}\M$ and it satisfies the analogue of Equation \ref{eq:differentiable0}. It can be shown that the gradient vector is unique, \textit{i.e.} that $\partial_{-} F(\mu) \cap \partial^{+} F(\mu)\cap T_{\mu}\M=\{\pi_\mu(\xi)\}$.

Finally, if the gradient vector exists for every $\mu\in\M$ we can define the \textit{differential} or \textit{exterior derivative} of $F$ to be the 1-form $dF$ determined, for any $\mu\in\M$ and $Y\in T_\mu\M$, by $dF(\mu)(Y):=\int_{\RD}\langle\nabla_\mu F,Y\rangle\,d\mu$. To simplify the notation we will sometimes write $Y(F)$ rather then $dF(Y)$.
\end{definition}

\begin{remark} \label{rem:Ftaylorexpansion}
Assume $F:\M\rightarrow\R$ is differentiable. Given $X\in \nabla C_c^\infty(\RD)$, let $\phi_t$ denote the flow of $X$. Fix $\mu\in\M$. 

\noindent(i) Set $\nu_t:= (\id+tX)_\# \mu.$ Then 
$$F(\nu_t)=F(\mu)+t\int_{\RD}\langle \nabla_{\mu} F,X\rangle d\mu + o(t).$$ 
(ii) Set $\mu_t:=\phi_{t\#}\mu$. If $||\nabla_\mu F(\mu)||_\mu$ is bounded on compact subsets of $\M$ then 
$$F(\mu_t)=F(\mu)+t\int_{\RD}\langle \nabla_{\mu} F,X\rangle d\mu + o(t).$$

\end{remark}
\proof{} The proof of (i) is a direct consequence of Equation \ref{eq:differentiable0} and of the fact that, if $r>0$ is small enough, $\bigl( \id \times (\id + tX) \bigr)_\#\mu \in \Gamma_o(\mu, \nu_t)$ for $t \in [-r,r]$. 

To prove (ii), set 
$$
A(s,t):=(1-s) (\id +t X) +s \phi_t.
$$
Notice that $||\phi_t -\id-t X||_\mu \leq t^2 ||(\nabla X) X||_\infty$ and that $(s,t) \rightarrow m(s,t):=A(s,t)_\# \mu$ defines a continuous map of the compact set $[0,1]\times [-r,r]$ into $\M.$ Hence the range of $m$ is compact so $||\nabla_\mu F(\mu)||_\mu$ is bounded there by a constant $C.$  One can use elementary arguments to conclude that $F$ is $C$-Lipschitz on the range of $m$, cf. \cite{hwakil:thesis} for details. Let $\gamma_t:= \bigl( (\id +t X) \times \phi_t\bigr)_\# \mu$. We have $\gamma_t \in \Gamma (\nu_t, \mu_t)$ so $W_2(\mu_t, \nu_t) \leq ||\phi_t -\id-t X||_\mu=O(t^2).$ We conclude that  
$$|F(\nu_t)-F(\mu_t)|\leq  C W_2(\mu_t, \nu_t) =O(t^2).$$
This, together with (i), yields (ii).\endproof

\begin{example} \label{e:linearisdifferentiable}
Fix $f\in\ccinf$ and let $F:\M\rightarrow\R$ be the corresponding linear function, as in Example \ref{e:deflinearforms}. Then $F$ is differentiable with gradient $\nabla_\mu F\equiv\nabla f$. Thus $dF$ is a linear 1-form on $\M$. Viceversa, according to our definitions every linear 1-form $\Lambda$ is \textit{exact}. In other words, if $\Lambda_\mu(X)=\int_{\R^D}\langle A,X\rangle d\mu$ for some $A=\nabla f$ then $\Lambda=dF$ for $F(\mu):=\intD f\,d\mu$. 
\end{example}

\begin{definition}\label{d:differential}
Let $\bar\Lambda$ be a pseudo 1-form on $\M$. We say that $\bar\Lambda$ is \textit{differentiable} if the following two conditions hold: 

(i) For all $X\in\cfields$, the function $\bar\Lambda(X):\M\rightarrow \R$ is differentiable. We can then define the \textit{exterior derivative} of $\bar\Lambda$ on pairs $X,Y\in\cfields$ by setting
\begin{equation}\label{eq:differential}
d\bar\Lambda(X,Y):=X\bar\Lambda(Y)-Y\bar\Lambda(X)-\bar\Lambda([X,Y]). 
\end{equation}
(ii) For all $\mu\in\M$, $d\bar\Lambda_\mu$ is continuous when restricted to $\gradfields$, \textit{i.e.} there exists $c_\mu\in\R$ such that 
\begin{equation*}
|d\bar\Lambda_\mu(\nabla f,\nabla g)|\leq c_\mu\|\nabla f\|_\mu\cdot\|\nabla g\|_\mu,\ \ \mbox{for all }\nabla f,\nabla g\in \gradfields.
\end{equation*}
Notice that condition (ii) implies that $d\bar\Lambda_\mu$ has a unique extension to $T_\mu\M\times T_\mu\M$.

Let $\Lambda$ be a 1-form on $\M$. Let $\bar\Lambda$ denote the associated pseudo 1-form, as in Equation \ref{eq:bar}. We say that $\Lambda$ is \textit{differentiable} if $\bar\Lambda$ is differentiable. We can then define its \textit{exterior derivative} by setting $d\Lambda:=d\bar\Lambda$.
\end{definition}

\begin{remark}\label{r:differentiability}
The assumption that $d\bar\Lambda$ satisfies the continuity assumption (ii) on $\gradfields$ implies that some form of cancelling occurs to eliminate first-order terms as in Equation \ref{eq:usuald}, cf. Remark \ref{r:dof1form}. Notice that $d\bar\Lambda$, restricted to $T_\mu\M\times T_\mu\M$, is a well-defined 2-form. On the other hand, Example \ref{e:linearimpliesdiff} shows that it is not natural to impose a continuity assumption on $\cfields$ so $d\bar\Lambda$ does not in general extend to a uniquely defined pseudo 2-form. 

If $\Lambda$ is differentiable in the above sense, it is natural to ask if $d\bar\Lambda_\mu(X,\cdot)=0$ for any $X\in \kerdiv$. It is not clear that this is the case.
\end{remark}
\begin{remark} One could also define a notion of differentiability for 1-forms by testing $\Lambda$ only against gradient vector fields $\nabla f\in\gradfields$. This is clearly a weaker condition, which would yield a very poor understanding of the differentiability of the associated pseudo 1-form $\bar\Lambda$. Indeed, assume $\Lambda$ is differentiable in the weaker sense and choose $X\in\cfields$. Then $\bar\Lambda(X)=\Lambda(\pi_\mu(X))$ and $\pi_\mu(X)$ depends on $\mu$. In particular, the differentiability of $\bar\Lambda$ is now related to the smoothness of the projection operators $\mu\mapsto\pi_\mu$. We will avoid this notion, using instead the stronger definition given in Definition \ref{d:differential}.
\end{remark}

\begin{example}\label{e:linearimpliesdiff} 
Assume $\bar\Lambda$ is a linear pseudo 1-form, \textit{i.e.} $\bar\Lambda(\cdot)=\intD \langle \bar A,\cdot\rangle d\mu$ for some $\bar A\in\cfields$. Then $\bar\Lambda$ is differentiable and, $\forall X,Y\in\cfields$, 
\begin{align*}
d\bar{\Lambda}(X,Y)
&= \int_{\RD}<\pi_\mu(\nabla \bar{A}^T\cdot Y+ \nabla Y^T\cdot
\bar{A}),X>-<\pi_\mu(\nabla \bar{A}^T\cdot X + \nabla X^T\cdot
\bar{A}),Y>d\mu\\\nonumber &\quad - \int_{\RD}<\nabla Y\cdot X-
\nabla X\cdot Y ,\bar{A}>d\mu.
\end{align*}
If $X,Y\in T_\mu\M$ then $d\bar{\Lambda}(X,Y) = \int_{\RD}<(\nabla\bar{A}-\nabla\bar{A}^T)X,Y>d\mu.$

\ 

\proof{}
Define $F_X: \M \rightarrow \R $  by
$$F_X(\mu):=\bar{\Lambda}(X)= \int_{\RD}<\bar{A},X>d\mu. $$
Let $\mu,\nu\in\M$ and $\gamma\in\Gamma_o(\mu,\nu)$. Then
\begin{align*}
F_X(\nu)-F_X(\mu) &= \int_{\RD}<\bar{A},X>d\nu-\int_{\RD}<\bar{A},X>d\mu\\
&=\iint_{\RD\times\RD}<\bar{A}(y),X(y)>-<\bar{A}(x),X(x)>d\gamma(x,y).
\end{align*}
Set $\phi:=<\bar{A},X> .$ Then $\phi\in C_c^\infty(\RD)$ so
\begin{align}\label{1}
F_X(\nu)-F_X(\mu) &=\iint_{\RD\times\RD}\phi(y)-\phi(x)d\gamma(x,y)\\
\nonumber & =\iint_{\RD\times\RD}(<\nabla\phi(x),y-x> + O(|x-y|^2) d\gamma(x,y)\\
\nonumber &= \iint_{\RD\times\RD}<\nabla\phi(x),y-x>d\gamma(x,y) + o(W_2(\mu,\nu)).
\end{align}
Equation \ref{1} shows that $F_X: \M \rightarrow \R$ is differentiable and that 
$\nabla_\mu F_X=\pi_\mu(\nabla \phi)$. Thus
\begin{align}\label{2}
Y\bar{\Lambda}(X) := d F_X(Y) &= \int_{\RD}<\nabla_\mu F_X,Y>d\mu\\
\nonumber &= \int_{\RD}<\pi_\mu(\nabla \bar{A}^T\cdot X + \nabla X^T\cdot \bar{A}),Y>d\mu.
\end{align}
Analogously,
\begin{align}\label{3}
X\bar{\Lambda}(Y) := d F_Y(X) &= \int_{\RD}<\nabla_\mu
F_Y,X>d\mu\\\nonumber &= \int_{\RD}<\pi_\mu(\nabla \bar{A}^T\cdot
Y+ \nabla Y^T\cdot \bar{A}),X>d\mu.
\end{align}
We combine Equations \ref{2} and \ref{3} to get
\begin{align}\label{4}
d\bar{\Lambda}(X,Y):&=
X\bar{\Lambda}(Y)-Y\bar{\Lambda}(X)-\bar{\Lambda}([X,Y])\nonumber\\
&= \int_{\RD}<\pi_\mu(\nabla \bar{A}^T\cdot Y+ \nabla Y^T\cdot
\bar{A}),X> d\mu\\
\nonumber &\quad -\int_{\RD}<\pi_\mu(\nabla \bar{A}^T\cdot X + \nabla X^T\cdot
\bar{A}),Y> d\mu\\ 
&\quad - \int_{\RD}<\nabla Y\cdot X-
\nabla X\cdot Y ,\bar{A}> d\mu.\nonumber
\end{align}
If $X,Y\in T_\mu\M$ then Equation \ref{4} simplifies to
\begin{align*}
d\bar{\Lambda}(X,Y) &= \int_{\RD}<(\nabla \bar{A}^T\cdot Y+ \nabla
Y^T\cdot \bar{A}),X>-<(\nabla \bar{A}^T\cdot X + \nabla X^T\cdot
\bar{A}),Y>d\mu\\
&\quad - \int_{\RD}<\nabla Y\cdot X-
\nabla X\cdot Y ,\bar{A}>d\mu\\ 
&= \int_{\RD}<(\nabla\bar{A}-\nabla\bar{A}^T)X,Y>d\mu.
\end{align*}\endproof

In Lemma \ref{le10.12.4} we will generalize this result to the class of \textit{regular} pseudo 1-forms.
\end{example}


\subsection{Discussion}\label{ss:section4discussion}
As explained in Section \ref{ss:tangentspacesbis}, we can think of $\M$ as the union of smooth manifolds $\mathcal{O}$. Each tangent space $T_\mu\M$ should then be thought of as the tangent space of $\mathcal{O}$ at the point $\mu$. Our notion of k-form $\Lambda$ is defined in terms of the dual tangent spaces, so each $\Lambda_{|\mathcal{O}}$ is, at least formally, a k-form on a smooth manifold in the usual sense. The logic behind our definition of the operator $d$ on 1-forms is as follows. As seen in Section \ref{ss:tangentspacesbis}, any $X\in\cfields$ defines a fundamental vector field on $\mathcal{O}$ (or on $\M)$. In particular we can think of the construction of fundamental vector fields as a canonical way of extending given tangent vectors $X$, $Y$ at any point $\mu\in\mathcal{O}$ to global tangent vector fields on $\mathcal{O}$. Combining Remark \ref{r:tangentspacevsalgebra} with Lemma \ref{l:adjointispushforward} shows that the construction of fundamental vector fields determines a Lie algebra homomorphism $\cfields\rightarrow \fields(\mathcal{O})$. Equation \ref{eq:differential} thus mimics Equation \ref{eq:usuald} for $k=1$. In Section \ref{ss:greenstheorem} we will study the corresponding first cohomology group. We can think of this as the \textit{de Rham cohomology} of the manifold $\mathcal{O}$. 

The notion of pseudo k-form is less standard, but also very natural. The finite-dimensional analogue of this notion is explained in Section \ref{ss:cohomology}. Roughly speaking, \textit{i.e.} up to $L^2_\mu$-closure, if we restrict our space of pseudo k-forms to any manifold $\mathcal{O}$ we obtain the space of maps $\mathcal{O}\rightarrow \Lambda^k(\mathfrak{g})$, where $\mathfrak{g}=\cfields$. Our definition of the operator $d$, given in Equation \ref{eq:differential}, should now be compared to Equation \ref{eq:orbitd}. Notice that the sign discrepancy between these equations is explained by the fact that the Lie bracket on $\mathfrak{g}$ is the opposite of the usual Lie bracket on $\cfields$, cf. Lemma \ref{l:adjointispushforward}. In this setting the key point is that each manifold $\mathcal{O}$ is actually the orbit of a group action. More specifically, we can identify it with a quotient of the group $\cdiffrD$. Proposition \ref{prop:2cohoms} then shows that pseudo k-forms are actually k-forms on the group, rather than on the manifold, endowed with a special ``invariance" property. In some sense the corresponding cohomology is more closely related to the orbit structure of the manifold $\mathcal{O}$ than to its topological structure. However Proposition \ref{prop:2cohoms} shows that, at least in finite dimensions, there is a simple relation between this \textit{invariant cohomology} and the usual de Rham cohomology of the manifold: for $k=1$, the latter is a subgroup of the former. Proposition \ref{prop:2cohoms} also shows that the operation of Equation \ref{eq:bar} is very natural from this point of view: up to the appropriate identifications, it coincides with the standard pull-back operation from k-forms on the quotient of the group to k-forms on the group.

It may be useful to emphasize that the identification between the orbit $\mathcal{O}$ and the quotient space is not canonical. The details involved in changing this identification are explained in Section \ref{ss:groupactions}.


\section{Calculus of pseudo differential 1-forms}\label{s:pseudocalculus}

Given a 1-form $\alpha$ on a finite-dimensional manifold, Green's formula relates the integral of $d\alpha$ along a surface to the integral of $\alpha$ along the boundary curves. In Section \ref{subse:smoothGreen} we show that an analogous result for $\M$ is rather simple if both the form and the surface satisfy certain regularity conditions. The conditions we need to impose on the form are rather mild: we investigate these in Sections \ref{ss:regularforms} and \ref{ss:morediff}, developing a general theory of \textit{regular} pseudo 1-forms. The conditions on the surface, instead, are very strong. In Section \ref{ss:greenstheorem} we thus prove a second version of Green's formula, valid only for certain surfaces we call \textit{annuli}. For these surfaces we need no extra regularity conditions, and Green's formula then leads to a proof that every closed regular 1-form is exact. 

\subsection{Green's formula for smooth surfaces and $1$-forms}\label{subse:smoothGreen}
Let $\Lambda$ be a differentiable 1-form on $\M$ in the sense of Definition \ref{d:differential}. Let $\bar\Lambda$ denote the associated pseudo 1-form in the sense of Equation \ref{eq:bar}. Set $||\Lambda_\mu||:= \sup_{v }\{ \Lambda_\mu(v): v \in T_\mu \M, ||v||_\mu \leq 1\}$. We assume that, for all compact subsets $\calK \subset \M$,
\begin{equation}\label{eq:assumpgreen1}
\sup_{\mu \in \calK }||\Lambda_\mu|| <\infty.
\end{equation}
We also assume that for all compact subsets $\calK \subset \M$ there exists a constant $C_{\calK}$ such that, for all $\mu, \nu \in \calK$  and $u \in C_b(\RD, \RD)$ such that $\nabla u$ is bounded, 
\begin{equation}\label{eq:assumpgreen1b}
|\bar\Lambda_\nu(u)-\bar\Lambda_\mu(u) | \leq C_{\calK} W_2(\mu, \nu) (||u||_\infty + ||\nabla u||_\infty).
\end{equation}
Set $||d \Lambda_\mu||$ to be the smallest nonnegative number $c_\mu$ such that 
$$|d\Lambda_\mu(\nabla f,\nabla g)| \leq  c_\mu ||\nabla f||_\mu\cdot ||\nabla g||_\mu,\ \ \mbox{for all }\nabla f, \nabla g \in \gradfields.$$
Now let $S:[0,1]\times[0,T]\rightarrow \M$ denote a map satisfying the following three regularity conditions: 

(i) For each $s \in [0,1],$ $S(s, \cdot) \in AC_2(0,T; \M)$ and, for each $t \in [0,T]$, $S(\cdot, t) \in AC_2(0,1; \M).$ 

(ii) Let $v(s, \cdot,\cdot)$ denote the velocity of minimal norm for $S(s, \cdot)$ and $w(\cdot, t,\cdot)$ denote the velocity of minimal norm for $S(\cdot,t).$ We assume that $v, w \in C^2( [0,1] \times [0,T] \times \RD, \RD)$ and that their derivatives up to third order are bounded. We further assume that $v$ and $w$ are gradient vector fields so that $\partial_s v$ and $\partial_t w$ are also gradients: this implies that $\Lambda$ and $\bar\Lambda$ coincide when evaluated on these fields.

(iii) $S$ takes values in the set of absolutely continuous measures. More specifically, $S(s,t)=\rho(s,t, \cdot) \calLD$ for some $\rho \in C^1( [0,1] \times [0,T] \times \RD)$ which is bounded with bounded derivatives.

\ 

Using Remark \ref{reac2curves}, Proposition \ref{prac2curves} and the bound on $v, w$ and on their derivatives, we find that $S$ is $1/2$-H\"older continuous. Hence its range is compact so $||\Lambda_{S(s,t)}||$ is bounded. We then use Equations \ref{eq:assumpgreen1},  \ref{eq:assumpgreen1b} and Taylor expansions for $w_{t+h}^s$ and $v_t^{s+h}$ to obtain that 
\begin{equation}\label{eq:conclusiongreen1}
\partial_t \Bigl( \Lambda_{S(s,t)} (w_t^s) \Bigr)_{|\bar s, \bar t} = v_{\bar t}^{\bar s}(\Lambda_{S(s,t)} (w_{\bar t}^{\bar s})) + \Lambda_{S(\bar s,\bar t)} (\partial_t w_t^s),
\end{equation} 
where we use the notation of Definition \ref{def:Fdifferentiable}. Analogously, 
\begin{equation}\label{eq:conclusiongreen2}
\partial_s \Bigl( \Lambda_{S(s,t)} (v_t^s) \Bigr)_{|\bar s, \bar t}= w_{\bar t}^{\bar s}(\Lambda_{S(s,t)} (v_{\bar t}^{\bar s})) + \Lambda_{S(\bar s,\bar t)} (\partial_s v_t^s).
\end{equation} 

\begin{lemma} \label{le:relation.velocities} For  $(s,t) \in (0,1) \times (0,T)$ we have $(\partial_t w_{t}^s-  \partial_s v_{t}^s\bigr) -[ w_{t}^s, v_{t}^s ] \in \mbox{Ker}(div_{S(s,t)}).$
\end{lemma}
\proof{} We have, in the sense of distributions,
\begin{equation}\label{eq:tangent2}
\partial_t \rho_t^s + \nabla \cdot (\rho_t^s v_t^s)=0, \quad \partial_s \rho_t^s + \nabla \cdot (\rho_t^s w_t^s)=0
\end{equation}
and so
$$
\nabla \cdot \partial_s (\rho_t^s v_t^s)=- \partial_s \partial_t \rho_t^s= \nabla \cdot (\partial_t \rho_t^s w_t^s).
$$
We use that $\rho$, $v$ and $w$ are smooth to conclude that
$$
\nabla \cdot
\Bigl(
v_t^s \partial_s \rho_t^s + \rho_t^s \partial_s v_t^s
\Bigr)
=
\nabla \cdot
\Bigl(
w_t^s \partial_t \rho_t^s + \rho_t^s \partial_t w_t^s
\Bigr).
$$
This implies that if $\varphi \in C_c^\infty(\RD)$ then
\begin{equation}\label{eq:tangent3}
\intD \langle \nabla \varphi, v_t^s \partial_s \rho_t^s + \rho_t^s \partial_s v_t^s \rangle
=
\intD \langle \nabla \varphi, w_t^s \partial_t \rho_t^s + \rho_t^s \partial_t w_t^s  \rangle.
\end{equation}
We use again that $\rho$, $v$ and $w$ are smooth to obtain that Equation \ref{eq:tangent2} holds pointwise. Hence, Equation \ref{eq:tangent3} implies
$$
\intD \langle \nabla \varphi, -v_t^s \nabla \cdot ( \rho_t^s w_t^s) + \rho_t^s \partial_s v_t^s \rangle
=
\intD \langle \nabla \varphi, -w_t^s \nabla \cdot ( \rho_t^s v_t^s) + \rho_t^s \partial_t w_t^s  \rangle.
$$
Rearranging, this leads to
$$
\intD \langle \nabla \varphi, \partial_s v_t^s - \partial_t w_t^s\rangle \rho_t^s d \calLD=
\intD \bigl \langle \nabla \varphi, v_t^s \bigr \rangle \nabla \cdot ( \rho_t^s w_t^s) -\bigl \langle \nabla \varphi, w_t^s \bigr \rangle \nabla \cdot ( \rho_t^s v_t^s).
$$
Integrating by parts and substituting $\rho_t^s \calLD$ with $S(s,t)$ we obtain
\begin{multline*}
\intD \langle \nabla \varphi, \partial_s v_t^s - \partial_t w_t^s\rangle dS(s,t)\\
\begin{split}
&=\intD\Bigl(
\bigl \langle \nabla^2 \varphi w_t^s+ (\nabla w_t^s)^T \nabla \varphi, v_t^s \bigr \rangle
-\bigl \langle \nabla^2 \varphi v_t^s+ (\nabla v_t^s)^T \nabla \varphi, w_t^s \bigr \rangle
\Bigr)dS(s,t)\\
&=\intD \Bigl \langle \nabla \varphi , [ v_t^s,  w_t^s] \Bigr\rangle dS(s,t).
\end{split}
\end{multline*}
Since $\varphi \in C_c^\infty(\RD)$ is arbitrary, the proof is finished. \endproof

\begin{prop}\label{surfacegreen1} For each $t \in (0,T)$ and $s \in (0,1)$ we have
\begin{equation*}
\partial_t \Bigl(  \Lambda_{S(s,t)} (w_t^s) \Bigr)- \partial_s \Bigl( \Lambda_{S(s,t)} (v_t^s) \Bigr) = d \Lambda_{S(s,t)}(v_t^s, w_t^s).
 \end{equation*} 
\end{prop}
\proof{}
We use Definition \ref{d:differential} and Equations \ref{eq:conclusiongreen1}, \ref{eq:conclusiongreen2} to find
\begin{align*}
d \Lambda_{S(\bar s,\bar t)}(v_{\bar t}^{\bar s}, w_{\bar t}^{\bar s}) &= d \bar\Lambda_{S(\bar s,\bar t)}(v_{\bar t}^{\bar s}, w_{\bar t}^{\bar s})\\ 
&= v_{\bar t}^{\bar s}(\bar\Lambda_{S(s,t)} (w_{\bar t}^{\bar s}))-
w_{\bar t}^{\bar s}(\bar\Lambda_{S(s,t)} (v_{\bar t}^{\bar s}))-\bar\Lambda_{S(\bar s,\bar t)}([v_{\bar t}^{\bar s}, w_{\bar t}^{\bar s}])\\
&= v_{\bar t}^{\bar s}(\Lambda_{S(s,t)} (w_{\bar t}^{\bar s}))-
w_{\bar t}^{\bar s}(\Lambda_{S(s,t)} (v_{\bar t}^{\bar s}))-\bar\Lambda_{S(\bar s,\bar t)}([v_{t}^{s}, w_{t}^{s}])\\
\begin{split}
&= \partial_t \Bigl( \Lambda_{S(s,t)} (w_t^s) \Bigr)_{|\bar s, \bar t}-\Lambda_{S(\bar s,\bar t)} (\partial_t w_t^s)-\partial_s \Bigl( \Lambda_{S(s,t)} (v_t^s) \Bigr)_{|\bar s, \bar t}\\
&\quad+\Lambda_{S(\bar s,\bar t)} (\partial_s v_t^s)-\bar\Lambda_{S(\bar s,\bar t)}([v_{t}^{s}, w_{t}^{s}])
\end{split}\\
\begin{split}
&= \partial_t \Bigl( \Lambda_{S(s,t)} (w_t^s) \Bigr)_{|\bar s, \bar t}-\partial_s \Bigl( \Lambda_{S(s,t)} (v_t^s) \Bigr)_{|\bar s, \bar t}\\
&\quad+\bar\Lambda_{S(\bar s,\bar t)}(\partial_s v_t^s-\partial_t w_t^s-[v_{t}^{s}, w_{t}^{s}]).
\end{split}
\end{align*}
We can now use Lemma \ref{le:relation.velocities} to conclude. 
\endproof
\begin{theorem}[Green's formula for smooth surfaces]\label{th:greengeneral} Let $S$ be a surface in $\M$ satisfying the above three conditions. Let $\partial S$ denote its boundary, defined as the union of the negatively oriented curves $S(0, \cdot),$ $S(\cdot,T)$ and the positively oriented curves $S(1, \cdot),$ $S(\cdot,0).$  Suppose that  $\mu \rightarrow ||d \Lambda_\mu||$ is also bounded on compact subsets of $\M.$ Then
$$
\int_{S} d\Lambda= \int_{\partial S} \Lambda.
$$
\end{theorem}
\proof{} Recall that $v_t^s$, $w_t^s$ and their derivatives are bounded. This, together with  Equations \ref{eq:assumpgreen1} and  \ref{eq:assumpgreen1b}, implies that the functions $(s,t) \rightarrow \Lambda_{S(s,t)}(v_t^s)$ and $(s,t) \rightarrow  \Lambda_{S(s,t)}(w_t^s)$ are continuous. Hence, by Proposition \ref{surfacegreen1}, $(s,t) \rightarrow d \Lambda_{S(s,t)}(v_t^s, w_t^s)$ is Borel measurable as it is a limit of quotients of continuous functions. The fact that $\mu \rightarrow ||d \Lambda_\mu||$ is  bounded on compact subsets of $\M$ gives that $(s,t) \rightarrow d \Lambda_{S(s,t)}(v_t^s, w_t^s)$ is bounded. The rest of the proof of this theorem is identical to that of Theorem \ref{th:greencurveannulus} when we use Proposition \ref{surfacegreen1} in place of Corollary \ref{le:lipschitzgreen3}. \endproof

The regularity conditions we have imposed on $S$ to obtain Theorem \ref{th:greengeneral} are very strong. In particular, given an $AC_2$ curve $\sigma$ of absolutely continuous measures in $\M$, it is not clear if there exists any surface $S$ satisfying these assumptions and whose boundary is $\sigma$. It is thus difficult to use Theorem \ref{th:greengeneral} to reach any conclusion about the de Rham cohomology of $\M$. Such conclusions will however be obtained in Section \ref{ss:greenstheorem}, based on a different version of Green's theorem.


\subsection{Regularity and differentiability of pseudo 1-forms}\label{ss:regularforms}
The goal of this section is to introduce a regularity condition for pseudo 1-forms which guarantees differentiability. It will also ensure the validity of assumptions such as Equations \ref{eq:assumpgreen1} and \ref{eq:assumpgreen1b}.

\begin{definition} \label{def:regularform}
Let $\bar\Lambda_\mu:=\int_{\R^D}\langle \bar A_\mu,\cdot\rangle d\mu$ be a pseudo 1-form on $\M$. We say that $\bar\Lambda$ is \textit{regular} if for each $\mu \in \M$ there exists a Borel field of $D \times D$ matrices $B_\mu \in L^\infty(\RD \times \RD, \mu)$ and a function $O_\mu\in C(\R)$ with $O_\mu(0)=0$ such that
\begin{multline}
\sup_{\gamma \in \Gamma_o(\mu, \nu)} \Bigl\{ \int_{\RD \times \RD} |\bar A_\nu(y)-\bar A_\mu(x)-B_\mu(x)(y-x)|^2 d\gamma(x,y)\Bigr\}\\
\leq W^2_2(\mu, \nu) \min\{ O_\mu(W_2(\mu, \nu)), c(\bar\Lambda)\}^2.\label{eqcontdiff}
\end{multline}
where as usual $\Gamma_o(\mu, \nu)$ denotes the set of minimizers in Equation \ref{eq:wasserstein} and $c(\bar\Lambda)>0$ is a constant independent of $\mu.$  We also assume that $||B_\mu||_\mu$ is uniformly bounded. Taking $c(\bar\Lambda)$ large enough, there is no loss of generality in assuming that
\begin{equation} \label{eq:contdiffbis}
\sup_{\mu \in \M} ||B_\mu||_\mu \leq c(\bar\Lambda).
\end{equation}
Let $\Lambda$ be a 1-form on $\M$. Let $\bar\Lambda$ denote the associated pseudo 1-form, as in Equation \ref{eq:bar}. We say that $\Lambda$ is \textit{regular} if $\bar\Lambda$ is regular.
\end{definition}

\begin{remark} Some of the assumptions in Definition \ref{def:regularform} could be weakened for the purposes of this paper. We make these choices simply to avoid introducing more notation and to shorten some computations.
\end{remark}

\begin{example}\label{e:linearisregular}
Every linear pseudo 1-form is regular. In other words, given $\bar A \in \cfields$, if we define $\bar \Lambda_\mu(Y):= \intD \langle \bar A, Y\rangle d\mu$ then $\bar\Lambda$ is regular. Indeed,  setting $B_\mu:= \nabla \bar A$ one can use Taylor expansion and the fact that the second derivatives of $A$ are bounded to obtain Equation \ref{eqcontdiff}. 
\end{example}

\begin{remark} Even if Equation \ref{eqcontdiff} holds for $\bar A_\mu$, it does not necessarily hold for $A_\mu:=\pi_\mu(\bar A)$. This implies that, in general, it is not clear what regularity properties might hold for the 1-form obtained by restricting a regular pseudo 1-form. This is true even in the simplest case where $\bar\Lambda$ is as in Example \ref{e:linearisregular}. The case in which $\bar\Lambda$ is related to $\Lambda$ as in Equation \ref{eq:bar} is an obvious exception: in this case, according to Definition \ref{def:regularform}, $\bar\Lambda$ is regular iff $\Lambda$ is regular.
\end{remark}

From now till the end of Section \ref{s:pseudocalculus} we assume $\bar\Lambda$ is a regular pseudo 1-form on $\M$ and we use the notation $\bar{A}_\mu$, $B_\mu$ as in Definition \ref{def:regularform}.

\begin{remark}\label{re:contdiffbis} 
If $\mu, \nu \in \M$, $X \in L^2(\mu),$ $Y\in L^2(\nu)$ and $\gamma \in \Gamma_o(\mu, \nu)$ then
\begin{multline}\label{eq:contdifflate0}
\bar \Lambda_{\nu}(Y) - \bar \Lambda_\mu(X)
- \intDD \Bigl( \langle \bar A_{\mu}(x) , Y (y)-X (x)\rangle+  \langle B_{\mu}(x)(y-x) , Y (y)\rangle \Bigr) d\gamma(x,y)\\
= \intDD \langle \bar A_{\nu}(y) -\bar A_\mu(x)-B_\mu(x)(y-x) , Y (y)\rangle d\gamma(x,y). 
\end{multline}
By Equation \ref{eqcontdiff} and H\"older's inequality,
\begin{equation} \label{eq:contdifflate1}
\Bigl | \intDD \langle \bar A_{\nu}(y) -\bar A_\mu(x)-B_\mu(x)(y-x) , Y (y)\rangle \Bigr |
\leq W_2(\mu, \nu) c(\bar\Lambda) \, ||Y||_\nu.
\end{equation}
Similarly, Equation \ref{eq:contdiffbis} and H\"older's inequality yield 
\begin{equation} \label{eq:contdifflate2}
\Bigl |  \intDD \langle B_\mu(x)(y-x) , Y (y)\rangle \Bigr |
\leq W_2(\mu, \nu) c(\bar\Lambda) \, ||Y||_\nu.
\end{equation}
We use Equations \ref{eq:contdifflate1} and  \ref{eq:contdifflate2} to obtain
\begin{equation} \label{eq:contdiff3}
\Bigl|\bar \Lambda_{\nu}(Y) - \bar \Lambda_\mu(X)- \intDD \langle \bar A_{\mu}(x) , Y (y)-X (x)\rangle d\gamma(x,y) \Bigr| \leq 2c(\bar\Lambda)W_2(\mu, \nu) \, ||Y||_\nu.
\end{equation}
\end{remark}

\begin{remark}\label{re:LipschitzofF} 
Let $Y \in C_c^1(\RD)$ and define $F(\mu):=\bar \Lambda_\mu(Y).$ Then
$$
|F(\nu)-F(\mu)| \leq W_2(\nu, \mu) \Bigl( ||\bar A_\nu||_\nu ||\nabla Y||_\infty+ 2c(\bar\Lambda) ||Y||_\infty \Bigr).
$$
\end{remark}
\proof{} By H\"older's inequality, 
$$
\Bigl| \intDD \langle \bar A_{\mu}(x) , Y (y)-Y (x)\rangle d\gamma(x,y) \Bigr| \leq || \bar A_{\mu}||_\mu ||\nabla Y||_\infty W_2(\nu, \mu).
$$
We apply Remark \ref{re:contdiffbis} with $Y=X$ and we exchange the role of $\mu$ and $\nu$ to conclude the proof. \endproof

\begin{lemma}\label{le:boundedA}  
The function
$$\M\rightarrow\R,\ \ \mu \mapsto ||\bar A_\mu||_\mu$$ 
is continuous on $\M$ and bounded on bounded subsets of $\M.$ Suppose $S:[r,1] \times [a,b]\rightarrow \M$ is continuous. Then
$$
\sup_{(s,t) \in [r,1] \times [a,b]} ||\bar A_{S(s,t)}||_{S(s,t)}<\infty.
$$
\end{lemma}
\proof{} Fix $\mu_0 \in \M.$ For each $\mu \in \M$ we choose $\gamma_\mu \in \Gamma_o(\mu_0, \mu).$  We have
$$
\bigl|\;  ||\bar A_\mu||_\mu- ||\bar A_{\mu_0}||_{\mu_0} \bigr|=\bigl|\; ||\bar A_\mu(y)||_{\gamma_\mu}- ||\bar A_{\mu_0}(x)||_{\gamma_\mu} \bigr| \leq ||\bar A_\mu(y)-\bar A_{\mu_0}(x)||_{\gamma_\mu}.
$$
This, together with Equations \ref{eqcontdiff} and \ref{eq:contdiffbis}, yields 
$$
\Bigl|  ||\bar A_\mu||_\mu- ||\bar A_{\mu_0}||_{\mu_0}  \Bigr|   \leq ||B_{\mu_0}(x)(y-x)||_{\gamma_\mu}+ c(\bar\Lambda)W_2(\mu_0, \mu)  \leq 2 c(\bar\Lambda)W_2(\mu_0, \mu).
$$
To obtain the last inequality we have used H\"older's inequality. This proves the first claim.

Notice that $(s,t) \rightarrow  ||\bar A_{S(s,t)}||_{S(s,t)}$ is the composition of two continuous functions and is defined on the compact set $[r,1] \times [a,b]$. Hence it achieves its maximum. \endproof

\begin{lemma}\label{rm10.12.4} 
Let $Y \in C_c^2(\RD)$ and define $F(\mu):=\bar \Lambda_{\mu}(Y).$ Then $F$ is differentiable with gradient $\nabla _\mu F=\pi_\mu(\nabla Y^T (x) \bar A_{\mu}(x) + B_{\mu}^T(x) Y(x))$.

Furthermore, assume $X \in \nabla C_c^2(\RD)$ and let $\varphi_t(x)=x+ tX(x) +t \bar O_t(x)$, where $\bar O_t$ is any continuous function on $\RD$ such that $||\bar O_t ||_\infty$ tends to $0$ as $t$ tends to $0.$ Set $\mu_t:=\varphi(t,\cdot)_\# \mu$. Then 
\begin{equation}\label{eq10.13.11a}
F(\mu_t)  =F(\mu) + t \intD \Bigl[ \langle \bar A_{\mu}(x), \nabla Y(x) X(x) \rangle +
\langle B_{\mu}(x) X(x), Y(x) \rangle \Bigr] d\mu(x) + o(t).
\end{equation}  
\end{lemma}
\proof{} Choose $\mu,\nu \in \mathcal{M}$ and $\gamma \in \Gamma_o(\mu,\nu)$. As in Remark \ref{re:contdiffbis},
\begin{multline*}
\bar{ \Lambda}_{\nu}(Y) - \bar {\Lambda}_\mu(Y) - \intDD \Bigl( \langle \bar A_{\mu}(x) , Y (y)-Y (x)\rangle+  \langle B_{\mu}(x)(y-x) , Y (y)\rangle \Bigr) d\gamma(x,y)\\
= \intDD \langle \bar A_{\nu}(y) -\bar A_\mu(x)-B_\mu(x)(y-x) , Y (y)\rangle d\gamma(x,y).
\end{multline*}
By Equation \ref{eqcontdiff} and H\"older's inequality,
\begin{equation*}
\Bigl | \intDD \langle \bar A_{\nu}(y) -\bar A_\mu(x)-B_\mu(x)(y-x) , Y (y)\rangle \Bigr |
\leq o(W_2(\mu, \nu))  \, ||Y||_\nu.
\end{equation*}
Since $Y\in C_c^2(\R^D)$ we can write $Y(y)=Y(x)+\nabla Y(x)(y-x)+R(x,y)(y-x)^2$, for some continuous field of vector-valued 2-tensors $R=R(x,y)$. In particular, $R$ has compact support and depends on the second derivatives of $Y$. Then
\begin{align}
\intDD  \langle \bar A_{\mu}(x) , Y (y)-Y(x)\rangle d\gamma(x,y) &= \intDD \langle \bar A_{\mu}(x) , \nabla Y(x)(y-x)\rangle  d\gamma(x,y)\nonumber\\
&\quad + \intDD \langle \bar A_\mu(x),R(y-x)^2\rangle d\gamma(x,y).\label{eq:estimate}
\end{align}
We now want to show that the term in Equation \ref{eq:estimate} is of the form $o(W_2(\mu,\nu))$ as $\nu$ tends to $\mu$. For any $\epsilon>0$, choose a smooth compactly supported vector field $Z=Z(x)$ such that $\|\bar A_\mu-Z\|_\mu<\epsilon$. Then by transposing the matrix $R(y-x)$ and using H\"older's inequality we obtain
\begin{align*}
|\intDD \langle \bar A_\mu(x),R(y-x)^2\rangle d\gamma(x,y)| &\leq \intDD |\langle (R(y-x))^T(\bar A_\mu(x)-Z(x)),y-x\rangle | d\gamma(x,y)\\ 
&\quad+ \intDD |\langle Z(x),R(y-x)^2\rangle | d\gamma(x,y)\\
&\leq \epsilon\,\|(R(y-x))\|_\infty W_2(\mu,\nu)+\|Z\|_\infty\|R\|_\infty W_2^2(\mu,\nu).
\end{align*} 
Since $\epsilon$ and $\|Z\|_\infty$ are independent of $\nu$, this gives the required estimate. Likewise,
\begin{align*}
\intDD  \langle B_{\mu}(x)(y-x) , Y (y)\rangle d\gamma(x,y) &= \intDD  \langle B_{\mu}(x)(y-x) , Y (y)- Y(x)\rangle d\gamma(x,y)\\
&\quad + \intDD  \langle B_{\mu}(x)(y-x) ,  Y(x)\rangle d\gamma(x,y)\\
&= \intDD  \langle B_{\mu}(x)(y-x) ,  Y(x)\rangle d\gamma(x,y) + o(W_2(\mu,\nu)).
\end{align*}
Combining these results shows that
\begin{equation*}
\bar{ \Lambda}_{\nu}(Y) = \bar {\Lambda}_\mu(Y)+ \intDD  \langle \nabla Y^T (x) \bar A_{\mu}(x) + B_{\mu}^T(x) Y(x) , y-x\rangle  d\gamma(x,y) + o(W_2(\mu,\nu)).
\end{equation*}
As in Definition \ref{def:Fdifferentiable}, this proves that $F$ is differentiable and that $\nabla_\mu F = \pi_\mu(\nabla Y^T (x) \bar A_{\mu}(x) + B_{\mu}^T(x) Y(x)).$

Now assume that $\phi_t$ is the flow of $X$. Notice that the curve $t \rightarrow \mu_t$ belongs to $AC_2(-r,r; \M)$ for $r>0$. We could choose for instance $r=1.$ Hence the curve is continuous on $[-1, 1].$ By Lemma \ref{le:boundedA}, the  composed   function $t \rightarrow ||\bar A_{\mu_t}||_{\mu_t}$ is also continuous. Hence its range is compact in $\R$, so there exists $\bar C>0$ such that $||\bar A_{\mu_t}||_{\mu_t} \leq \bar C$ for all $t \in [-1,1].$ We may now use Remark \ref{rem:Ftaylorexpansion} to conclude. 

The general case of $\phi_t$ as in the statement of Lemma \ref{rm10.12.4} can be studied using analogous methods.
\endproof

\begin{lemma}\label{le10.12.4} Any regular pseudo 1-form $\bar\Lambda$ is differentiable in the sense of Definition \ref{d:differential}. Furthermore, $\forall X, Y \in \TmuM$,
\begin{equation} \label{eq10.13.9}
d \bar \Lambda_\mu(X,Y)= \int_{\RD} \langle (B_\mu -B^T_\mu) X, Y \rangle d\mu.
\end{equation}
\end{lemma}

\proof{} The fact that, for each $Y\in\cfields$, $\bar\Lambda(Y)$ is a differentiable function on $\M$ follows from Lemma \ref{rm10.12.4}. Lemma \ref{rm10.12.4} also gives an expression for the gradient of this function. Using this expression it is simple to check that, for $X,Y\in\gradfields$,
\begin{equation}\label{eq:itisdiff}
X\bar\Lambda(Y)-Y\bar\Lambda(X)-\bar\Lambda([X,Y])= \int_{\RD} \langle (B_\mu -B^T_\mu) X, Y \rangle d\mu.
\end{equation}
This proves that $\bar\Lambda$ is differentiable. By continuity, the same expression holds for any $X,Y\in T_\mu\M$.
\endproof


\subsection{Regular forms and absolutely continuous curves} \label{ss:morediff}

The goal of this section is to study the regularity and integrability properties of regular pseudo 1-forms evaluated along curves in $AC_2(a,b;\M)$. 

\begin{lemma}\label{le10.12.2}
Assume $\{\mu_\epsilon\}_{\epsilon \in E} \subset \M$ and $v_\epsilon \in L^2( \mu_\epsilon)$ are such that $C:=\sup_{\epsilon \in E} ||v_\epsilon||_{L^2(\mu_\epsilon)}$ is finite. Assume $\{\mu_\epsilon\}_{\epsilon \in E}$ converges to $\mu$ in $\M$ as $\epsilon$ tends to $0$ and that there exists $v \in \lfields$ such that $\{v_\epsilon \mu_\epsilon\}_{\epsilon \in E}$ converges weak-$\ast$  to $v \mu$ as $\epsilon\rightarrow 0.$  If  $\gamma_\epsilon \in \Gamma_o(\mu, \mu_\epsilon)$ then $\lim_{\epsilon \rightarrow 0} a_\epsilon=0$, where  $a_\epsilon=\int_{\RDD} \langle \bar A_{\mu}(x) ,  v_\epsilon(y) -v (x) \rangle d \gamma_\epsilon(x,y).$
\end{lemma}
\proof{} It is easy to obtain that $||v||_{\lfields} \leq C.$ Let $\gamma_\epsilon \in \Gamma_o(\mu, \mu_\epsilon)$ and $\xi \in \cfields.$ Then there exists  a bounded function $C_\xi \in C(\RDD)$ and a real number $M$ such that
\begin{equation} \label{10.13.2}
\xi(x)- \xi(y)= \nabla \xi(y) (x-y) + |x-y|^2 C_\xi(x,y), \quad |C_\xi(x,y)| \leq M,
\end{equation}
for $x, y \in \RD.$ We use the first equality in Equation \ref{10.13.2} to obtain that
\begin{align*}
\langle \bar A_{\mu}(x) , v_\epsilon (y)-v (x)\rangle
&= \langle \bar  A_{\mu}(x) -\xi(x) , v_\epsilon (y)-v(x) \rangle
+ \langle \xi(y) ,v_\epsilon (y) \rangle-\langle \xi(x) , v(x) \rangle\\
&\quad + \langle \nabla \xi(y) (x-y)+ |x-y|^2 C_\xi(x,y)  , v_\epsilon (y)\rangle.
\end{align*}
Hence,
\begin{align}
|a_\epsilon|
& \leq || \bar A_{\mu}(x)- \xi(x)||_{L^2(\gamma_\epsilon)} ||v_\epsilon(y)-  v(x)||_{L^2(\gamma_\epsilon)}
+ b_\epsilon \nonumber\\
&\quad+
| \intDD \langle\bigl( \nabla \xi(y)(x-y) + |x-y|^2C_\xi(x,y) \bigr),v_\epsilon(y)\rangle d\gamma_\epsilon(x,y) |
\label{eq10.13.4}.
\end{align}
Above, we have set $b_\epsilon:= | \intDD \bigl( \langle \xi(y), v_\epsilon (y)\rangle-\langle \xi(x) ,v (x) \rangle \bigr) d\gamma_\epsilon(x,y)| .$ By the second inequality in Equation \ref{10.13.2} and by Equation \ref{eq10.13.4}, 
\begin{equation} \label{eq10.13.5}
|a_\epsilon|   \leq 2C || \bar A_{\mu}- \xi||_{\lfields} + b_\epsilon + ||\nabla \xi||_\infty W_2(\mu, \mu_\epsilon) + M W^2_2(\mu, \mu_\epsilon).
\end{equation}
By assumption $\{W_2(\mu, \mu_\epsilon) \}_{\epsilon \in E}$ tends to $0$ and $\{b_\epsilon \}_{\epsilon \in E}$ tends to $0$ as $\epsilon$ tends to $0$. These facts, together with Equation \ref{eq10.13.5}, yield $\limsup_{\epsilon \rightarrow 0} |a_\epsilon| \leq 2C || \bar A_{\mu}- \xi||_{\lfields}$ for arbitrary $\xi \in \cfields.$ We use that  $\cfields$ is dense in  $\lfields$ to conclude that $\lim_{\epsilon \rightarrow 0} a_\epsilon=0.$
\endproof

\begin{corollary}\label{co10.12.3} Assume $\{\mu_\epsilon\}_{\epsilon \in E} \subset \M$, $\mu$,  $v_\epsilon \in L^2( \mu_\epsilon)$ and $v$ satisfy the assumptions of Lemma \ref{le10.12.2}. Then $\lim_{\epsilon \rightarrow 0} \bar \Lambda_{\mu_\epsilon}(v_\epsilon)= \bar \Lambda_{\mu}(v).$
\end{corollary}

\proof{} Let $\gamma_\epsilon \in \Gamma_o(\mu, \mu_\epsilon).$ Observe that
\begin{align}
\langle \bar A_{\mu_\epsilon}(y), v_\epsilon (y)\rangle - \langle \bar  A_{\mu}(x) , v (x) \rangle
&=\langle \bar A_{\mu}(x) , v_\epsilon (y)-v (x)\rangle +\langle B_\mu(x) (y-x) ,v_\epsilon (y)\rangle \nonumber\\
&\quad+ \Bigl\langle \bar A_{\mu_\epsilon}(y) -\bar A_{\mu}(x) -B_\mu(x) (y-x) , v_\epsilon (y) \Bigr \rangle
\label{eq10.13.1}.
\end{align}
We now integrate Equation \ref{eq10.13.1} over $\RDD$ and use Equations \ref{eqcontdiff}, \ref{eq:contdiffbis} and the fact that $\gamma_\epsilon \in \Gamma_o(\mu, \mu_\epsilon).$ We obtain
\begin{align}
|\bar \Lambda_{\mu_\epsilon}(v_\epsilon)- \bar \Lambda_{\mu}(v) |
& \leq |a_\epsilon| +||B_\mu||_{L^\infty(\mu)} W_2(\mu, \mu_\epsilon) || v_\epsilon||_{\mu_\epsilon}  + o(W_2(\mu, \mu_\epsilon))  || v_\epsilon||_{\mu_\epsilon} \nonumber\\
& \leq
|a_\epsilon| +C ||B_\mu||_{L^\infty(\mu)} W_2(\mu, \mu_\epsilon) + C \, o(W_2(\mu, \mu_\epsilon))
\label{eq10.13.7}.
\end{align}
Letting $\epsilon$ tend to $0$ in Equation \ref{eq10.13.7} we conclude the proof of the corollary.
\endproof

\begin{lemma}[continuity of $\bar \Lambda_{\sigma_t}(X_t)$] \label{le10.27.7}
Suppose $\sigma \in AC_2(a,b;\M)$. If $X \in C((a,b) \times \RD, \RD)$ then  $\lambda(t):=\bar \Lambda_{\sigma_t}(X_t)$ is continuous on $(a,b).$
\end{lemma}

\proof{} Fix $t \in (a,b)$ so that $t$ belongs to the interior of a compact set $K^* \subset (a,b).$  Let $\varphi \in C_c(\RD,\RD)$ and denote by $K$ a compact set containing its  support.  Observe that $X$ is uniformly continuous on $K^* \times K$ so
\begin{equation}\label{01.18.1}
\limsup_{h \rightarrow 0} |\intD \langle \varphi(x), X_{t+h}(x)- X_{t}(x)  \rangle d\sigma_{t+h} (x)| \leq \limsup_{h \rightarrow 0} || \varphi||_\infty \sup_{x \in K} |X_{t+h}(x)- X_{t}(x)| =0.
\end{equation}
Since  $\langle X_t, \varphi\rangle \in C_c$ and $\sigma$ is continuous at $t$ by Remark  \ref{reac2curves},   we also see that
\begin{equation}\label{01.18.2}
\lim_{h \rightarrow 0}  \intD \langle \varphi(x), X_{t}(x)  \rangle d\sigma_{t+h} (x)
= \intD \langle \varphi(x), X_{t}(x)  \rangle d\sigma_{t} (x).
\end{equation}
Since $\varphi \in C_c(\RD,\RD)$ is arbitrary, Equations \ref{01.18.1} and \ref{01.18.2} show that $\{X_{t+h} \sigma_{t+h}\}_{h>0}$ converges weak-$\ast$ to $\sigma_t X_t$ as $h$ tends to zero. Corollary \ref{co10.12.3} yields that $\lambda$ is continuous at $t.$ 

\endproof

\begin{remark}
Using the techniques of Lemma \ref{le:lipschitzgreen2} one could further prove that if $\sigma \in AC_2(a,b; \M)$  and $X$ is sufficiently regular then $\lambda(t):=\bar \Lambda_{\sigma_t}(X_t)$ is Lipschitz and $\mathcal{L}^1$-almost everywhere differentiable.
\end{remark}

We now assume that $\eta_D^\epsilon \in C^\infty(\RD)$ is a mollifier :  $\eta_D^\epsilon(x)= 1/\epsilon^D \eta(x / \epsilon)$, for some bounded symmetric function $\eta \in C^\infty( \RD)$ whose derivatives of all orders are bounded. We also impose that $\eta > 0,$ $\intD |x|^2 \eta(x)dx<\infty$ and $\intD \eta =1.$ We fix $\mu \in \M$ and  define $f^\epsilon(x):= \intD \eta^\epsilon_D (x-y) d\mu(y).$ Observe that $f^\epsilon \in C^\infty(\RD)$ is bounded, all its derivatives are bounded and $\intD f^\epsilon =1.$

We suppose that $\eta^\epsilon_1 \in C^\infty(\R)$ is a standard mollifier: $\eta_1^\epsilon(t)= 1/\epsilon \eta_1(t / \epsilon)$, for some bounded symmetric function $\eta_1 \in C^\infty( \R)$ which is positive on $(-1,1)$ and vanishes outside $(-1,1)$.   We also impose that $\intR  \eta_1 =1$ and assume that $|\epsilon| <1.$

Suppose $\sigma \in AC_2(a,b; \M)$ and $v:(a,b) \times \RD \rightarrow \RD$ is a velocity associated to $\sigma$ so that $t \rightarrow ||v_t||_{\sigma_t} \in L^\infty(a,b).$ Suppose that for each $t \in (a,b)$ there exists $\rho_t>0$ such that $\sigma_t= \rho_t \calLD.$

We can extend $\sigma$ and $v$ in time on an interval larger than $[a,b].$ For instance,  set $\tilde \sigma_t=\sigma_a$ for $t \in (a-1,a)$ and set $\tilde \sigma_t=\sigma_b$ for $t \in (b,b+1)$. Observe that $\tilde \sigma \in AC_2(a-1,b+1; \M)$ and we have a velocity $\tilde v$ associated to $\tilde \sigma$ such that $\tilde v_t=v_t$ for $t \in [a,b].$ We can choose $\tilde v$ such that $||\tilde v_t||^2_{\tilde \sigma_t}=0$ for $t $ outside $(a,b).$ In particular, $\int_{a-1}^{b-1} ||\tilde v_t||^2_{\tilde \sigma_t}dt = \int_{a}^{b} || v_t||^2_{ \sigma_t}dt.$ In the sequel we won't distinguish between $\sigma$, $\tilde \sigma$ on the one hand and $v$, $\tilde v$ on the other hand. This extension becomes useful when we try to define $\rho_t^\epsilon$ as it appears in Equation \ref{eqoct27.5}. The new density functions are meaningful if we substitute $\sigma$ by $\tilde \sigma$ and impose that $\epsilon \in (0,1).$

For $\epsilon \in (0,1)$, set
\begin{equation}\label{eqoct27.5}
\rho_t^\epsilon(x):= \intR \eta^\epsilon_1(t -\tau) \rho_\tau(x) d\tau, \; \sigma^\epsilon_t:= \rho_t^\epsilon \calLD, \; \rho_t^\epsilon(x) v_t^\epsilon(x) := \intR \eta^\epsilon_1(t -\tau) \rho_\tau(x) v_\tau(x) d\tau.
\end{equation}
Note that $\rho_t^\epsilon(x)>0$ for all $t \in (a,b)$ and $x \in \RD$ and $\rho_t^\epsilon$ is a probability density. Also, $v^\epsilon:(a,b) \times \RD \rightarrow \RD$ is a velocity associated to $\sigma^\epsilon.$ In the sequel we set 
$$C^2:= \intD |x|^2 \eta(x) dx, \quad C_1= \intR \eta_1(\tau) \tau d\tau, \quad C_v:=\sup_{\tau \in (a-1, b+1)} ||v_\tau||_{\sigma_\tau}.$$
\begin{lemma}\label{leapproximation1} We assume that for each $t \in (a,b)$ there exists $\rho_t>0$ such that $\sigma_t= \rho_t \calLD.$ Then $\sigma^\epsilon \in AC_2(a,b; \M)$.  For   $a<s<t<b$,
$$
 (i)\;\; W_2(\mu, f^\epsilon \calLD) \leq \epsilon C, \;\;\;\; (ii)\; ||v^\epsilon_t||_{\sigma^\epsilon_t} \leq  C_v
 \;\;\;\hbox{and} \;\;\;\;
(iii)\; W_2(\sigma^\epsilon_t, \sigma_t) \leq \epsilon C_1 C_v .
$$
\end{lemma}
\proof{} We denote by ${\mathcal{U}}$ the set of pairs $(u,v)$ such that $u, v \in C(\RD)$ are bounded and $u(x) + v(y) \leq |x-y|^2$ for all $x, y \in \RD.$ Fix $(u,v) \in {\mathcal{U}}$. By Fubini's theorem one gets the well known identity
\begin{equation}\label{oct20.1}
\intD u(x) f^\epsilon(x)dx =\intD  d\mu(y) \intD u(x) \eta_\epsilon(x-y)dx.
\end{equation}
Since $v(y)=\intD v(y) \eta_\epsilon(x-y)dx$, Equation \ref{oct20.1} yields that
\begin{align}
\intD u(x) f^\epsilon(x)dx + \intD v(y)d\mu(y)
& =\intD d\mu(y) \intD \eta_\epsilon(x-y) \bigl(u(x) +v(y)  \bigr) dx
\nonumber\\
& \leq \intD d\mu(y) \intD \eta_\epsilon(x-y) |x-y|^2 dx
\label{oct20.2f}\\
& =\intD d\mu(y) \intD {1 \over \epsilon^D} \eta({z \over \epsilon}) |z|^2 dz= C^2 \epsilon^2 .\nonumber
\end{align}
To obtain Equation \ref{oct20.2f} we have used that $(u,v)\in\mathcal{U}$. We have proven that $\intD u(x) f^\epsilon(x)\,dx + \intD v(y)\,d\mu(y) \leq  C^2 \epsilon^2$ for arbitrary $(u,v) \in\mathcal{U}$. Thanks to the dual formulation of the Wasserstein distance Equation \ref{eq:dualwasserstein}, we conclude the proof of (i).

Notice that for each $t\in (a,b)$ and $x \in \RD$, $\eta^\epsilon_1(t -\tau) \rho_\tau(x)/ \rho_t^\epsilon(x)$ is a probability density on $\R.$ Hence, by Jensen's inequality,
$$|v_t^\epsilon(x)|^2= \Bigl|   1/  \rho_t^\epsilon(x)  \intR \eta^\epsilon_1(t -\tau) \rho_\tau(x) v_\tau(x) d\tau   \Bigr|^2 \leq 1 / \rho_t^\epsilon(x)  \intR \eta^\epsilon_1(t -\tau) \rho_\tau(x) |v_\tau(x)|^2 d\tau.$$
We multiply both sides of the previous inequality by $\rho_t^\epsilon(x)$. We then integrate the subsequent inequality over $\R^D$ and use Fubini's theorem to conclude the proof of (ii).

We use (ii) and Remark \ref{reac2curves} (i) to obtain that $\sigma^\epsilon \in AC_2(a,b; \M).$ We have
$$\intD u(x) d\sigma^\epsilon_t(x)=\intD u(x) dx \intR \eta^\epsilon_1(\tau) \rho_{t -\tau}( x) d\tau= \intR \eta^\epsilon_1(\tau) d\tau \intD u(x)d \sigma_{t- \tau}( x).$$
Hence, using that $v(y)= \intR \eta^\epsilon_1(\tau) v(y) d\tau,$ we obtain
\begin{align}
\intD u(x) d\sigma^\epsilon_t(x) + \intD v(y)d\sigma_t(y)
& =\intR  \eta^\epsilon_1(\tau) d\tau  \Bigl( \intD  u d \sigma_{t -\tau} + \intD vd\sigma_t
\Bigr)
 \nonumber\\
& \leq \intR  \eta^\epsilon_1(\tau) W_2^2(\sigma_{t -\tau}, \sigma_{t} )d\tau
\label{oct20.4b}\\
& \leq \intR  \eta^\epsilon_1(\tau) \tau^2 C_v^2 d\tau  = \epsilon^2 C_1 C_v^2.
  \label{oct20.4c}
\end{align}
To obtain Equation \ref{oct20.4b} we have used the dual formulation of the Wasserstein distance Equation \ref{eq:dualwasserstein} and the fact that $(u,v) \in {\mathcal{U}}$. We have used Remark  \ref{reac2curves} to obtain Equation \ref{oct20.4c}. Since $\intD u\,d\sigma^\epsilon_t + \intD v\,d\sigma_t \leq  \epsilon C C_v$ for arbitrary $(u,v) \in {\mathcal{U}}$, we conclude that (iii) holds. 

\endproof

\begin{remark}\label{reweakstarconv} Assume that for each $t \in (a,b)$ there exists $\rho_t>0$ such that $\sigma_t= \rho_t \calLD.$ Let $\phi \in C_c(\RD)$. Setting $ I_\phi(t):=\intD \langle \phi, v_t \rangle \rho_t d \calLD,$ we have
\begin{equation}\label{oct26.4}
|\intD \langle \phi, v^\epsilon_t\rangle \rho^\epsilon_td\calLD| =|\eta^\epsilon_1 \ast I_\phi(t) |
 \leq ||\phi||_\infty \; C_v.
\end{equation}
\end{remark}

\begin{corollary}\label{co10.26.1} Suppose that for each $t \in (a,b)$ there exists $\rho_t>0$ such that $\sigma_t= \rho_t \calLD.$ Then, for each $t \in [a,b]$, $\{\sigma^\epsilon_t\}_{\epsilon>0}$ converges to $\sigma_t$ in $\M$ as $\epsilon$ tends to zero. For $\calLone$-almost every $t \in [a,b]$, $\{\sigma^\epsilon_t v^\epsilon_t\}_{\epsilon>0}$ converges weak-$\ast$ to $\sigma_t v_t$ as $\epsilon$ tends to zero.
\end{corollary}

\proof{} By Lemma \ref{leapproximation1} (iii), $\{\sigma^\epsilon_t\}_{\epsilon>0}$ converges to $\sigma_t$ in $\M$ as $\epsilon$ tends to zero.

Let $\mathcal{C}$ be a countable family in $C_c(\RD)$. For each $\phi \in C_c(\RD)$, the set of Lebesgue points of $I_\phi$ is a set of full measure in $[a,b].$ For these points $ \eta^\epsilon_1 \ast I_\phi(t)$ tends to $ I_\phi(t)$ as $\epsilon$ tends to zero. Thus there is a set $S$ of full measure in $[a,b]$ such that for all $\phi \in \mathcal{C}$ and all $t \in S$, $ \eta^\epsilon_1 \ast I_\phi(t)$ tends to $ I_\phi(t)$ as $\epsilon$ tends to zero. Fix $\varphi \in C_c(\RD)$ and choose $\delta>0$ arbitrary. Let $\phi \in \mathcal{C}$ be such that $||\varphi- \phi||_\infty \leq \delta.$ Note that
$$|\eta^\epsilon_1 \ast I_\varphi(t) -I_\varphi(t)| \leq |\eta^\epsilon_1 \ast I_\phi(t) -  I_\phi(t)| + |\eta^\epsilon_1 \ast I_{\phi-\varphi}(t)| + | I_{\phi-\varphi}(t)|.$$
We use the inequality in Equation \ref{oct26.4} to conclude that
$$|\eta^\epsilon_1 \ast I_\varphi(t) -I_\varphi(t)| \leq |\eta^\epsilon_1 \ast I_\phi(t) -  I_\phi(t)| +2 \delta C_v.$$
If $t \in S$, the previous inequality gives that $\limsup_{\epsilon \rightarrow 0} |\eta^\epsilon_1 \ast I_\varphi(t) -I_\varphi(t)| \leq 2 \delta C_v.$ Since $\delta>0$ is arbitrary we conclude that $\lim_{\epsilon \rightarrow 0} |\eta^\epsilon_1 \ast I_\varphi(t) -I_\varphi(t)|=0.$ \endproof

\begin{corollary}\label{co:10.27.1} Suppose that $\sigma\in AC_2(a,b;\M)$ for all $a<b$, ${v}$ is a velocity associated to ${\sigma}$ and $C:=\sup_{t \in [a,b]} ||{ v}_t||_{\sigma_t}<\infty$. Define
$$
f^r_t(x):= \intD \eta_D^r(x-y) d \sigma_t(y), \quad
\sigma_t^r:= f^r_t \mathcal{L}^D, \quad
f^r_t(x) { v}^r_t(x):=\intD
\eta_D^r(x-y) { v}_t(y) d \sigma_t(y).
$$
As in Equation \ref{eqoct27.5}, for $0<\epsilon<1$ we define 
$$
\rho_t^{\epsilon, r} (x):= \intR \eta^\epsilon_1(t -\tau)f^r_\tau(x) d\tau, \;\;\;
\sigma^{\epsilon, r}_t:=
\rho^{\epsilon,r} \mathcal{L}^D, \;\;\; \rho_t^{\epsilon, r}(x)
v_t^{\epsilon,r}(x) := \intR \eta^\epsilon_1(t -\tau) f^r_\tau(x)
v^r_\tau(x) d\tau.
$$
Then:

(i)  ${v}^r$ is a velocity associated to $ \sigma^r$ and, for each $t \in (a,b)$, $\{ \sigma^r_t \}_{r}$ converges to
$ \sigma_t$ in $\M$ as $r$ tends to zero. For
$\mathcal{L}^1$-almost every $t \in (a,b),$  $||{ v}^r_t ||_{\sigma^r_t} \leq C$ and $\{{v}^r_{t} {\sigma}^r_t\}_{r>0}$ converges weak-$\ast$ to ${v}_t {\sigma}_t$ as $r$ tends to zero.

(ii)  ${ v}^{\epsilon, r}$ is a velocity associated to $
\sigma^{\epsilon, r}$ and, for each $t \in (a,b)$, $\{
\sigma^{\epsilon, r}_t \}_{\epsilon}$ converges to $
\sigma^r_t$ in $\M$ as $\epsilon$ tends to zero. For every $t \in
(a,b),$  $||{ v}^{\epsilon, r}_t ||_{ \sigma^{\epsilon,
r}_t} \leq C$ while for $\mathcal{L}^1$-almost every $t \in (a,b),$
$\{   { v}^{\epsilon, r}_{t} { \sigma}^{\epsilon,
r}_t\}_{r>0}$ converges weak-$\ast$ to ${ v}^r_t {
\sigma}_t^r$ as $\epsilon$ tends to zero.

(iii) The function $t \rightarrow \bar \Lambda_{\sigma^{\epsilon,
r}_t}(v_t^{\epsilon,r})$ is continuous  while $t \rightarrow \bar
\Lambda_{ \sigma_t}({ v}_t)$ is measurable on $(a,b).$

(iv) Suppose in addition that $\sigma$ is time-periodic in the sense that $\sigma_a=\sigma_b$. Then $\sigma^{r}_a=\sigma^{r}_b$.
\end{corollary}
\proof{} It is well known that $||{ v}^r_t ||_{\sigma^r_t} \leq ||{ v}_t||_{ \sigma_t} \leq C$ (cf. \cite{ags:book} Lemma 8.1.10) so, by Remark \ref{reac2curves} (i), $\sigma \in AC_2(a,b;\M).$ One can readily check that ${v}^r$ is a velocity associated to $\sigma^r.$ Lemma \ref{leapproximation1} shows that, for each $t \in (a,b)$, $\{\sigma^r_t \}_{r}$ converges to $\sigma_t$ in $\M$ as $r$
tends to zero. Let $\varphi \in C_c(\RD, \RD).$ Set $\varphi^r:=\eta_D^r \ast \varphi.$ Since $\{\varphi^r\}_{r>0}$ converges
uniformly to $\varphi$,
$$
\lim_{r \rightarrow 0} \intD \langle \varphi, { v}^r_t\rangle
d{ \sigma}^r_t= \intD \langle { v}_t, \varphi \rangle
d \sigma_t.
$$
Thus $\{{ v}^r_{t} {\sigma}^r_t\}_{r>0}$ converges weak-$\ast$ to ${v}_t { \sigma}_t$ as $r$ tends to zero. This proves (i).

We next fix $r>0.$ For a moment we won't display the dependence in $r.$ For instance we write $v^\epsilon$ instead of
$v_t^{\epsilon,r}$ as in Equation \ref{eqoct27.5}. Notice that
$\rho^\epsilon \in C^1([a,b] \times \RD),$ $\rho^\epsilon>0$ and $\rho^\epsilon_t$ is a probability density. Also $v_t^\epsilon \in C^1([a,b] \times \RD, \RD)$ and $v^\epsilon$ is a velocity associated to $\sigma^\epsilon.$ Fix $t \in [\bar a, \bar b]\subset (a,b).$ Lemma \ref{leapproximation1} gives that $||v^\epsilon_t||_{\sigma^\epsilon_t} \leq C$ for all $\epsilon>0$ small enough. By Corollary \ref{co10.26.1}, $\{v^\epsilon_{t}\sigma^\epsilon_t\}_{\epsilon>0}$ converges weak-$\ast$ to $v_t\sigma_t$ as $\epsilon$ tends to zero. This proves (ii).

By Lemma \ref{le10.27.7}, $t \rightarrow \bar\Lambda_{\sigma^\epsilon_t}(v^\epsilon_t)$ is continuous in $(a,b).$ Hence by (ii) $t \rightarrow \bar \Lambda_{{\sigma}^{r}_t}({ v}_t^r)$ is measurable as a pointwise limit of measurable functions. We then use (i) to conclude that $t\rightarrow \bar \Lambda_{{ \sigma}_t}({ v}_t)$ is measurable as a pointwise limit of measurable functions.  This proves (iii). The proof of (iv) is straightforward. 

\endproof

We can now prove that regular pseudo 1-forms can be integrated along absolutely continuous curves, as follows.

\begin{corollary}\label{cointegrability} Let $\sigma \in AC_2(a,b;\M)$ and let $v$ be a velocity associated to $\sigma.$ Suppose $t \rightarrow ||v_t||_{\sigma_t}$ is square integrable on $(a,b).$ Then $t \rightarrow \bar \Lambda_{\sigma_t}(v_t)$ is measurable and square integrable on $(a,b).$
\end{corollary}

\proof{} Let $\bar \sigma$ be the reparametrization of $\sigma$ as introduced in Remark \ref{reparametrization1} and let $\bar v$ be the associated velocity. By Corollary \ref{co:10.27.1} (iii), because $\sup_{s \in [0,L]}||\bar v_s||_{\bar \sigma_s} \leq 1,$ we have that
$ s \rightarrow \bar \Lambda_{\bar \sigma_{s}} (\bar v_{s})$ is measurable. But $\bar \Lambda_{ \sigma_{t}}(v_{t})=\dot S(t) \bar \Lambda_{\bar \sigma_{S(t)}} (\bar v_{S(t)}) .$ Thus $t \rightarrow \bar \Lambda_{ \sigma_{t}}(v_{t})$ is measurable.

By Corollary \ref{le:boundedA} there exists a constant $C_\sigma$ independent of $t$ such that $||\bar A_{\sigma_t}||_{\sigma_t} \leq C_\sigma$ for all $t \in [a,b].$ Thus
$$
|\bar \Lambda_{\sigma_t}(v_t)|=\Bigl| \intD\langle \bar A_{\sigma_t}, v_t\rangle d\sigma_t \Bigr| \leq ||\bar A_{\sigma_t}||_{\sigma_t} ||v_t||_{\sigma_t} \leq  C_\sigma ||v_t||_{\sigma_t}.
$$
Since $t \rightarrow ||v_t||_{\sigma_t}$ is square integrable, the previous inequality yields the proof. \endproof

\begin{corollary}\label{co:convergenceform} Suppose $\{\sigma^r\}_{0 \leq r \leq c} \subset AC_2(a,b;\M),$ $ v^r$ is a velocity associated to $\sigma^r$ and $\infty> C:=\sup_{(t,r) \in E} ||v^r_t||_{\sigma^r_t}$ where $E:=[a,b] \times [0,c].$ Suppose that, for $\calLone$-almost every $t \in (a,b),$ $\{   { v}^r_{t} { \sigma}^r_t\}_{r>0}$ converges weak-$\ast$ to ${ v}_t { \sigma}_t$ and $\{ \sigma^r_t\}_{r>0}$ converges in $\M$ to $\sigma_t$ as $r$ tends to zero. If $(t,r) \rightarrow \sigma^r_t$ is continuous at every $(t,0) \in [a,b] \times \{0\}$ then $\lim_{r \rightarrow 0} \int_a^b \bar \Lambda_{\sigma^r}(v^r) dt= \int_a^b \bar \Lambda_\sigma(v) dt.$ Here we have set $\sigma_t:=\sigma_t^0.$
\end{corollary}

\proof{} By Lemma \ref{le:boundedA} we may assume without loss of generality that $||\bar A_{\sigma^r_t}||_{\sigma^r_t}$ is bounded on $E$ by a constant $\bar C$ independent of $(t,r) \in E.$
We obtain
\begin{equation}\label{eq:boundsigmart}
\sup_{(t,r) \in E} |\bar \Lambda_{\sigma^r_t} (v^r_t)| \leq \sup_{(t,r) \in E} ||\bar A_{\sigma^r_t}||_{\sigma^r_t} ||v^r_t||_{\sigma^r_t} \leq  \bar C C.
\end{equation}

Corollary \ref{co10.12.3} ensures that $\lim_{r \rightarrow 0} \bar \Lambda_{\sigma^r_t} (v^r_t)=\bar \Lambda_{\sigma_t} (v_t)$ for $\calLone$-almost every $t \in [a,b].$ This, together with Equation \ref{eq:boundsigmart} shows that, as $r$ tends to $0$, the sequence of functions $t \rightarrow \bar \Lambda_{\sigma^r_t} (v^r_t)$ converges to the function $t \rightarrow \bar \Lambda_{\sigma_t} (v_t)$ in $L^1(a,b).$ This proves the corollary. \endproof

\begin{definition}\label{de10.15.1} 
Let $\sigma \in AC_2(a,b; \M)$ and let $v$ be a velocity associated to $\sigma.$ Suppose $t \rightarrow ||v_t||_{\sigma_t}$ is square integrable on $(a,b).$ By Corollary \ref{cointegrability},  $t \rightarrow \bar \Lambda_{\sigma_t}(v_t)$ is also square integrable on $(a,b).$  It is thus meaningful to calculate the integral $ \int_a^b \bar \Lambda_{\sigma_t}(v_t)dt$. 

We will call $\int_a^b \bar \Lambda_{\sigma_t}(v_t)dt$ the \textit{integral of $\bar \Lambda$ along $(\sigma,v)$}. 
When $v$ is the velocity of minimal norm we will write this simply as $\int_\sigma\bar\Lambda$ and call it the \textit{integral of $\bar \Lambda$ along $\sigma$}.

\end{definition}

\begin{remark}\label{re:integralcurve1} Suppose that $ r:[c,d] \rightarrow [a,b]$ is invertible and Lipschitz. Define $\bar \sigma_s=\sigma_{r(s)}$. Then $\bar \sigma \in AC_2(c,d; \M)$ and $ \bar v_s(x)= \dot r(s) v_{r(s)}(x)$ is a velocity for $\bar \sigma.$ Furthermore, $ \int_c^d \bar \Lambda_{\bar \sigma_t} (\bar v_t)dt = \int_a^b \bar \Lambda_{\sigma_t}(v_t)dt.$
\end{remark}

\proof{} Let $\beta \in L^2(a,b)$ be as in Definition \ref{deac2curves}. Then
$$
W_2(\sigma_{r(s+h)}, \sigma_{r(s)}) \leq \Bigl | \int_{r(s)}^{r(s+h)} \beta(t) dt \Bigr | = \Bigl | \int_s^{s+h} \bar \beta(\tau) d\tau \Bigr | \quad \hbox{where}  \quad \bar \beta  (s):= |\dot r(s) | \beta(r(s)).
$$
Because $\bar \beta \in L^2(c,d)$ we conclude that $\bar \sigma \in AC_2(c,d; \M)$. Direct computations give that, for ${\mathcal{L}}^1-\hbox{a.e.} \;\; s \in (c,d)$,
$$\lim_{h \rightarrow 0} W_2(\sigma_{r(s+h)}, \sigma_{r(s)}) / |h|= |\dot r(s) | \, |\sigma'|(r(s)).
$$
Thus $|\bar \sigma'|(s)= |\dot r(s) | \; |\sigma'|(r(s)).$ Let $\phi \in C_c^\infty(\RD)$ and let $v$ be a velocity for $\sigma$ (see Proposition \ref{prac2curves}). The chain rule shows that, in the sense of distributions,
$$
 {d \over ds}  \int_{\RD} \phi d \sigma_{r(s)}= \dot r(s) \langle \nabla \phi, v_{r(s)} \rangle_{\sigma_{r(s)}}= \langle \nabla \phi, \bar v_{s} \rangle_{\bar \sigma_s},$$
where $\bar v_s(x)= \dot r(s) v_{r(s)}(x).$ Thus $\bar v$ is a velocity for $\bar \sigma.$ Using the linearity of $\bar \Lambda$ we have
$$
\int_c^d \bar \Lambda_{\bar \sigma_s}(\bar v_s) ds = \int_c^d \dot r(s)  \bar \Lambda_{ \sigma_{r(s)}}(v_{r(s)}) ds=\int_a^b \bar \Lambda_{ \sigma_{t}} (v_t) dt.
$$
This concludes the proof. \endproof


\subsection{Green's formula for annuli}\label{ss:greenstheorem}
Let $\sigma \in AC_2(a,b; \M)$. Given $r\in (0,1)$ and $s\in [r,1]$, let $D_s:\R^D\rightarrow \R^D$ denote the map defined by $D_s(x):=sx$. Using this map we can canonically associate to $\sigma$ the surface 
\begin{equation}
\sigma(s,t):[r,1]\times(a,b)\rightarrow \M, \ \ \sigma(s,t)=\sigma_t^s:=D_{s \#} \sigma_t.
\end{equation}
We call such a surface an \textit{annulus}. We now want to study its properties. 

Let $v$ denote the velocity of $\sigma$ of minimal norm. Set
$$
w(s,t, \cdot)=w_t^s(x):={x \over s}=D^{-1}_s(x), \ \ v(s,t, \cdot)=v_t^s:=D_{s*} v_t.
$$
According to Lemma \ref{le:pushforwardvelocity}, for each $s  \in  [r,1]$, $\sigma(s, \cdot) \in AC_2(a,b; \M)$ admits $v(s, \cdot)$ as a velocity. For each $t$ and $\phi \in C_c^\infty(\RD)$, in the sense of distributions,
$$
{d \over ds} \intD \phi\,d\sigma_t^s= {d \over ds} \intD \phi (sx)\,d\sigma_t(x)= \intD d\phi (sx)(x)\,d\sigma_t(x)=\intD d\phi(w_t^s)\,d\sigma^s_t.
$$
Thus $w(\cdot, t)$ is a velocity for $\sigma(\cdot, t).$

We assume that
$$
||\sigma'||_\infty:= \sup_{t \in [a,b]} ||v_t||_{\sigma_t}<\infty.
$$
By Remark \ref{reac2curves},
$$
c^0_\sigma:=\sup_{t \in [a,b]} W_2(\sigma_t, \delta_0) <\infty.
$$
By the fact that $D_{s \#} \sigma_t= \sigma_t^s$ we have
\begin{equation} \label{eq:greenannulus6}
W^2_2(\sigma^s_t, \delta_0)= s^2 W^2_2(\sigma_t, \delta_0) \leq s^2 c^0_\sigma \leq \bar C_\sigma,
\end{equation}
where we are free to choose $\bar C_\sigma$ to be any constant greater than $c^0_\sigma.$

\begin{remark} \label{re:greenannulus} Note that $(1+ h/s)\id$ pushes $\sigma_t^s$ forward to $\sigma_t^{s+h}$ and is the gradient of a convex function.  Thus
$$
\gamma^{h}:= \Bigl(\id \times (1+ h/s)\id  \Bigr)_\# \sigma_t^s  \in \Gamma_o(\sigma_t^s , \sigma_t^{s+h}).
$$
For $\gamma^{h}$-almost every $(x,y) \in \RDD$ we have $y= (1 +h/s) x$, so
\begin{equation} \label{eq:greenannulus1}
v_t^{s+h}(y)=(s+h) v_t({y \over s+h})= (1 +{h \over s}) v^s_t({sy \over s+h})= (1 +{h \over s}) v_t^s(x).
\end{equation}
Using the definition of $\sigma_t^s$ and $v_t^s$ we obtain the identities
\begin{equation} \label{eq:greenannulus3}
||\id||_{\sigma_t^s}=s ||\id||_{\sigma_t} \leq s \bar C_\sigma, \quad ||v_t^s||_{\sigma_t^s}= s ||v_t||_{\sigma_t} \leq s ||\sigma'||_\infty.
\end{equation}
We use the first identity in Equation \ref{eq:greenannulus3} and the fact that $(1+ h/s)\id$ pushes $\sigma_t^s$ forward to $\sigma_t^{s+h}$ to obtain
\begin{equation} \label{eq:greenannulus4}
W^2_2(\sigma_t^s, \sigma_t^{s+h})={h^2 \over s^2} ||\id||^2_{\sigma_t^s}=h^2 ||\id||^2_{\sigma_t}= h^2 W_2^2(\sigma_t, \delta_0) \leq h^2 \bar C^2_\sigma .
\end{equation}
\end{remark}

\begin{lemma}\label{coboundedA} 
There exists a constant $C_\sigma(r)$ depending only on $\sigma$ and $r$ such that $||\bar A_{\sigma^s_t}||_{\sigma^s_t} \leq C_\sigma(r)$, for all $(s, t) \in [r, 1] \times [a,b].$
\end{lemma}

\proof{} By Remark \ref{reac2curves} (i), $\sigma: [a,b] \rightarrow \M$ is $1/2$-H\"older continuous: there exists a constant $c>0$ such that $W^2_2(\sigma_{t_2}, \sigma_{t_1}) \leq c |t_2- t_1|. $ Together with Lemma \ref{l:dualwasserstein2} and the fact that $Lip(D_s)=s \leq 1$, this gives that $t \rightarrow \sigma_t^s$ is uniformly $1/2$-H\"older continuous:
$$
W^2_2(\sigma_{t_2}^s, \sigma_{t_1}^s) \leq W^2_2(\sigma_{t_2}, \sigma_{t_1}) \leq  c |t_2- t_1|.
$$
Remark \ref{reac2curves} (ii) ensures that $\{ \sigma_t | \; t \in [a,b]\}$ is bounded and so there exists $\bar c>0$ such that $W_2(\sigma_t, \delta_0) \leq \bar c$ for all $t \in [a,b].$
One can readily check that $\gamma:=\bigl( D_{s_1}  \times D_{s_2} \bigr)_\# \sigma_t \in \Gamma(\sigma_{t}^{s_1}, \sigma_{t}^{s_2})$, so
\begin{align*}
W^2_2(\sigma_{t}^{s_1}, \sigma_{t}^{s_2}) &\leq \intDD |x-y|^2 d \gamma= \intD |D_{s_1}  x - D_{s_2} x|^2 d\sigma_t(x)= |s_2 -s_1|^2 \intD |x|^2 d\sigma_t(x)\\  
&\leq \bar c |s_2 -s_1|^2 .
\end{align*}
Thus $s \rightarrow \sigma_t^s$ is 1-Lipschitz. Consequently $(t,s) \rightarrow \sigma_t^s$ is 1/2-H\"older continuous. This, together with Lemma \ref{le:boundedA},  yields the proof.
\endproof

Set
$$
V(s,t):= \bar \Lambda_{\sigma_t^s}(v_t^s), \qquad W(s,t):= \bar \Lambda_{\sigma_t^s}(w_t^s).
$$
The following proposition is extracted from \cite{ags:book} Theorem 8.3.1 and Proposition 8.4.6.

\begin{prop} \label{pr:geodesic} Let $\sigma\in AC_2(a,b;\M)$ and let $v$ be its velocity of minimal norm. Let $\mathcal{N}_1$ be the set of $t$ such that $v_t$ fails to be in $T_{\sigma_t} \M$. Let $\mathcal{N}_2$ be the set of $t \in [a,b]$ such that $\Bigl( \pi^1 \times (\pi^2 -\pi^1 )/h     \Bigr)_\# \eta_h$ fails to converge to $(\id \times v_t)_{\#}{\sigma_t}$ in the Wasserstein space $\M(\RD \times \RD)$, for some $\eta_h \in \Gamma_o(\sigma_{t}, \sigma_{t+h}).$ Let $\mathcal{N}$ be the union of $\mathcal{N}_1$ and $\mathcal{N}_2.$ Then $\calLone(\mathcal{N})=0.$
\end{prop}

We can now study the derivatives of $V$ and $W$, as follows.

\begin{lemma} \label{le:lipschitzgreen1} For each $t \in (a,b) \setminus \mathcal{N}$, the function $V(\cdot,t)$ is differentiable on $(r,1)$ and its derivative is bounded by  a constant $L_1(r)$ depending only on $\sigma$ and $r.$ Furthermore
$$
\partial_sV(s, t)=\intD \langle \bar A_{\sigma_t^{s}}(x) ,{v_t^s (x) \over s} \rangle d\sigma_t^{s} (x) +
\intD \langle B_{\sigma_t^{s}}(x)w_t^s(x) , v_t^{s} (x)\rangle d\sigma_t^{s} (x).
$$
\end{lemma}
\proof{} Let $C_\sigma(r)$ be as in Lemma \ref{coboundedA} and let $\bar C_\sigma$ be as in Equation \ref{eq:greenannulus6}. We use Equations \ref{eq:contdiff3}, \ref{eq:greenannulus1} and then H\"older's inequality to obtain
\begin{equation} \label{eq:greenannulus8}
|V(s+h,t) - V(s,t)|  \leq
{h \over s}  ||\bar A_{\sigma_t^s}||_{\sigma_t^s} \, || v_t^s|| _{\sigma_t^s} +
2 c(\bar\Lambda)W_2(\sigma_t^s, \sigma_t^{s+h}) ||v_t^{s+h}||_{\sigma_t^{s+h}} .
\end{equation}
We combine Equations \ref{eq:greenannulus3}, \ref{eq:greenannulus4} and \ref{eq:greenannulus8} to conclude that
\begin{equation} \label{eq:greenannulus10}
|V(s+h,t) - V(s,t)| \leq
h C_\sigma(r)  ||\sigma'||_\infty + 2 h c(\bar\Lambda) \bar C_\sigma (s+ h) ||\sigma'||_\infty .
\end{equation}
This proves that $V(\cdot, t)$ is Lipschitz on $(r,1)$ and that its derivative is bounded by  a constant $L_1(r).$ As in Remark \ref{re:contdiffbis},
\begin{align}
\lim_{h \rightarrow 0} {V(s+h,t) - V(s,t) \over h}
&=\lim_{h \rightarrow 0}
 \intDD \langle \bar A_{\sigma_t^{s}}(x) ,{ v_t^{s+h} (y)-v_t^s (x) \over h} \rangle d\gamma^{h}\nonumber \\
&\quad +
 \intDD \langle B_{\sigma_t^{s}}(x){ y-x \over h} , v_t^{s+h} (y)\rangle d\gamma^{h}\nonumber \\
&\quad + {1\over h} \intDD \langle \bar A_{\sigma_t^{s+h}}(y) -\bar A_{\sigma_t^{s}}(x)-B_{\sigma_t^{s}}(x)(y-x) , v_t^{s+h} (y)\rangle d\gamma^{h}, \label{eq:greenannulus11}
\end{align}
where $\gamma_h\in\Gamma_o(\sigma_t^{s+h},\sigma_t^s)$. By Equation \ref{eq:contdifflate1}, the last inequality in Equation \ref{eq:greenannulus3} and Equation \ref{eq:greenannulus4} we have
\begin{equation} \label{eq:greenannulus12}
\lim_{h \rightarrow 0} {1\over h} \intDD \langle \bar A_{\sigma_t^{s+h}}(y) -\bar A_{\sigma_t^{s}}(x)-B_{\sigma_t^{s}}(x)(y-x) , v_t^{s+h} (y)\rangle \gamma^{h} (x,y)=0.
\end{equation}
We use  Equations \ref{eq:greenannulus1}, \ref{eq:greenannulus11},  \ref{eq:greenannulus12}  and the fact that, for $\gamma^h$-almost every $(x,y) \in \RDD$, $y= (1 +h/s) x$ to conclude that
\begin{align*}
\lim_{h \rightarrow 0} {V(s+h,t) - V(s,t)  \over h}
&=\intD \langle \bar A_{\sigma_t^{s}}(x) ,{v_t^s (x) \over s} \rangle d\sigma_t^{s}
+ \lim_{h \rightarrow 0}
 \intD \langle B_{\sigma_t^{s}}(x){ x \over s} , (1+{h \over s}) v_t^{s} (x)\rangle d\sigma_t^{s}
\nonumber \\
&=
\intD \langle \bar A_{\sigma_t^{s}}(x) ,{v_t^s (x) \over s} \rangle d\sigma_t^{s} +
\intD \langle B_{\sigma_t^{s}}(x)w_t^s(x) , v_t^{s} (x)\rangle d\sigma_t^{s}
. \label{eq:greenannulus13}
\end{align*}
This proves the lemma.  \endproof

\begin{lemma} \label{le:lipschitzgreen2} For each $s \in [r,1]$ and $t \in (a,b) \setminus \mathcal{N},$ the function $W(s, \cdot)$ is differentiable at $t$ and its derivative is bounded by  a constant $L_2(r)$ depending only on $\sigma$ and $r.$ Furthermore,
$$
\partial_t W(s, t)=\intD \langle \bar A_{\sigma_t^{s}}(x) ,{v_t^s (x) \over s} \rangle d\sigma_t^{s} (x) +
\intD \langle w_t^s(x) , B_{\sigma_t^{s}}(x) v_t^{s} (x)\rangle d\sigma_t^{s} (x).
$$
\end{lemma}

\proof{} As in the proof of Lemma 5.30, we have
\begin{equation*}
|W(s,t+h) - W(s,t)|  \leq {||\bar
A_{\sigma_t^s}||_{\sigma_t^s} \over s}  \, W_2(\sigma_{t+h}^s, \sigma_t^s) + 2
c(\bar\Lambda)W_2(\sigma_{t+h}^s, \sigma_t^s)
||w_{t+h}^s||_{\sigma_{t+h}^s} .
\end{equation*}
This gives
\begin{equation*}
|W(s,t+h) - W(s,t)| \leq
h \bigl (1+2c(\bar\Lambda)\bigr) \bar C_\sigma,
\end{equation*}
proving that $W(s, \cdot)$ is Lipschitz on $(a,b)$ and that its derivative is bounded by a constant $L_2(r).$ For fixed $s \in [r,1]$ and $t \in (a,b) \setminus \mathcal{N}$, let $\gamma_h^s \in\Gamma_o(\sigma_t^s,\sigma_{t+h}^s)$. Then
\begin{align*}
\lim_{h \rightarrow 0} &{W(s,t+h) - W(s,t) \over h}\\
&=\lim_{h
\rightarrow 0}
 \intDD \langle \bar A_{\sigma_t^s}(x) ,{w_{t+h}^s (y)-w_t^s (x) \over h} \rangle d\gamma_h^s (x,y)\\
&\quad +
 \intDD \langle B_{\sigma_t^{s}}(x){ y-x \over h} , w_{t+h}^{s} (y)\rangle d\gamma_h^s (x,y)\\
&\quad + {1\over h} \intDD \langle \bar A_{\sigma_{t+h}^{s}}(y) -\bar
A_{\sigma_t^{s}}(x)-B_{\sigma_t^{s}}(x)(y-x) , w_{t+h}^{s}
(y)\rangle \gamma_h^s (x,y).
\end{align*}
However,
\begin{equation*}
\lim_{h \rightarrow 0} {1\over h} \intDD \langle \bar
A_{\sigma_t^{s}}(y) -\bar
A_{\sigma_t^{s}}(x)-B_{\sigma_t^{s}}(x)(y-x) , w_{t+h}^{s}
(y)\rangle \gamma_h^s (x,y)=0.
\end{equation*}
We then use the fact that $w_t^s(z)=z/s$ to get
\begin{align}\label{eq1}
\lim_{h \rightarrow 0} {W(s,t+h) - W(s,t) \over h}
&= \lim_{h\rightarrow 0} \intDD \langle \bar A_{\sigma_t^s}(x) ,{y-x \over sh}
\rangle d\gamma_h^s (x,y) \\\nonumber
&\quad +\intDD \langle B_{\sigma_t^{s}}(x){ y-x \over h} , {y\over s}\rangle d\gamma_h^s(x,y).
\end{align}
To conclude the lemma it suffices to show that if  $t \in (a,b) \setminus \mathcal{N}$ and $\gamma_h^s \in \Gamma_o(\sigma_t^s, \sigma_{t+h}^s)$ then $\Bigl( \pi^1 \times (\pi^2 -\pi^1 )/h     \Bigr)_\# \gamma_h^s$ converges to $(\id \times v^s_t)_{\sigma^s_t}$ in $\M(\RD \times \RD)$ as $h$ tends to $0.$ Set
$$
\gamma_h:= \bigl(  D_s^{-1} \times D_s^{-1} \bigr)_\#  \gamma_h^s.
$$
Since
$$
\pi^1 \circ \bigl(  D_s^{-1} \times D_s^{-1} \bigr)= D_s^{-1}
\circ \pi^1 \quad \hbox{and} \quad \pi^2 \circ \bigl(  D_s^{-1}
\times D_s^{-1} \bigr)= D_s^{-1} \circ \pi^2,
$$
we conclude that $\gamma_h \in \Gamma(\sigma_t, \sigma_{t+h}).$ By the fact that the support of $\gamma_h^s$ is cyclically monotone
we have that the support of $\gamma_h$ is also cyclically monotone. Hence $\gamma_h \in \Gamma_o(\sigma_t, \sigma_{t+h}).$
We have
\begin{align*}
\bigl( \pi^1 \times {\pi^2 -\pi^1 \over h}     \bigr)_\#
\gamma_h^s &= \bigl(  D_s \times D_s \bigr)_\# \Bigl( ( \pi^1 \times {\pi^2
-\pi^1 \over h} )_\# \gamma_h \Bigr)\rightarrow\bigl(  D_s
\times D_s \bigr) \circ (\id \times v_t)_\# \sigma_t\\
&=(\id \times v^s_t)_\# \sigma^s_t.
\end{align*}
This, together with Equation \ref{eq1}, gives
\begin{align*}
\lim_{h \rightarrow 0} &{W(s,t+h) - W(s,t) \over h}\\
&= \intD \langle \bar A_{\sigma_t^s}(x) ,{v_t^s(x) \over s}\rangle d\sigma_t^s (x)
+\intD \langle B_{\sigma_t^{s}}(x)v_t^s(x) , {x\over s}\rangle d\sigma_t^s(x)\\
&= \intD \langle \bar A_{\sigma_t^s}(x) ,{v_t^s(x) \over s}
\rangle d\sigma_t^s (x)
+\intD \langle B_{\sigma_t^{s}}(x)v_t^s(x) , w_t^s(x)\rangle d\sigma_t^s(x).
\end{align*}
\endproof


\begin{corollary}\label{le:lipschitzgreen3} For each $s \in (r,1)$ and $t \in (a,b) \setminus \mathcal{N}$ we have
$$
\partial_t \Bigl( \bar \Lambda_{\sigma_t^s}(w_t^s) \Bigr)- \partial_s \Bigl(  \bar \Lambda_{\sigma_t^s}(v_t^s) \Bigr)= d\bar \Lambda_{\sigma_t^s}(v_t^s, w_t^s).
$$
\end{corollary}

\proof{} This corollary is a direct consequence of Lemmas \ref{le10.12.4}, \ref{le:lipschitzgreen1} and \ref{le:lipschitzgreen2}. \endproof

\begin{remark} 
Proposition \ref{surfacegreen1} was a direct consequence of Equations \ref{eq:conclusiongreen1} and \ref{eq:conclusiongreen2}. Those equations depended strongly on the smoothness of $v^s_t$ and $w^s_t$ in all variables. In this section we have removed all smoothness assumptions on $v^s_t$. Specifically, now we do not know that $\nabla v^s_t$ nor $\partial_t v^s_t$ exist. However, Equation \ref{eq:greenannulus1} ensures that $v^s_t$ is differentiable with respect to $s$, in particular establishing the inequality $||v_t^{s+h} \circ \pi^2 -v_t^s \circ \pi^1||_{\gamma^h} \leq h ||\sigma'||_\infty.$ This inequality is crucial for the proof of Lemma \ref{le:lipschitzgreen1}. 

Lemma \ref{le:lipschitzgreen1} is the analogue of Equation \ref{eq:conclusiongreen2}. Notice that the first term (respectively, the second term) on the right hand side of Lemma \ref{le:lipschitzgreen1} corresponds to the second term (respectively, the first term) on the right hand side of Equation \ref{eq:conclusiongreen2}. Likewise, Lemma \ref{le:lipschitzgreen2} is the analogue of Equation \ref{eq:conclusiongreen1}.  
\end{remark}

\begin{theorem}[Green's formula on the annulus] \label{th:greencurveannulus} Consider in $\M$ the annulus $S(s,t)=D_{s \# } \sigma_t$ for $(s,t ) \in [r, 1 ] \times [0,T]$. Let $\partial S$ denote its boundary, defined as the union of the negatively oriented curves $S(r, \cdot),$ $S(\cdot,T)$ and the positively oriented curves $S(1, \cdot),$ $S(\cdot,0).$  Then
$$
\int_{S} d\bar \Lambda= \int_{\partial S} \bar \Lambda.
$$
\end{theorem}
\proof{} We use Corollary \ref{le:lipschitzgreen3} to obtain
\begin{align}
\int_{S} d\bar \Lambda&= \int_0^T dt \int_r^1 d\bar \Lambda_{S(s,t)}(v_t^s,w_t^s)ds =
 \int_0^T dt \int_r^1 \Bigl[ \partial_t \Bigl( \bar \Lambda_{S(s,t)}(w_t^s) \Bigr)- \partial_s \Bigl(  \bar \Lambda_{S(s,t)}(v_t^s) \Bigr) \Bigr] ds
\nonumber \\
&= \int_r^1 \Bigl( \bar \Lambda_{S(s,T)}(w_T^s) -\bar \Lambda_{S(s,0)}(w_0^s) \Bigr) ds -
\int_0^T  \Bigl(  \bar \Lambda_{S(1,t)}(v_t^1) -\bar \Lambda_{S(r,t)}(v_t^r) \Bigr)dt
= \int_{\partial S} \bar \Lambda.
\label{eq:greencurveannulus1}
\end{align}  \endproof

\begin{corollary}\label{co:closecurve} If we further assume that $\bar \Lambda$ is a closed pseudo $1$-form and that $\sigma_0=\sigma_T$, then $\int_\sigma\bar \Lambda=0.$
\end{corollary}
\proof{} For $s \in [r,1]$ define
$$
l(s)=\int_0^T \bar \Lambda_{S(s,t)}(v_t^s)  dt, \qquad \bar l(t)=\int_r^1  \bar \Lambda_{S(s,t)}(w_t^s)ds.
$$
Since $w_T^s=w_0^s$ and $\sigma_T^s=D_{s \#} \sigma_T=D_{s \#} \sigma_0=\sigma_0^s$, we have  $\bar l(T)=\bar l(0)$. This, together with Equation \ref{eq:greencurveannulus1} and the fact that $d\bar \Lambda=0$, yields $\int_0^T \bar \Lambda_{\sigma_t}(v_t)dt=l(1)=l(r).$  But
\begin{equation} \label{eq:closecurve2}
|l(r)|\leq \int_0^T | \bar \Lambda_{S(s,t)}(v_t^r)| dt \leq \int_0^T || \bar A_{S(s,t)}||_{S(s,t)} ||v_t^r||_{S(s,t)} dt \leq r  ||\sigma'||_\infty
 \int_0^T || \bar A_{S(s,t)}||_{S(s,t)} dt,
\end{equation}
where we have used the last inequality in Equation \ref{eq:greenannulus3}. The first inequality in Equation \ref{eq:greenannulus6} shows that, for $r$ small enough, $\{S(s,t) \}_{t \in \times [0,T])}$ is contained in a small ball centered at $\delta_0$. But Lemma \ref{le:boundedA} gives that $\mu \rightarrow ||\bar A_\mu||_\mu$ is continuous at $\delta_0.$ Thus there exist constants $c$ and $r_0$ such that   $ || \bar A_{S(s,t)}||_{S(s,t)} \leq c$ for all $t \in [0,T]$ and all $r<r_0.$  We can now exploit Equation \ref{eq:closecurve2} to obtain
$$
|l(1)| =\liminf_{r \rightarrow 0} |l(r)| \leq \liminf_{r \rightarrow 0} rT c ||\sigma'||_\infty=0.
$$
\endproof

\begin{corollary}\label{co:closedexact} Let $\bar \Lambda$ be a regular pseudo $1$-form on $\M$. Let $\Lambda$ denote the corresponding $1$-form on $\M$, defined by restriction. Assume $\bar\Lambda$ is closed, \textit{i.e.} $d\bar\Lambda=0$. Then $\Lambda$ is exact, \textit{i.e.} there exists a differentiable function $F$ on $\M$ such that $d F=\Lambda.$
\end{corollary}
\proof{}
Fix $\mu \in \mathcal{M}.$  Let $\sigma$ be any curve in $AC_2(a,b; \mathcal{M})$ such that $\sigma_a=\delta_0$ and $\sigma_b=\mu.$ Assume that $v$ is its velocity of minimal norm and that $\sup_{(a,b)} ||v_t||_{\sigma_t}<\infty.$ By Corollary \ref{co:closecurve}, $\int_{\sigma} \bar \Lambda$ depends only on $\mu$, \textit{i.e.} it is independent of the path $\sigma.$ Also, Remark \ref{re:integralcurve1} ensures that $\int_{\sigma} \bar \Lambda$ is independent of $a, b.$ It is thus meaningful to define
$$F(\mu):= \int_{\sigma} \bar \Lambda.$$
We now want to show that $F$ is differentiable. Fix $\mu, \nu \in \mathcal{M}$ and $\gamma \in \Gamma_o(\mu,\nu).$ Define $\sigma_t := ((1-t)\pi^1 + t \pi^2)_{\#}\gamma$. Then $\sigma:[0,1] \rightarrow \mathcal{M}$ is a constant speed geodesic between $\mu$ and $\nu$. Let $v_t$ denote its velocity of minimal norm. Clearly,
\begin{equation}\label{eq:Fdifference}
F(\nu) - F(\mu)= \int_0^1 \bar \Lambda_{\sigma_t}(v_t)dt.
\end{equation}
Let $\bar\gamma:\RD\rightarrow\RD$ denote the \textit{barycentric projection} of $\gamma$, cf. \cite{ags:book} Definition 5.4.2. Set $v:= \bar{\gamma} - Id$. Then $\gamma_t:=(\pi^1,(1-t)\pi^1 +t\pi^2)_{\#}\gamma \in \Gamma_o(\sigma_0,\sigma_t)$ and
\begin{align*}
\bar \Lambda_{\sigma_t}(v_t) - \bar \Lambda_{\sigma_0}(v)&= \int_{\RDD} \langle  \bar{A}_{\sigma_0}(x), v_t(y)-v(x)\rangle +
\langle B_{\sigma(0)}(x)(y-x),v_t(y)\rangle d\gamma_t(x,y)\\
 &\quad + \int_{\RDD} \langle  \bar{A}_{\sigma_t}(y) -  \bar{A}_{\sigma_0}(x) -  B_{\sigma_0}(x)(y-x),v_t(y)\rangle d\gamma_t(x,y).
\end{align*}
By Equation (\ref{eqcontdiff}) and H\"older's inequality,
\begin{equation*}
\bigl | \int_{\RDD} \langle  \bar{A}_{\sigma_t}(y) -  \bar{A}_{\sigma_0}(x) -  B_{\sigma_0}(x)(y-x),v_t(y)\rangle d\gamma_t(x,y)\bigr | \leq  o(W_2(\sigma_0, \sigma_t))  \, ||v_t||_{\sigma_t}.
\end{equation*}
It is well known (cf. \cite{ags:book} Lemma 7.2.1) that if $0<t \leq 1$ then there exists a unique optimal transport map $T_t^1$ between $\sigma_t$ and $\sigma_1$, \textit{i.e.} $\Gamma_o(\sigma_t,\sigma_1)= \{(Id \times T_t^1)_{\#}\sigma_t \}$. One can check that $v_t(y)= \frac{T_t^1(y)-y}{1-t}$ and $||v_t||_{\sigma_t} = W_2(\sigma_t,\sigma_1)/(1-t) =  W_2(\sigma_0,\sigma_1).$ Thus
\begin{align*}
\int_{\RDD} \langle  \bar{A}_{\sigma_0}(x), v_t(y)-v(x)\rangle d\gamma_t(x,y) &=  \int \langle  \bar{A}_{\sigma_0}(x), \frac{T_t^1(y)-y}{1-t}-(\bar{\gamma}(x)-x)\rangle d\gamma_t(x,y)\\
&= \int\langle  \bar{A}_{\sigma_0}(x), \frac{z-((1-t)x + t z)}{1-t}- (z-x)\rangle d\gamma(x,z) \\
&= 0.
\end{align*}
Similarly, 
\begin{align*}
\int_{\RDD}   \langle B_{\sigma_0}(x)(y-x),v_t(y)\rangle d\gamma_t(x,y) &=  t\int_{\RDD}   \langle B_{\sigma_0}(x)(z-x),z-x\rangle d\gamma(x,y)\\
& =  o(W_2(\sigma_0,\sigma_1)) = o(W_2(\mu,\nu)).
\end{align*}
Combining these equations shows that
\begin{align}\label{eq:estimatebis}
\bar \Lambda_{\sigma_t}(v_t) - \bar \Lambda_{\sigma_0}(v)  &= o(W_2(\mu,\nu)).
\end{align}
Notice that (\ref{eq:estimatebis}) is independent of $t.$ Combining (\ref{eq:Fdifference}) and (\ref{eq:estimatebis}) we find
\begin{align*}
F(\nu) &= F(\mu) + \bar \Lambda_{\sigma_0}(v) + \int_0^1 \bar \Lambda_{\sigma_t}(v_t)- \bar \Lambda_{\sigma_0}(v) dt\\
&= F(\mu) + \bar \Lambda_{\sigma_0}(v)+ o(W_2(\mu,\nu))\\
&= F(\mu) + \int_{\RDD}\langle \bar{A}_{\sigma_0}(x),y-x\rangle d\gamma(x,y)+ o(W_2(\mu,\nu)).
\end{align*}
As in Definition \ref{def:Fdifferentiable}, this proves that $F$ is differentiable and that $\nabla_\mu F= \pi_\mu(\bar{A}_\mu).$ Thus $d F=\Lambda.$
\endproof

\begin{example} 
Assume $\bar\Lambda$ is a linear pseudo 1-form, \textit{i.e.} $\bar\Lambda(\cdot)=\int_{\RD}\langle\bar A,\cdot\rangle d\mu$ for some $\bar A\in\cfields$. According to Example \ref{e:linearimpliesdiff}, if $d\bar\Lambda=0$ then $\int_{\RD}\langle(\nabla\bar A-\nabla\bar A^T)\,\cdot,\cdot\rangle d\mu=0$ on $T_\mu\M$. Restricting to Dirac measures proves that $\nabla \bar A$ is symmetric so $\bar A$ is a gradient vector field. In other words, any closed linear pseudo 1-form is actually a linear 1-form. 
\end{example}

\subsection{Example: 1-forms on the space of discrete measures} Fix an integer $n \geq 1$. Given $x_1, \cdots, x_n \in \RD$, set $\x:=(x_1, \cdots, x_n)$ and $\mu_\x:= 1/n \sum_{i=1}^n \delta_{x_i}$. Let $M$ denote the set of such measures and $TM$ denote its tangent bundle, cf.  Examples \ref{e:discrete1} and  \ref{ex:discreteextension}. Choose a regular pseudo $1$-form $\bar \Lambda$ on $\M$. By restriction we obtain a $1$-form $\alpha$ on $M$, defined by $\alpha_\x:=\bar \Lambda_{\mu_\x}.$
Let $A: \RnD \rightarrow \RnD$ be defined by
$$A(\x)=(A_1(\x), \cdots, A_n(\x)):=\Bigl(\bar A_{\mu_\x}(x_1), \cdots, \bar A_{\mu_\x}(x_n) \Bigr).$$
Notice that if $X=(X_1, \cdots, X_n) \in \RnD$ satisfies $X_i=X_j$ whenever $x_i=x_j$ then $\alpha_\x(X)={1 \over n} \langle A(\x), X \rangle.$ Now define a $nD \times nD$ matrix $B(\x)$ by setting
\begin{equation}\label{eq:discretecontdiffbis}
B_{k+i, k+j}:= \Bigl( B_{\mu_\x} (x_{k+1}) \Bigr)_{ij}, \qquad \mbox{for } k=0, \cdots n-1, \quad i,j=1, \cdots, D,
\end{equation}
\begin{equation}\label{eq:discretecontdiffter}
B_{l, m}:=0 \quad \hbox{if} \quad (l,m) \not \in \{ (k+i, k+j): k=0, \cdots n-1, \quad i,j=1, \cdots, D\}.
\end{equation}

\begin{prop}\label{pr:diffconA} The map $A: \RnD \rightarrow \RnD$ is differentiable and $\nabla A(\x)=B(\x)$ for $\x \in \RnD.$
\end{prop}
\proof{} Let $\x=(x_1, \cdots, x_n) \in \RnD$. Set $r:= \min_{x_i \not=x_j} |x_i-x_j|.$ If $\y=(y_1, \cdots, y_n) \in \RnD$ and $|\y -\x| <r/2$ then $\Gamma_o(\mu_\x, \mu_\y)$ has a single element $\gamma_\y=1/n \sum_{i=1}^n \delta_{(x_i,y_i)}$ and $n W_2^2(\mu_\x, \mu_\y)= |\y -\x|^2.$ By Equation \ref{eqcontdiff},
\begin{equation}\label{eq:discretecontdiff}
|    A(\y)-A(\x) -B(\x)(\y -\x)|^2=n \; o({|\y-\x|^2 \over n}).
\end{equation}
This concludes the proof. \endproof
\begin{lemma} \label{le:discretetwoform} Suppose $\x=(x_1, \cdots, x_n) \in \RnD$ and $X=(X_1, \cdots, X_n)$, $Y=(Y_1, \cdots, Y_n) \in \RnD$ are such that $X_i=X_j$, $Y_i=Y_j$ whenever $x_i=x_j.$ Then
$$
\quad d\bar \Lambda_{\mu_\x}(X,Y)= d\alpha_\x(X,Y).
$$
\end{lemma}
\proof{} We use Lemma \ref{eq10.13.9} and Equations \ref{eq:discretecontdiffbis}, \ref{eq:discretecontdiffter} to obtain
$$
d\bar \Lambda_{\mu_\x}(X,Y)= \sum_{k=1}^n \Bigl\langle  (B_{\mu_\x}(x_k) -B_{\mu_\x}(x_k)^T) X_k, Y_k\Bigr\rangle= d\alpha_\x(X,Y).
$$\endproof

\begin{corollary}\label{co:greendiscrete} Suppose that $\rb=(r_1, \cdots, r_n) \in C^2([0,T], \RnD)$  and set $\sigma_t:=1/n \sum_{i=1}^n \delta_{r_i(t)}.$ If $\bar\Lambda$ is closed and $\sigma_0=\sigma_T$ then $\int_{\sigma} \alpha=0.$
\end{corollary}
\proof{} This is a direct consequence of Corollary \ref{co:closecurve}.
\endproof

\begin{remark}\label{re:greenformula} One can check by direct computation that, for a surface $\x=\x(s,t)$ in $M$, the familiar identity
$\partial_t ( \alpha_{\x}(\partial_s \x) )- \partial_s ( \alpha_{\x}(\partial_t \x) ) =d\alpha_{\x}(\partial_t \x,  \partial_s \x)$ holds. Together with Lemma \ref{le:discretetwoform}, this is the analogue of Corollary \ref{le:lipschitzgreen3} which we used to prove Theorem \ref{th:greencurveannulus}. 
\end{remark}

\begin{remark} Notice that the assumption $\sigma_0=\sigma_T$ is weaker than $\rb(0)=\rb(T)$.
\end{remark}


\subsection{Discussion}\label{ss:section5discussion}
As mentioned in Section \ref{ss:topology}, the space $\M$ is convex. This is true not only for probability measures on $\R^D$ but also for probability measures on any manifold $M$. It is thus trivially true that all cohomology groups $H^k(\M;\R)$ vanish if these groups are defined purely in terms of the topology on $\M$.

In this paper, however, we are concerned with the differentiable structure on $\M$. As seen in Section \ref{ss:tangentspacesbis}, from this point of view $\M$ is the union of smooth manifolds $\mathcal{O}$, defined as the orbits of the $\cdiffrD$-action on $\M$. Given any such orbit $\mathcal{O}$, it is then interesting to define and calculate the first cohomology group $H^1(\mathcal{O};\R)$. Notice that these orbits are in general not convex, so the above reasoning does not apply. In the case $M=\R^D$, however, Corollary \ref{co:closedexact} can intuitively be interpreted as a vanishing result for $H^1(\mathcal{O};\R)$ if we think of cohomology in the differentiable sense, \textit{i.e.} in the sense of de Rham, as follows.

Let $\Lambda$ be a regular 1-form on $\M$. Recall from Definition \ref{d:differential} that $\Lambda$ is closed iff $\bar\Lambda$ is closed. Then Corollary \ref{co:closedexact} shows that if $\Lambda$ is closed then it is exact, \textit{i.e.} $\Lambda=dF$ for some differentiable $F:\M\rightarrow\R$. Now choose an orbit $\mathcal{O}$. Given any $\mu\in\mathcal{O}$, recall from Section \ref{ss:tangentspacesbis} that $T_\mu\M$ can be thought of as the tangent space to $\mathcal{O}$ at the point $\mu$. In this sense the restriction of $F$ to $\mathcal{O}$ is still differentiable and $d(F_{|\mathcal{O}})=(dF)_{|\mathcal{O}}$ so the restriction of $\Lambda$ to $\mathcal{O}$ is also exact. Roughly speaking, Corollary \ref{co:closedexact} thus shows that the first \textit{de Rham cohomology} group $H^1(\mathcal{O};\R)$ of each orbit vanishes. It may be useful to point out that if $\mu$ is a Dirac measure then $\mathcal{O}_\mu=\R^D$, so at least in this case the above vanishing result makes sense. 

It is reasonable to expect that most of the theory of \cite{ags:book} can be extended to probability measure spaces on other manifolds $M$. In this case, many of the results of this paper should also extend. However the above example, where $\mathcal{O}$ is the space of Dirac measures on $M$, shows that one should not expect $H^1(\mathcal{O};\R)$ to vanish in general. In this sense our results are specific to the choice $M=\R^D$. Another way to see this is as follows. The proof of Corollary \ref{co:closedexact} relied on the construction of certain ``annuli" built using maps $D_s$ of $\R^D$. These maps exist only because $\R^D$ is contractible. Such a construction would not be possible on other manifolds.

The following considerations also support the above interpretation of Corollary \ref{co:closedexact}. Recall that, for a finite-dimensional manifold $M$, the first de Rham cohomology group is closely related to the topology of $M$, as follows: $H^1(M;\R)=\mbox{Hom}(\pi_1(M),\R)$, where the latter denotes the space of group homomorphisms from the first fundamental group $\pi_1(M)$ to $\R$. In our case, an orbit $\mathcal{O}$ is generally not a manifold in any rigorous sense so it is not \textit{a priori} clear that there exists any relationship between our $H^1(\mathcal{O};\R)$ and $\pi_1(\mathcal{O})$. However, we can formally prove the topological counterpart of our vanishing result as follows. 

Let $G$ be a finite-dimensional Lie group and $H$ be a closed subgroup. Recall that there exists a \textit{homotopy long exact sequence}
$$\dots\rightarrow\pi_1(H)\rightarrow \pi_1(G)\rightarrow \pi_1(G/H) \rightarrow \pi_0(H) \rightarrow \pi_0(G)\dots,$$
cf. \textit{e.g.} \cite{bredon}, VII.5. Now assume $G$ is connected, \textit{i.e.} $\pi_0(G)=1$. We can then dualize the final part of this sequence obtaining a new exact sequence
\begin{equation}\label{eq:dualhomotopysequence}
1\rightarrow\mbox{Hom}(\pi_0(H),\R)\rightarrow \mbox{Hom}(\pi_1(G/H),\R)\rightarrow\mbox{Hom}(\pi_1(G),\R).
\end{equation}
Now choose $\mu\in\mathcal{O}$ and set $G:=\cdiffrD$ and $H:=\cdiffmurD$ so that $G/H\simeq\mathcal{O}$. In many cases it is known that $\pi_1(G)$ is finite: specifically, this is true at least for $D=1,2,3$ and $D\geq 12$, cf. \cite{antonellietal} for related results. Let us assume that $H$ has a finite number of components and that the homotopy long exact sequence is still valid in this infinite-dimensional setting.  Sequence \ref{eq:dualhomotopysequence} then becomes 
$$1\rightarrow 1\rightarrow \mbox{Hom}(\pi_1(\mathcal{O}),\R)\rightarrow 1,$$
so by exactness $\mbox{Hom}(\pi_1(\mathcal{O}),\R)$ must also be trivial.

To conclude, it is also interesting to examine the relationship between Corollary \ref{co:closedexact} and \textit{invariant cohomology}, in the sense of Section \ref{ss:cohomology}. Recall from Proposition \ref{prop:2cohoms} that the first cohomology group of an orbit is a subgroup of the corresponding first invariant cohomology group. The statement that the first invariant cohomology vanishes is thus stronger than the statement that the first cohomology of the orbit vanishes. Now choose any orbit $\mathcal{O}$ in $\M$. According to Section \ref{ss:section4discussion}, the corresponding invariant cohomology should be defined in terms of regular pseudo 1-forms. To prove that the invariant cohomology vanishes would thus require showing that $d\bar\Lambda=0$ implies $\bar\Lambda=dF$, for some differentiable function $F:\M\rightarrow \R$. Since $dF$ is a 1-form, such a result would imply that any closed pseudo 1-form is a 1-form. Corollary \ref{co:closedexact} does not achieve this. On the other hand, it is not clear that such a result should even be expected.


\section{A symplectic foliation of $\M$}\label{s:symplectic}
In Section \ref{ss:tangentspacesbis} we used the action of the group of diffeomorphisms $\cdiffrD$ to build a foliation of $\M$: this allowed us to formally reconstruct the differential calculus on $\M$. We now specialize to the case $D=2d$. In this case the underlying manifold $\Rtwod$ has a natural extra structure, the \textit{symplectic structure} $\omega$. The goal of this section is to use this extra data to build a second, finer, foliation of $\M$; we then prove that each leaf of this foliation admits a symplectic structure $\Omega$. The foliation is obtained via a smaller group of diffeomorphisms defined by $\omega$, the \textit{Hamiltonian diffeomorphisms}. Section \ref{ss:hamdiffs} provides an introduction to this group, cf. \cite{mcduffsalamon:book} or \cite{marsdenratiu:book} for details.

\subsection{The group of Hamiltonian diffeomorphisms} \label{ss:hamdiffs} Recall that a \textit{symplectic structure} on a (possibly infinite-dimensional) vector space $V$ is a 2-form $\omega:V\times V\rightarrow \R$
such that 
\begin{equation} \label{eq:defsymplectic}
\omega^\flat:V\rightarrow V^*,\ \ v\mapsto i_v\omega:=\omega(v,\cdot)
\end{equation}
is injective. If $V$ is finite-dimensional then $\omega^\flat$ is an isomorphism; we will denote its inverse by $\omega^\sharp$. 

Let $M$ be a smooth manifold of dimension $D:=2d$. A \textit{symplectic structure} on $M$ is a smooth closed 2-form $\omega$ satisfying Equation \ref{eq:defsymplectic} at each tangent space $V=T_xM$; equivalently, such that $\omega^d$ is a volume form on $M$. Notice that, since $d\omega=0$, Cartan's formula \ref{eq:cartan} shows that $\mathcal{L}_X\omega=di_X\omega$. Throughout this section, to simplify notation, we will drop the difference between compact and noncompact manifolds but the reader should keep in mind that in the latter case we always silently restrict our attention to maps and vector fields with compact support.

Consider the set of \textit{symplectomorphisms} of $M$, \textit{i.e.}
$$\sympm:=\{\phi\in\diffm:\phi^*\omega=\omega\}.$$
This is clearly a subgroup of $\diffm$. Using the methods of Section \ref{ss:diffM} (see in particular Remark \ref{r:charts}) one can show that it has a Lie group structure. Its tangent space at $Id$, thus its Lie algebra, is by definition isomorphic to the space of closed 1-forms on $M$. Via $\omega^\sharp$ and Formula \ref{eq:cartan} this space is isomorphic to the space of \textit{symplectic} or \textit{locally Hamiltonian} vector fields, \textit{i.e.}
$$\sympx:=\{X\in\xm:\mathcal{L}_X\omega=0\}.$$
\begin{remark}
Equation \ref{eq:liederivativebracket} confirms that $\sympx$ is closed under the bracket operation, \textit{i.e.} that it is a Lie subalgebra of $\xm$. Equation \ref{eq:liederivativepullback} confirms that $\sympx$ is closed under the push-forward operation, \textit{i.e.} under the adjoint representation of $\sympm$ on $\sympx$, cf. Lemma \ref{l:adjointispushforward}.
\end{remark}

We say that a vector field $X$ on $M$ is \textit{Hamiltonian} if the corresponding 1-form $\xi:=\omega(X,\cdot)$ is exact: $\xi=df$. We then write $X=X_f$. This defines the space of \textit{Hamiltonian vector fields} $\hamfields$. It is useful to rephrase this definition as follows. Consider the map
\begin{equation} \label{eq:poissonmap}
\cm\rightarrow\xm,\ \ f\mapsto df\simeq X_f:=\omega^\sharp(df).
\end{equation}
The Hamiltonian vector fields are the image of this map. This map is linear. It is not injective: its kernel is the space of functions constant on $M$. In Section \ref{ss:poisson} we will start referring to these functions as the \textit{Casimir functions} for the map of Equation \ref{eq:poissonmap}.

\begin{remark}\label{r:normalization}
We can rephrase the properties of the map of Equation \ref{eq:poissonmap} by saying that there exists a \textit{short exact sequence}
\begin{equation}\label{eq:hamsequence}
0\rightarrow\R\rightarrow\cm\rightarrow\hamfields\rightarrow 0.
\end{equation}
As already mentioned, the function corresponding to a given Hamiltonian vector field is well-defined only up to a constant. In some cases we can fix this constant via a \textit{normalization}, \textit{i.e.} we can build an inverse map $\hamfields\rightarrow \cm$. We then obtain an isomorphism between $\hamfields$ and the space of normalized functions. For example, if $M$ is compact we can fix this constant by requiring that $f$ have integral zero, $\int_M f \omega^d=0$. If instead $M=\Rtwod$ and we restrict our attention as usual to Hamiltonian diffeomorphisms with compact support, we should restrict Equation \ref{eq:poissonmap} to the space $\R\oplus C^\infty_c(\Rtwod)$ of functions which are constant outside of a compact set; by restriction we then get an isomorphism $C^\infty_c(\Rtwod)\simeq\chamfields$.
\end{remark}

More generally, we say that a time-dependent vector field $X_t$ is \textit{Hamiltonian} if $\omega(X_t,\cdot)=df_t$ for some curve of smooth functions $f_t$. We then say that the diffeomorphism $\phi\in \diffm$ is \textit{Hamiltonian} if it can be obtained as the time $t=1$ flow of a time-dependent Hamiltonian vector field $X_{f_t}$, \textit{i.e.} if $\phi=\phi_1$ and $\phi_t$ solves Equation \ref{eq:tflow}. 

Let $\hamm$ denote the set of Hamiltonian diffeomorphisms. It follows from Lemma \ref{l:preservedform} that all such maps are symplectomorphisms. It is not immediately obvious that $\hamm$ is closed under composition but it is not hard to prove that this is indeed true, cf. \cite{mcduffsalamon:book} Proposition 10.2 and Exercise 10.3. Once again, the methods of Section \ref{ss:diffM} and Remark \ref{r:charts} show that $\hamm$ has a Lie group structure. Its tangent space at $Id$, thus its Lie algebra, is isomorphic to the space of exact 1-forms, which via $\omega^\sharp$ corresponds to the space of Hamiltonian vector fields.

It is a fundamental fact of Symplectic Geometry that $\omega$ defines a Lie bracket on $\cm$, as follows:
$$\{f,g\}:=\omega(X_f,X_g)=df(X_g)=\mathcal{L}_{X_g}f.$$

This operation is clearly bilinear and anti-symmetric. The fact that it satisfies the Jacobi identity, cf. Definition \ref{def:liealgebra}, follows from the following standard result.

\begin{lemma} \label{l:hamfieldproperties}
Let $\phi\in \sympm$. Then $\phi^*X_f=X_{\phi^*f}$ and $\phi^*\{f,g\}=\{\phi^*f,\phi^*g\}$. Applying this to $\phi_t\in \sympm$ and differentiating, it implies:
\begin{equation}\label{eq:bracketonfunctions}
\mathcal{L}_{X_h}\{f,g\}=\{\mathcal{L}_{X_h}f,g\}+\{f,\mathcal{L}_{X_h}g\}.
\end{equation}
\end{lemma}

\begin{lemma} \label{l:mapislie}
The map $f\mapsto X_f$ has the following property:
$$X_{\{f,g\}}=-[X_f,X_g].$$
\end{lemma}

\proof{} It is enough to check that $dh(X_{\{f,g\}})=-dh([X_f,X_g])$, for all $h\in\cm$. As usual, it will simplify the notation to set $X(f):=df(X)$. In particular $X_f(h)=\{h,f\}$ and $dh([X,Y])=X(Y(h))-Y(X(h))$. Then:
\begin{align*}
X_{\{f,g\}}(h) &= \{h,\{f,g\}\}=-\{f,\{g,h\}\}-\{g,\{h,f\}\}\\
&=-\{\{h,g\},f\}+\{\{h,f\},g\}= -X_f(X_g(h))+X_g(X_f(h))\\
&=-[X_f,X_g](h).
\end{align*}
\endproof

Recall from Section \ref{ss:diffM} the negative sign appearing in the Lie bracket $[\cdot,\cdot]_\mathfrak{g}$ on vector fields. It follows from Lemma \ref{l:mapislie} that the map of Equation \ref{eq:poissonmap} is a Lie algebra homomorphism between $\cm$ and the space of Hamiltonian vector fields, endowed with that Lie bracket.

\begin{remark} Lemma \ref{l:mapislie} confirms that $\hamfields$ is a Lie subalgebra of $\xm$. Lemma \ref{l:hamfieldproperties} confirms that it is closed under symplectic push-forward, so in particular it is closed under the adjoint representation of $\hamm$.
\end{remark}

\begin{remark} \label{r:hamvssymp}
Notice that $\hamm$ is connected by definition. If $M$ satisfies $H^1(M,\R)=0$, \textit{i.e.} every closed 1-form is exact, then every symplectic vector field is Hamiltonian. Now assume that $\phi\in\sympm$ is such that there exists $\phi_t\in\sympm$ with $\phi_0=Id$ and $\phi_1=\phi$. It then follows from Lemma \ref{l:preservedform} that $\phi$ is Hamiltonian, \textit{i.e.} that the connected component of $\sympm$ containing the identity coincides with $\hamm$. In particular this applies to $M=\R^{2d}$, so in later sections we could just as well choose to work with (the connected component containing $Id$ of) $\symprtwod$ rather than with $\hamrtwod$. We choose however not to do this, so as to emphasize the fact that for general $M$ the two groups are indeed different and that generalizing our constructions would require working with $\hamm$ rather than with $\sympm$.
\end{remark}

\begin{remark} 
In many cases it is known that $\sympm$ is closed in $\diffm$ and that $\hamm$ is closed in $\sympm$, see \cite{mcduffsalamon:book} and \cite{ono:flux} for details. 
\end{remark}


\subsection{A symplectic foliation of $\M$}  \label{ss:symplecticfoliation} The manifold $\R^{2d}$ has a natural symplectic structure defined by $\omega:=dx^i\wedge dy^i$. Let $J$ denote the natural \textit{complex structure} on $\R^{2d}$, defined with respect to the basis $\partial x^1,\dots,\partial x^d,\partial y^1,\dots,\partial y^d$ by the matrix
\begin{align*}
J &= \left(\begin{array}{cc} 
0 & -I\\
I & 0
\end{array}\right).
\end{align*}
Notice that $\omega(\cdot,\cdot)=g(J\cdot,\cdot)$. It follows from this that Hamiltonian vector fields on $\R^{2d}$ satisfy the identity 
\begin{equation}\label{eq:stdhamfields}
X_f=-J\nabla f. 
\end{equation}
Set $\G:=\hamrtwod$, the group of compactly-supported Hamiltonian diffeomorphisms on $\R^{2d}$. Let $\chamfields$ denote the corresponding Lie algebra, \textit{i.e.} the space of compactly supported Hamiltonian vector fields on $\Rtwod$. The push-forward action of $\cdiffrtwod$ on $\M$ restricts to an action of $\G$. The corresponding orbits and stabilizers are
$$\mathcal{O}_\mu:=\{\nu\in\M:\nu=\phi_\#\mu, \ \ \mbox{for some }\phi\in \G\}, \ \ \Gmu:=\{\phi\in\G:\phi_\#\mu=\mu\}.$$
Notice that this action provides a second foliation of $\M$, finer than the one of Section \ref{ss:tangentspacesbis}.

\begin{example} \label{e:hamorbits}
As in Example \ref{e:discrete1}, let $a_i$ $(i=1,\dots,n)$ be a fixed collection of positive numbers such that $\sum a_i=1$ and $x_1, \dots, x_n \in \Rtwod$ be $n$ distinct points. Set $\mu= \sum_{i=1}^n a_i \, \delta_{x_i} \in \M$ and
$$\calO=\biggl\{ \sum_{i=1}^n a_i \, \delta_{\bar x_i} \; : \; \bar x_1, \dots, \bar x_n \in \Rtwod \quad \hbox{are distinct}\biggr\}.$$
Since smooth Hamiltonian diffeomeorphisms are one-to-one maps of $\Rtwod$ it is clear that $\calOmu \subseteq \calO$. Given any $\bar x_1 \in \Rtwod \setminus \{x_2, \dots, x_n \}$ one can show that there exists a Hamiltonian diffeomorphism $\phi$ with compact support such that $\phi(x_1)=\bar x_1$ and $\phi(x_i)=x_i$ for $i\neq 1$. Thus, setting $\bar \mu:= a_1 \, \delta_{\bar x_1} + \sum_{i=2}^n a_i \, \delta_{x_i}$, we see that $\bar \mu \in \calOmu.$ Repeating the argument $n-1$ times we conclude that $\calO \subseteq \calOmu$, so $\mathcal{O}=\mathcal{O}_\mu$.
\end{example}

\begin{definition} \label{d:symplectictangent}
Let $\mu\in\M$. Consider the $\lfields$-closure $\lhamfields$ of $\chamfields$. We can restrict the operator $\divmu$ to this space; we will continue to denote its kernel $\kerdiv$. We define the \textit{symplectic tangent subspace} at $\mu$ to be the space
$$T_\mu\mathcal{O}:=\lhamfields/\kerdiv\subseteq \lfields/\kerdiv.$$
Recall from Remark \ref{rem:altertangentspace} the identification $\pi_\mu:\lfields/\kerdiv\rightarrow T_\mu\M$. By restriction this allows us to identify $T_\mu\mathcal{O}$ with the subspace $\pi_\mu(\lhamfields)\subseteq T_\mu\M$. We define the \textit{pseudo symplectic distribution} on $\M$ to be the union of all spaces $\lhamfields$, for $\mu\in\M$. It is a subbundle of $\mathcal{TM}$. We define the the \textit{symplectic distribution} on $\M$ to be the union of all spaces $T_\mu\mathcal{O}$, for $\mu\in\M$. Up to the above identification, it is a subbundle of $T\M$.
\end{definition}

\begin{remark}\label{rem:symplectictangent}
Recall that in general a Hilbert space projection will not necessarily map closed subspaces to closed subspaces. Thus it is not clear that $\pi_\mu(\lhamfields)$ is closed in $T_\mu\M$. On the other hand, the space $T_\mu\mathcal{O}$ of Definition \ref{d:symplectictangent} has a natural Hilbert space structure. In other words, from the Hilbert space point of view the two notions of $T_\mu\mathcal{O}$ introduced in Definition \ref{d:symplectictangent} are not necessarily equivalent. This is in contrast with the two notions of $T_\mu\M$, cf. Definition \ref{def:tangentspaces} and Remark \ref{rem:altertangentspace}.
\end{remark}

\begin{remark} 
Formally speaking the symplectic distribution is \textit{integrable} because it is the set of tangent spaces of the smooth foliation defined by the action of $\G$.
\end{remark}

\begin{example} \label{e:hamfieldsvsdecomp}
It is interesting to compare the space $\lhamfields$ to the subspaces defined by Decomposition \ref{eq10.11.1}. For example, let $\mu=\delta_x$. Recall from Example \ref{ex:discreteextension} that for any $\xi\in L^2(\mu)$ there exists $\bar{\varphi}\in\ccinf$ such that $\xi(x)=\nabla\bar{\varphi}(x)$. Thus $\lgradfields=L^2(\mu)$. Now choose any $X\in L^2(\mu)$ and apply this construction to $\xi:=JX$. Then $X(x)=-J\nabla\bar{\varphi}(x)$, so $\lhamfields=L^2(\mu)$. This is the infinitesimal version of Example \ref{e:hamorbits}. In particular, $\lhamfields=\lgradfields$. 

The ``opposite extreme'' is represented by the absolutely continuous case $\mu=\rho\mathcal{L}$, for some $\rho>0$. In this case if a Hamiltonian vector field is a gradient vector field, \textit{e.g.} $-J\nabla v=\nabla u$, then the function $u+iv$ is holomorphic on $\C^d$, so $u$ and $v$ are \textit{pluriharmonic} functions on $\R^{2d}$ in the sense of the theory of several complex variables. This is a very strong condition: in particular, it implies that $u$ and $v$ are harmonic. Thus $\chamfields\cap\gradfields=\{0\}$.

We can also compare $\lhamfields$ with $\kerdiv$. When $\mu=\delta_x$ we saw in Example \ref{ex:discreteextension} that $\kerdiv=\{0\}$, so $\lhamfields\cap\kerdiv=\{0\}$. On the other hand, assume $\mu=\rho\mathcal{L}$ for some $\rho>0$. Then $div_\mu(X)=\rho div(X)+\langle\nabla\rho,X\rangle$. Choose $X=-J\nabla f$. Then $div(X)=0$ so $div_\mu(X)=0$ iff $\langle\nabla \rho, -J\nabla f\rangle=0$. Choosing in particular $f=\rho$ shows that $\chamfields\cap\kerdiv\neq\{0\}$.
\end{example}

We now want to show that each $T_\mu\mathcal{O}$ has a natural symplectic structure; this will justify the terminology of Definition \ref{d:symplectictangent}. We rely on the following general construction.

\begin{definition}\label{def:symplecticquotient}
Let $(V,\omega)$ be a symplectic vector space. Let $W$ be a subspace of $V$. In general the restriction of $\omega$ to $W$ will not define a symplectic structure on $W$ because $\omega^\flat:W\rightarrow W^*$ will not be injective. However, set $Z:=\{w\in W:\omega(w,\cdot)_{|W}\equiv 0\}$. Then $\omega$ \textit{reduces} to a symplectic structure on the quotient space $W/Z$, defined by
$$\omega([w],[w']):=\omega(w,w').$$
\end{definition}

In our case we can set $V:=\lfields$ and $W:=\lhamfields$. The 2-form
\begin{equation}\label{def:2form}
\hat{\Omega}_\mu(X,Y):=\int_{\Rtwod}\omega(X,Y)\,d\mu
\end{equation}
defines a symplectic structure on $\lfields$. The restriction of $\hat{\Omega}_\mu$ defines a 2-form
$$\bar{\Omega}_\mu:\lhamfields\times\lhamfields\rightarrow\R.$$
Notice that $\hat{\Omega}_\mu$ is continuous in the sense of Definition \ref{d:kforms}, so $\bar{\Omega}_\mu$ can also be defined as the unique continuous extension of the 2-form
\begin{equation}\label{eq:pseudosymplectic}
\bar{\Omega}_\mu:\chamfields\times\chamfields\rightarrow\R, \ \ \bar{\Omega}_\mu(X_f,X_g):=\int\omega(X_f,X_g)\,d\mu.
\end{equation}
Notice also that, for any $X\in \lfields$,
\begin{equation}\label{eq:degeneracy}
\int\omega(X,X_f)\,d\mu=-\int df(X)\,d\mu=\langle\divmu(X),f\rangle
\end{equation} 
so $\int\omega(X,\cdot)\,d\mu\equiv 0$ on $\lhamfields$ iff $X\in \kerdiv$. This calculation shows that the space $Z$ of Definition \ref{def:symplecticquotient} coincides with the space $\kerdiv\cap\lhamfields$. We can now define $\Omega_\mu$ to be the reduced symplectic structure on the space $T_\mu\mathcal{O}=W/Z$. In terms of the identification $\pi_\mu$, this yields
\begin{equation}\label{eq:symplectic}
\Omega_\mu:T_\mu\mathcal{O}\times T_\mu\mathcal{O}\rightarrow\R,\ \ \Omega_\mu(\pi_\mu(X_f),\pi_\mu(X_g)):=\int\omega(X_f,X_g)\,d\mu.
\end{equation}
Notice that $\Omega_\mu$ is not necessarily continuous in the sense of Definition \ref{d:kforms}.

Using Equation \ref{eq:stdhamfields} we can also write this as
$$\Omega_\mu(\pi_\mu(X_f),\pi_\mu(X_g))=\int\omega(J\nabla f,J\nabla g)\,d\mu=\int g(J\nabla f,\nabla g)\,d\mu.$$ 
We now want to understand the geometric and differential properties of $\bar{\Omega}$. It is simple to check that $\bar{\Omega}$ is $\G$-invariant, in the sense that $\phi^*\bar{\Omega}=\bar{\Omega}$, for all $\phi\in\G$. Indeed, using Definition \ref{d:pullback} and Lemma \ref{l:hamfieldproperties},
\begin{align*}
(\phi^*\bar{\Omega})_\mu(X_f,X_g) &= \bar{\Omega}_{\phi_{\#}\mu}(\phi_*(X_f),\phi_*(X_g))=\int_{\Rtwod}\omega(X_{f\circ\phi^{-1}},X_{g\circ\phi^{-1}})\,d\phi_\#\mu\\
&= \int_{\Rtwod}\{f\circ\phi^{-1},g\circ\phi^{-1}\}\,d\phi_\#\mu = \int_{\Rtwod}\{f,g\}\circ\phi^{-1}\,d\phi_\#\mu\\
&= \bar{\Omega}_\mu(X_f,X_g).
\end{align*}
It then follows that $\Omega$ is also $\G$-invariant. 

\begin{lemma}
Given any $X,Y,Z\in\chamfields$,
\begin{align*}
X\bar{\Omega}(Y,Z)-Y\bar{\Omega}(X,Z)+Z\bar{\Omega}(X,Y)-\bar{\Omega}([X,Y],Z)+\bar{\Omega}([X,Z],Y)-\bar{\Omega}([Y,Z],X)&=0.
\end{align*}
\end{lemma}
\proof{}
Notice that $\bar{\Omega}(Y,Z)$ is a linear function on $\M$ in the sense of Example \ref{e:deflinearforms}. It is thus differentiable, cf. Example \ref{e:linearisdifferentiable}, and $X\bar{\Omega}(Y,Z)=\int X\omega(Y,Z)\,d\mu$. It follows that the left hand side of the above equation reduces to $\int d\omega(X,Y,Z)\,d\mu$, which vanishes because $\omega$ is closed.

\endproof

This shows that $\bar{\Omega}$ is differentiable and closed in the sense analogous to Definition \ref{d:differential}, \textit{i.e.} using Equation \ref{eq:orbitd} (or Equation \ref{eq:usuald}) with $k=2$ instead of $k=1$. Using the terminology of Section \ref{ss:forms} we can say that $\bar{\Omega}$ is a closed pseudo linear 2-form defined on the pseudo distribution $\mu\rightarrow \lhamfields$ of Definition \ref{d:symplectictangent}. 

\begin{remark} \label{rem:symplectictangentbis}
As in Remark \ref{rem:symplectictangent}, it may again be useful to emphasize a possible misconception related to the identification $\pi_\mu:\lhamfields/\kerdiv\simeq \pi_\mu(\lhamfields)$. One could also restrict $\hat{\Omega}_\mu$ to the subspace $W':=\pi_\mu(\lhamfields)$, obtaining a 2-form
$$\Omega'_\mu(\pi_\mu(X_f),\pi_\mu(X_g))=\int\omega(\pi_\mu(J\nabla f),\pi_\mu(J\nabla g))\,d\mu.$$
It is important to realize that $\Omega'_\mu$ does \textit{not} coincide, under $\pi_\mu$, with $\Omega_\mu$. Specifically, 
$\Omega'_\mu$ differs from $\Omega_\mu$ in that it does not take into account the divergence components of $X_f$, $X_g$.

In the framework of \cite{ags:book} it is more natural to work in terms of the space $\pi_\mu(\lhamfields)\subseteq T_\mu\M$ than in terms of $\lhamfields/\kerdiv$. From this point of view, the choice of $\Omega_\mu$ as a symplectic structure on $T_\mu\mathcal{O}$ may seem less natural than the choice of $\Omega'_\mu$. The fact that $\Omega_\mu$ is even well-defined on $T_\mu\mathcal{O}$ follows only from Equation \ref{eq:degeneracy}. Our reasons for preferring $\Omega_\mu$ are based on its geometric and differential properties seen above. Together with Remark \ref{rem:symplectictangent}, this shows that from a symplectic viewpoint the identification $\pi_\mu$ is not natural.
\end{remark}

We can now define the concept of a Hamiltonian flow on $\M$ as follows.

\begin{definition}\label{def:hamflow}
Let $F:\M\rightarrow \R$ be a differentiable function on $\M$ with gradient $\nabla F$. We define the \textit{Hamiltonian vector field} associated to $F$ to be $X_F(\mu):=\pi_\mu(-J\nabla F)$. A \textit{Hamiltonian flow} on $\M$ is a solution to the equation 
$$\frac{\partial\,\mu_t}{\partial t}=-div_{\mu_t}(X_F).$$
\end{definition}

Let $X_F$ be a Hamiltonian vector field. Choose any $Y\in T_\mu\mathcal{O}$. We can then write $Y=\pi_\mu(-J\tilde{Y})$, for some $\tilde{Y}\in\lgradfields$. Notice that 
\begin{equation*}
\Omega_\mu(X_F,Y)=\hat{\Omega}_\mu(\nabla F,\tilde{Y})=\hat{G}_\mu(\nabla F,-J\tilde{Y})=\hat{G}_\mu(\nabla F,Y)=dF(Y).
\end{equation*}
This shows that $X_F$ satisfies the analogue of Equation \ref{eq:poissonmap} along $\mathcal{O}$, justifying the terminology of Definition \ref{def:hamflow}. We refer to \cite{ambrosiogangbo:hamiltonian} and to \cite{hwakil:thesis} for specific results concerning Hamiltonian flows on $\M$.


\subsection{Algebraic properties of the symplectic distribution}\label{ss:algebraic} 
Regardless of Remarks \ref{rem:symplectictangent} and \ref{rem:symplectictangentbis}, from the point of view of \cite{ags:book} it is interesting to understand the linear-algebraic properties of the symplectic spaces $(T_\mu\mathcal{O},\Omega_\mu)$, viewed as subspaces $\pi_\mu(\lhamfields)\subseteq T_\mu\M$. Throughout this section we will use this identification. We will mainly work in terms of the complex structure $J$ on $\Rtwod$ and of certain related maps. This will also serve to emphasize the role played by $J$ within this theory. The key to this construction is of course the peculiar compatibility between the standard structures $g:=\langle\cdot,\cdot\rangle$, $\omega$ and $J$ on $\R^{2d}$, which we emphasize as follows.

\begin{definition}\label{def:compatible}
Let $V$ be a vector space endowed with both a metric $g$ and a symplectic structure $\omega$. Notice that there exists a unique injective $A\in End(V)$ such that $\omega(\cdot,\cdot)=g(A\cdot,\cdot)$. Using the isomorphism $g^\sharp:V^*\rightarrow V$ induced by $g$, $A=g^\sharp\circ\omega^\flat$ so $A$ is surjective iff $\omega^\flat$ is an isomorphism.

The fact that $\omega$ is anti-symmetric implies that $A$ is anti-selfadjoint, \textit{i.e.} $A^*=-A$. We say that $(\omega,g)$ are \textit{compatible} if $A$ is an isometry, \textit{i.e.} $A^*=A^{-1}$. In this case $A^2=-Id$, \textit{i.e.} $A$ is a \textit{complex structure} on $V$. A subspace $W\subseteq V$ is \textit{symplectic} if the restriction of $\omega$ to $W$ defines a symplectic structure on $W$. In particular, if $g$ and $\omega$ are compatible than any complex subspace of $V$ is symplectic.

The analogous definitions hold for a smooth manifold $M$ endowed with a Riemannian metric $g$ and a symplectic structure $\omega$. In general, given any function $f$ on $M$, the Hamiltonian vector field $X_f$ is related to the gradient field $\nabla f$ as follows: $X_f=A^{-1}\nabla f$. If $g$ and $\omega$ are compatible then $X_f=-A\nabla f$.
\end{definition}

The standard structures $g$ and $\omega$ on $\R^{2d}$ are of course the primary example of compatible structures. Given any $\mu\in\M$, $\hat{G}_\mu$ and $\hat{\Omega}_\mu$ (defined in Equations \ref{eq:innerproduct} and \ref{def:2form}) are compatible structures on $L^2(\mu)$. In this case the corresponding automorphism is the isometry
$$J:L^2(\mu)\rightarrow L^2(\mu),\ \ (JX)(x):=J(X(x)).$$

\begin{example} 
Notice that $\lhamfields=-J(T_\mu\M)$. Thus $\lhamfields$ is $J$-invariant iff $T_\mu\M$ is $J$-invariant iff $T_\mu\mathcal{O}=T_\mu\M$. In this case, $\bar{\Omega}_\mu=\Omega_\mu=\Omega'_\mu$. Example \ref{e:hamfieldsvsdecomp} shows that this is the case if $\mu$ is a Dirac measure. Example \ref{e:hamfieldsvsdecomp} also shows that if $\mu=\rho\mathcal{L}$ for some $\rho>0$ then the space $\gradfields$ is \textit{totally real}, \textit{i.e.} $J(\gradfields)\cap\gradfields=\chamfields\cap\gradfields=\{0\}$.
\end{example}

Our first goal is to characterize the orthogonal complement of the closure of $T_\mu\mathcal{O}$ in $T_\mu\M$. To this end, recall that any continuous map $P:H\rightarrow H$ on a Hilbert space $H$ satisfies $\mbox{Im} (P)^\perp=\mbox{Ker} (P^*)$, where $P^*:H\rightarrow H$ is the adjoint map. This yields an orthogonal decomposition $H=\overline{\mbox{Im}(P)}\oplus \mbox{Ker} (P^*)$. 

\begin{example} One example of such a decomposition is Decomposition \ref{eq10.11.1}, corresponding to the map $P:=\pi_\mu$ defined on $H:=\lfields$: in this case $\mbox{Im}(P)$ is closed and $\pi_\mu$ is self-adjoint so $\mbox{Ker} (P^*)=\mbox{Ker} (\pi_\mu)$. 

Another example is provided by the map $P:=\pi_\mu\circ J$, again defined on $\lfields$. In this case $P^*=-J\circ\pi_\mu$ and $\mbox{Im}(P)=\mbox{Im}(\pi_\mu)$, $\mbox{Ker} (P^*)=\mbox{Ker}(\pi_\mu)$ so the decomposition corresponding to $P$ again coincides with Decomposition \ref{eq10.11.1}. On the other hand, $\mbox{Im}(P^*)=-J(\mbox{Im}(\pi_\mu))$ and $\mbox{Ker} (P)=J^{-1}(\mbox{Ker}(\pi_\mu))=-J(\mbox{Ker}(\pi_\mu))$ so the decomposition corresponding to $P^*$ is the $(-J)$-rotation of Decomposition \ref{eq10.11.1}, \textit{i.e.}
\begin{equation*}\label{eq:decompbis}
\lfields=\mbox{Im}(P^*)\oplus \mbox{Ker}(P)=-J(\lgradfields)\oplus -J(\kerdiv).
\end{equation*}
In other words, there exists an orthogonal decomposition
\begin{equation}
\lfields = \lhamfields\oplus Ker(\pi_\mu\circ J).
\end{equation}
\end{example}

\begin{lemma}\label{l:orthogonalspaces}For any $\mu\in\M$ there exists an orthogonal decomposition
\begin{equation}
T_\mu\M = \overline{T_\mu\mathcal{O}}\oplus \mbox{Ker}((\pi_\mu\circ (-J))_{|T_\mu\M}).
\end{equation}
In particular, the restriction of $\pi_\mu\circ (-J)$ to $\overline{T_\mu\mathcal{O}}$ is injective.
\end{lemma}
\proof{}
We introduce the following notation: given any map $P$ defined on $L^2(\mu)$, let $P'$ denote its restriction to the closed subspace $T_\mu\M=\mbox{Im}(\pi_\mu)$. 

Set $P:=\pi_\mu\circ (-J)$. Then $\mbox{Im}(P')\subseteq \mbox{Im}(\pi_\mu)$ so we can think of $P'$ as a map $P':T_\mu\M\rightarrow T_\mu\M$, yielding a decomposition $T_\mu\M=\mbox{Im}(P')\oplus \mbox{Ker} (P'^*)$. It is simple to check that 
$$P'^*=(\pi_\mu \circ P^*)'=(\pi_\mu\circ J\circ\pi_\mu)'=(\pi_\mu\circ J)'.$$ 
In particular, $P'$ is anti-selfadjoint. This implies that $\mbox{Ker} (P'^*)=\mbox{Ker}(P')$ so 
\begin{equation*}
T_\mu\M=\overline{\mbox{Im}(P')}\oplus \mbox{Ker} (P')=\overline{T_\mu\mathcal{O}}\oplus \mbox{Ker}(P').
\end{equation*}
\endproof

It follows from Lemma \ref{l:orthogonalspaces} that any element $X\in T_\mu\mathcal{O}$ can be uniquely written as $X=\pi_\mu(-J\tilde{X})$, for some $\tilde{X}\in (\mbox{Ker}((\pi_\mu\circ (-J))_{|T_\mu\M}))^\perp$.

\begin{remark} 
Let $F$ be a differentiable function on $\M$ and let $X_F$ be the corresponding Hamiltonian vector field. It follows from Lemma \ref{l:orthogonalspaces} that $X_F$ depends only on the component of $\nabla F$ tangent to $\overline{T_\mu\mathcal{O}}$. Intuitively, this suggests that the corresponding Hamiltonian flow should depend only on the restriction of $F$ to $\mathcal{O}$.
\end{remark}

\begin{example}
Let $P'$ denote the restriction of $\pi_\mu\circ (-J)$ to $T_\mu\M$. It follows from Example \ref{e:hamfieldsvsdecomp} that if $\mu=\delta_x$ then the map $P'$ is an isomorphism of $T_\mu\M$. If instead $\mu=\rho\mathcal{L}$ for some $\rho>0$ then $P'$ is neither injective nor surjective.
\end{example}

\begin{lemma}\label{l:Omegaflat}
For any $\mu\in \M$ the map
\begin{equation*}
\Omega_\mu^\flat:T_\mu\mathcal{O}\rightarrow (T_\mu\mathcal{O})^*,\ \ X\mapsto\Omega_\mu(X,\cdot)
\end{equation*}
is a (non-continuous) isomorphism.
\end{lemma}
\proof{}
Given any $X,Y\in T_\mu\mathcal{O}$, we can write $X=\pi_\mu(-J\tilde{X})$, $Y=\pi_\mu(-J\tilde{Y})$ for some $\tilde{X},\tilde{Y}\in T_\mu\M$. Then
\begin{equation*}
\Omega_\mu(X,Y)=\hat{\Omega}_\mu(\tilde{X},\tilde{Y})=\hat{G}_\mu(\tilde{X}, -J\tilde{Y})=\hat{G}_\mu(\tilde{X}, Y)\leq\|\tilde{X}\|_\mu\cdot\|Y\|_\mu.
\end{equation*}
This proves that $\Omega_\mu^\flat$ is well-defined, \textit{i.e.} that $\Omega_\mu^\flat(X)\in (T_\mu\mathcal{O})^*$.

Assume $\Omega_\mu(X,Y)=0$ for all $Y\in T_\mu\mathcal{O}$. Then, as above, $\hat{G}_\mu(\tilde{X},Y)=0$ for all $Y\in T_\mu\mathcal{O}$ so $\tilde{X}\in(T_\mu\mathcal{O})^\perp$. It follows from Lemma \ref{l:orthogonalspaces} that $X=0$ so $\Omega_\mu^\flat$ is injective.

Now choose $\Lambda\in (T_\mu\mathcal{O})^*$. Since $(T_\mu\mathcal{O})^*=(\overline{T_\mu\mathcal{O}})^*$, there exists $\tilde{X}\in\overline{T_\mu\mathcal{O}}$ such that $\Lambda=\hat{G}_\mu(\tilde{X},\cdot)$. Set $X:=\pi_\mu(-J\tilde{X})\in T_\mu\mathcal{O}$. Then
\begin{equation*}
\Lambda(Y)=\hat{G}_\mu(\tilde{X},\pi_\mu(-J\tilde{Y}))=\hat{G}_\mu(\tilde{X},-J\tilde{Y})=\hat{\Omega}_\mu(\tilde{X},\tilde{Y})=\Omega_\mu(X,Y).
\end{equation*} 
This shows that $\Omega_\mu^\flat$ is surjective.
\endproof

\begin{remark} Lemma \ref{l:orthogonalspaces} shows that the map $\pi_\mu\circ(-J):\overline{T_\mu\mathcal{O}}\rightarrow T_\mu\mathcal{O}$ is invertible. Its inverse is related to the notions introduced in Definition \ref{def:compatible} as follows. We can use the isomorphism $(\overline{T_\mu\mathcal{O}})^*\simeq\overline{T_\mu\mathcal{O}}$ induced by $\hat{G}_\mu$ to define a (non-continuous) isomorphism $A:T_\mu\mathcal{O}\rightarrow\overline{T_\mu\mathcal{O}}$ such that $\Omega_\mu(\cdot,\cdot)=\hat{G}_\mu(A\cdot,\cdot)$ on $T_\mu\mathcal{O}$. Notice that 
\begin{equation*}
\hat{G}_\mu(AX,Y)=\Omega_\mu(X,Y)=\hat{\Omega}_\mu(\tilde{X},\tilde{Y})=\hat{G}_\mu(\tilde{X},Y),
\end{equation*}
where we use the notation introduced in the proof of Lemma \ref{l:Omegaflat}. This shows that $AX=\tilde{X}$ so $A^{-1}=\pi_\mu\circ (-J)$.

If $\mu$ is a Dirac measure it is clearly the case that $G_\mu$ and $\Omega_\mu$ are a compatible pair in the sense of Definition \ref{def:compatible}. This amounts to stating that $(\pi_\mu\circ (-J))^2=-Id$ on $T_\mu\mathcal{O}$. It is not clear if this is true in general, even under the additional assumption that $T_\mu\mathcal{O}=\overline{T_\mu\mathcal{O}}$.
\end{remark}


\section{The symplectic foliation as a Poisson structure} \label{s:foliationispoisson}
Most naturally occurring symplectic foliations owe their existence to an underlying \textit{Poisson structure}. The symplectic foliation described in Section \ref{ss:symplecticfoliation} is no exception. The existence of a Poisson structure on a certain space of distributions was pointed out in \cite{marsdenweinstein:maxwellvlasov}. It stems from the fact that the symplectic structure on $\Rtwod$ adds new structure into the framework of Section \ref{ss:embeddingbis}. The goal of this section is to show that, reduced to $\M$, this Poisson structure coincides with the symplectic structure $\Omega$ defined in Section \ref{ss:symplecticfoliation}. We achieve this in Section \ref{ss:poissonissymplectic} after briefly reviewing the relevant notions. We refer to \cite{marsdenratiu:book} for many details.

\subsection{Review of Poisson geometry}\label{ss:poisson}
Recall from Section \ref{ss:hamdiffs} that any symplectic structure $\omega$ on a manifold $M$ induces a Lie bracket on the space of functions $\cm$. Using the Liebniz rule for the derivative of the product of two functions, we see that the corresponding operators $\{\cdot,h\}$ have the following property:
$$\{fg,h\}=d(fg)(X_h)=df(X_h)g+dg(X_h)f=\{f,h\}g+\{g,h\}f.$$
This leads to the following natural ``weakening" of Symplectic Geometry.

\begin{definition}
Let $M$ be a smooth manifold. A \textit{Poisson structure} on $M$ is a Lie bracket $\{\cdot,\cdot\}$ on $\cm$ such that each operator $\{\cdot,h\}$ is a \textit{derivation} on functions, \textit{i.e.}
$$\{fg,h\}=\{f,h\}g+\{g,h\}f.$$
A \textit{Poisson manifold} is a manifold endowed with a Poisson structure.
\end{definition}

On any finite-dimensional manifold it is known that the space of derivations on functions is isomorphic to the space of vector fields. Thus on any Poisson manifold any function $h$ defines a vector field which we denote $X_h$: it is uniquely defined by the property that
$$df(X_h)=\{f,h\},\ \ \forall f\in\cm.$$
We call $X_h$ the \textit{Hamiltonian vector field} defined by $h$. As in Section \ref{ss:hamdiffs}, this process defines a map
\begin{equation}\label{eq:poissonmapbis}
\cm\rightarrow\xm,\ \ f\mapsto X_f.
\end{equation}
The kernel of this map includes the space of constant functions, but in general it will be larger. We call these the \textit{Casimir functions} of the Poisson manifold. Its image defines the space $\hamm$ of \textit{Hamiltonian vector fields}.  Lemma \ref{l:mapislie} applies with the same proof to show that the map of Equation \ref{eq:poissonmapbis} is a Lie algebra homomorphism (up to sign).

At each point $x\in M$, the set of Hamiltonian vector fields evaluated at that point define a subspace of $T_xM$. The union of such subspaces is known as the \textit{characteristic distribution} of the Poisson manifold. This distribution is \textit{integrable} in the sense that $M$ admits a smooth foliation such that each subspace is the tangent space of the corresponding leaf. In particular each leaf has a well-defined dimension, but this dimension may vary from leaf to leaf. Each leaf has a symplectic structure defined by setting
\begin{equation}\label{eq:symplecticonleaf}
\omega(X_f,X_g):=\{f,g\}.
\end{equation}
\begin{remark}
Notice that for any given Hamiltonian vector field $X_f$, the corresponding function $f$ is well-defined only up to Casimir functions. It is however simple to check that $\omega$ is a well-defined 2-form on each leaf, \textit{i.e.} it is independent of the particular choices made for $f$ and $g$. It is also non-degenerate. The fact that $\omega$ is closed follows from the Jacobi identity for $\{\cdot,\cdot\}$.
\end{remark}

\begin{remark} Notice that the definition of a Poisson manifold does not include a metric. Thus there is in general no intrinsic way to extend $\omega$ to a 2-form on $M$. 
\end{remark}

The following result is standard.

\begin{prop} \label{prop:smoothleaves}
Any Poisson manifold admits a symplectic foliation, of varying rank. Each leaf is preserved by the flow of any Hamiltonian vector field. Any Casimir function restricts to a constant along any leaf of the foliation.
\end{prop}

Poisson manifolds are of interest in Mechanics because they provide the following generalization of the standard symplectic notion of Hamiltonian flows.
\begin{definition}\label{def:stdhamflow}
A \textit{Hamiltonian flow} on $M$ is a solution of the equation $d/dt (x_t)=X_f(x_t)$, for some function $f$ on $M$.
\end{definition}
It follows from Proposition \ref{prop:smoothleaves} that if the initial data belongs to a specific leaf, then the corresponding Hamiltonian flow is completely contained within that leaf. It is simple to check that if $x_t$ is Hamiltonian then $f$ is constant along $x_t$.

The theory of Lie algebras provides one of the primary classes of examples of Poisson manifolds. To explain this we introduce the following notation, once again restricting our attention to the finite-dimensional case. Let $V$ be a finite-dimensional vector space, whose generic element will be denoted $v$. Let $V^*$ be its dual, with generic element $\phi$. Let $V^{**}$ be the bidual space, defined as the space of linear maps $V^*\rightarrow\R$. We will think of this as a subspace of the space of smooth maps on $V^*$, with generic element $f=f(\phi)$. We can identify $V$ with $V^{**}$ via the map
\begin{equation}\label{eq:bidual}
V\rightarrow V^{**}, \ \ v\mapsto f_v \mbox{ where }f_v(\phi):=\phi(v).
\end{equation}
Now assume $V$ is a Lie algebra. We will write $V=\mathfrak{g}$. Consider the vector space $\mathfrak{g}^*$ dual to $\mathfrak{g}$. We want to show that the Lie algebra structure on $\mathfrak{g}$ induces a natural Poisson structure on $\mathfrak{g}^*$. Let $f$ be a smooth function on $\mathfrak{g}^*$. Its linearization at $\phi$ is an element of the bidual: $df_{|\phi}\in \mathfrak{g}^{**}$. It thus corresponds via the map of Equation \ref{eq:bidual} to an element $\delta f/\delta\phi_{|\phi}\in\mathfrak{g}$. We can now define a Lie bracket on $\mathfrak{g}^*$ by setting:
\begin{equation}\label{eq:dualliepoisson}
\{f,g\}(\phi):=\phi([\delta f/\delta\phi_{|\phi},\delta g/\delta\phi_{|\phi}]),
\end{equation}
where $[\cdot,\cdot]$ denotes the Lie bracket on $\mathfrak{g}$. One can show that this operation satisfies the Jacobi identity and defines a Poisson structure on $\mathfrak{g}^*$.

\begin{example} \label{e:linearpoisson}
Assume $f$ is a linear function on $\mathfrak{g}^*$, $f=f_v$ (as in Equation \ref{eq:bidual}). Then $\delta f/\delta\phi\equiv v$, so $\{f_v,f_w\}(\phi)=\phi([v,w])$.
\end{example}

We now want to characterize the Hamiltonian vector fields and symplectic leaves of $\mathfrak{g}^*$. Unsurprisingly, this is best done in terms of Lie algebra theory. Every finite-dimensional Lie algebra is the Lie algebra of a (unique connected and simply connected) Lie group $G$. Recall from Section \ref{ss:groupactions} the adjoint representation of $G$ on $\mathfrak{g}$,
$$G\rightarrow Aut(\mathfrak{g}),\ \ g\mapsto Ad_g.$$
Differentiating this defines the \textit{adjoint representation} of $\mathfrak{g}$ on $\mathfrak{g}$,
\begin{equation}\label{eq:adjointrep}
ad:\mathfrak{g}\rightarrow End(\mathfrak{g}),\ \ v=d/dt (g_t)_{|t=0}\mapsto ad_v:=d/dt(Ad_{g_t})_{|t=0}.
\end{equation}
It follows from Lemma \ref{l:bracketviaad} that $ad_v(w)=[v,w]$.

By duality we obtain the \textit{coadjoint representation} of $G$ on $\mathfrak{g}^*$,
$$G\rightarrow Aut(\mathfrak{g}^*),\ \ g\mapsto (Ad_{g^{-1}})^*.$$
Notice that once again we have used inversion to ensure that this remains a group homomorphism, cf. Remark \ref{r:rightactions}. We can differentiate this to obtain the \textit{coadjoint representation} of $\mathfrak{g}$ on $\mathfrak{g}^*$, which can also be written in terms of the duals of the maps in Equation \ref{eq:adjointrep}:
\begin{equation}
ad^*:\mathfrak{g}\rightarrow End(\mathfrak{g}^*),\ \ v\mapsto (-ad_v)^*.
\end{equation}
The following result is standard, cf. \textit{e.g.} \cite{marsdenratiu:book} Formula 10.7.4.
\begin{lemma} \label{l:leaforbit}
The Hamiltonian vector field corresponding to a smooth function $f$ on $\mathfrak{g}^*$ is
$$X_f(\phi):=(-ad_{\delta f/\delta\phi_{|\phi}})^*(\phi).$$
Thus the leaves of the symplectic foliation of $\mathfrak{g}^*$ are the orbits of the coadjoint representation.
\end{lemma}


\subsection{The symplectic foliation of $\M$, revisited}\label{ss:poissonissymplectic}
Following \cite{marsdenweinstein:maxwellvlasov} we now apply the ideas of Section \ref{ss:poisson} to the case where $\mathfrak{g}$ is the Lie algebra of $\hamrtwod$. Since this is an infinite-dimensional algebra, the following discussion will be purely formal.

We saw in Remark \ref{r:normalization} that $\mathfrak{g}$ can be identified with the space of compactly-supported functions:
\begin{equation}\label{eq:lieisomorphism}
C^\infty_c(\Rtwod)\simeq\chamfields(\Rtwod), \ \ f\mapsto X_f.
\end{equation}
Its dual is then the distribution space $\N$. Section \ref{ss:poisson} suggests that $\N$ has a canonical Poisson structure, defined as in Equation \ref{eq:dualliepoisson}. We can identify the Poisson bracket, Hamiltonian vector fields and symplectic leaves on $\N$ very explicitly, as follows.

For simplicity let us restrict our attention to the linear functions on $\N$ defined by functions
$f\in C^\infty_c$ as follows:
\begin{equation}\label{eq:deflinear}
F_f:\N\rightarrow\R,\ \ F_f(\mu):=\langle\mu,f\rangle.
\end{equation}
Example \ref{e:linearpoisson} shows that the Poisson bracket of two such functions $F_f$ and $F_g$ can be written in terms of the Lie bracket on $C^\infty_c$:
\begin{equation}\label{eq:poissonmeasures}
\{F_f,F_g\}_{\N}(\mu)=\langle\mu,\{f,g\}_{\Rtwod}\rangle=\langle\mu,\omega(X_f,X_g)\rangle.
\end{equation}
Lemma \ref{l:leaforbit} gives an explicit formula for the corresponding Hamiltonian vector fields $X_{F_f}$: at $\mu\in\N$, $X_{F_f}(\mu)\in T_\mu\N=\N$ is given by
\begin{align*}
\langle X_{F_f}(\mu),g\rangle &=\langle(-ad_f)^*(\mu),g\rangle=\langle\mu,-ad_f(g)\rangle=\langle\mu,-\{f,g\}_{\Rtwod}\rangle=\langle\mu,dg(X_f)\rangle\\
&=-\langle div_\mu(X_f),g\rangle.
\end{align*}
In other words, $X_{F_f}(\mu)=-div_\mu(X_f)$.

Lemma \ref{l:leaforbit} also shows that the leaves of the symplectic foliation are the orbits of the coadjoint representation of $\hamrtwod$ on $\N$. Let us identify the coadjoint representation explicitly. Recall from Lemma \ref{l:adjointispushforward} that the adjoint representation of $\hamrtwod$ on $\chamfields$ is the push-forward operation. Lemma \ref{l:hamfieldproperties} shows that, under the isomorphism of Equation \ref{eq:lieisomorphism}, push-forward becomes composition. Thus the adjoint representation of $\hamrtwod$ on $\chamfields$ corresponds to the following representation of $\hamrtwod$ on $\ccinf(\R^{2d})$:
\begin{equation}
Ad:\hamrtwod\rightarrow Aut(C^\infty_c(\Rtwod)),\ \ Ad_\phi(f):=f\circ\phi^{-1}.
\end{equation}
The following calculation then shows that the coadjoint representation of $\hamrtwod$ on $\N$ is simply the natural action of $\hamrtwod$ introduced in Section \ref{ss:embeddingbis}:
$$\langle (Ad_{\phi^{-1}})^*(\mu),f\rangle=\langle\mu,Ad_{\phi^{-1}}(f)\rangle=\langle\mu,f\circ\phi\rangle=\langle\phi\cdot\mu, f\rangle.$$
The symplectic structure on each leaf is given by Equation \ref{eq:symplecticonleaf}:
\begin{equation}\label{eq:symplecticonN}
\omega_\mu(-div_\mu(X_f),-div_\mu(X_g)):=\{F_f,F_g\}_{\N}(\mu)=\langle\mu,\omega(X_f,X_g)\rangle.
\end{equation}

\begin{remark} \label{r:generalhamfields}
Notice that Poisson brackets and Hamiltonian vector fields are of first order with respect to the functions involved. We can use this fact to reduce the study of general functions $F:\N\rightarrow\R$ to the study of linear functions on $\N$, as presented above. For example if $\nabla_\mu F=\nabla_\mu F_f$, for some linear $F_f$ as above, then $X_F(\mu)=X_{F_f}(\mu)$.
\end{remark}
 
Let us now restrict our attention to $\M\subset\N$. We want to show that the data defined by the Poisson structure on $\N$ restricts to the objects defined in Section \ref{ss:symplecticfoliation}. Firstly, $\M$ is $\hamrtwod$-invariant and the action of $\hamrtwod$ on $\N$ restricts to the standard push-forward action on $\M$. This shows that the leaves defined above, passing through $\M$, coincide with the $\G$-orbits of Section \ref{ss:symplecticfoliation}. Now recall from Section \ref{ss:embeddingbis} that, given $\mu\in\M$, the operator $-\divmu$ is the natural isomorphism relating the tangent planes of Definition \ref{def:tangentspaces} to the tangent planes of $\M\subset\N$. Equation \ref{eq:symplecticonN} can thus be re-written as
$$\omega_\mu(\pi_\mu(X_f),\pi_\mu(X_g)):=\int_{\Rtwod}\omega(X_f,X_g)\,d\mu,$$
showing that the symplectic structure defined this way coincides with the symplectic form $\Omega_\mu$ defined in Equation \ref{eq:symplectic}. To conclude, we want to show that the Hamiltonian vector fields introduced in Definition \ref{def:hamflow} formally coincide with the Hamiltonian vector fields of the restricted Poisson structure. Let $F:\M\rightarrow\R$ be a  differentiable function on $\M$. Fix $\mu\in\M$. Up to $L^2_\mu$-closure, we can assume that $\nabla_\mu F=\nabla f$, for some $f\in\ccinf(\R^{2d})$. Example \ref{e:linearisdifferentiable} shows that $\nabla f=\nabla_\mu F_f$, where $F_f$ is the linear function defined in Equation \ref{eq:deflinear}. Using Remark \ref{r:generalhamfields}, the Hamiltonian vector of $F$ at $\mu$ defined by the Poisson structure is thus $X_F(\mu)=X_{F_f}(\mu)=-\divmu(X_f)$. In terms of the tangent space $T_\mu\M$, we can write this as
\begin{equation}
X_F(\mu)=\pi_\mu(X_f)=\pi_\mu(-J\nabla f)=\pi_\mu(-J\nabla_\mu F).
\end{equation}
It thus coincides with the vector field given in Definition \ref{def:hamflow}.

\begin{remark}
The identification of $\N$ with the dual Lie algebra of $\hamrtwod$ relied on the normalization introduced in Remark \ref{r:normalization}. In turn, this was based on our choice to restrict our attention to diffeomorphisms with compact support. In some situations one might want to relax this assumption. This would generally mean losing the possibility of a normalization so Equation \ref{eq:hamsequence} would leave us only with an identification $\hamfields\simeq C^\infty(M)/\R$. Dualizing this space would then, roughly speaking, yield the space of measures of integral zero: we would thus get a Poisson structure on this space but not on $\M$. However this issue is purely technical and can be avoided by changing Lie group, as follows.

Consider the group $G$ of diffeomorphisms on $\Rtwod\times\R$ preserving the \textit{contact form} $dz-y^idx^i$. It can be shown that its Lie algebra is isomorphic to the space of functions on $\Rtwod\times\R$ which are constant with respect to the new variable $z$: it is thus isomorphic to the space of functions on $\Rtwod$, so the dual Lie algebra is, roughly, the space of measures on $\Rtwod$; in particular, it contains $\M$ as a subset. This group has a one-dimensional center $Z\simeq\R$, defined by translations with respect to $z$. The center acts trivially in the adjoint and coadjoint representations, so the coadjoint representation reduces to a representation of the group $G/Z$, which one can show to be isomorphic to the group of Hamiltonian diffeomorphisms of $\Rtwod$. The coadjoint representation of $G$ reduces to the standard push-forward action of Hamiltonian diffeomorphisms, and the theory can now proceed as before.
\end{remark}


\appendix

\section{Review of relevant notions of Differential Geometry} \label{s:geometryreview}
The goal of the first two sections of this appendix is to summarize standard facts concerning Lie groups and calculus on finite-dimensional manifolds, thus laying out the terminology, notation and conventions which we use throughout this paper. We refer to \cite{kobayashinomizu:book} and \cite{marsdenratiu:book} for details. The third section introduces the notion of invariant cohomology. The point of view adopted here might be new. It provides useful analogies for the notion of ``pseudo forms" introduced in Section \ref{ss:forms}. The fourth section provides some basic facts concerning the infinite-dimensional Lie groups relevant to this paper. 

\subsection{Calculus of vector fields and differential forms} \label{ss:calculus}
Let $M$ be a connected differentiable manifold of dimension $D$, not necessarily compact. Let $\diffm$ denote the group of diffeomorphisms of $M$. Let $\cm$ denote the space of smooth functions on $M$. Let $TM$ denote the tangent bundle of $M$ and $\xm$ the corresponding space of sections, \textit{i.e.} the space of smooth vector fields. Let $T^*M$ denote the cotangent bundle of $M$. To simplify notation, $\Lambda^kM$ will denote both the bundle of k-forms on $M$ and the space of its sections, \textit{i.e} the space of smooth k-forms on $M$. Notice that $\Lambda^0(M)=\cm$ and $\Lambda^1M=T^*M$ (or the space of smooth 1-forms).

Let $\phi\in\diffm$. Taking its differential yields linear maps
\begin{equation}\label{eq:nablaphi}
\nabla\phi:T_xM\rightarrow T_{\phi(x)}M,\ \ v\mapsto \nabla\phi\cdot v,
\end{equation}

thus a bundle map which we denote $\nabla\phi: TM\rightarrow TM$. We will call $\nabla\phi$ the \textit{lift} of $\phi$ to $TM$.

By duality we obtain linear maps
$$(\nabla\phi)^*:T_{\phi(x)}^*M\rightarrow T_x^*M,\ \ \alpha\mapsto \alpha\circ\nabla\phi,$$
and more generally k-multilinear maps
\begin{equation}\label{eq:dualnablaphi}
(\nabla\phi)^*:\Lambda^k_{\phi(x)}M\rightarrow \Lambda^k_xM,\ \ \alpha\mapsto\alpha(\nabla\phi\,\cdot,\dots,\nabla\phi\,\cdot).
\end{equation}
This defines bundle maps $(\nabla\phi)^*:\Lambda^kM\rightarrow\Lambda^kM$ which we call the \textit{lift} of $\phi$ to $\Lambda^kM$.

\begin{remark}
Notice the different behaviour under composition of diffeomorphisms: $\nabla(\phi\circ\psi)=\nabla\phi\circ\nabla\psi$ while $(\nabla(\phi\circ\psi))^*=(\nabla\psi)^*\circ(\nabla\phi)^*$. We will take this into account and generalize it in Section \ref{ss:groupactions} via the notion of left versus right group actions.
\end{remark}

We can of course apply these lifted maps to sections of the corresponding bundles. In doing so one needs to ensure that the correct relationship between $T_xM$ and $T_{\phi(x)}M$ is maintained; we emphasize this with a change of notation, as follows.

The \textit{push-forward} operation on vector fields is defined by
\begin{equation}\label{eq:pushforward}
\phi_*:\xm\rightarrow \xm, \ \ \phi_*X:=(\nabla\phi\cdot X)\circ\phi^{-1}.
\end{equation}
The corresponding operation on k-forms is the \textit{pull-back}, defined by
\begin{equation}\label{eq:pullback}
\phi^*:\Lambda^k(M)\rightarrow\Lambda^k(M), \ \ \phi^*\alpha:=((\nabla\phi)^*\alpha)\circ\phi.
\end{equation}
\begin{definition}\label{def:liealgebra}
Let $V$ be a vector space. A bilinear anti-symmetric operation
$$V\times V\rightarrow V, \ \ (v,w)\mapsto [v,w]$$
is a \textit{Lie bracket} if it satisfies the \textit{Jacobi identity}
$$[u,[v,w]]+[v,[w,u]]+[w,[u,v]]=0.$$
A \textit{Lie algebra} is a vector space endowed with a Lie bracket.
\end{definition}

The space of smooth vector fields has a natural Lie bracket. Given two vector fields $X$, $Y$ on $M$, we define $[X,Y]$ in local coordinates as follows:
$$[X,Y]:=\nabla Y\cdot X-\nabla X\cdot Y.$$
It is simple to show that this operation indeed satisfies the Jacobi identity. Let $\phi_t$ denote the \textit{flow} of $X$ on $\RD$, \textit{i.e.} the 1-parameter group of diffeomorphisms obtained by integrating $X$ as follows:
\begin{equation} \label{eq:flow}
d/dt(\phi_t(x))=X(\phi_t(x)),\ \ \phi_0(x)=x.
\end{equation}
It is then simple to check that
\begin{equation} \label{eq:liebracket}
[X,Y]=-d/dt(\phi_{t*}Y)_{|t=0}=d/dt(\phi_{-t*}Y)_{|t=0}=d/dt((\phi_t^{-1})_*Y)_{|t=0}.
\end{equation}
Equation \ref{eq:liebracket} gives a coordinate-free expression for the Lie bracket. It also suggests an analogous operation for more general tensor fields. We will restrict our attention to the case of differential forms.

Let $\alpha$ be a smooth k-form on $M$. Let $X$, $\phi_t$ be as above. We define the \textit{Lie derivative} of $\alpha$ in the direction of $X$ to be the k-form defined as follows:
\begin{equation} \label{eq:liederivative}
\mathcal{L}_X\alpha:=d/dt(\phi_t^*\alpha)_{|t=0}.
\end{equation}
The fact that $t\mapsto \phi_t$ is a homomorphism leads to the fact that $d/dt(\phi_t^*\alpha)_{|t=t_0}=\phi_{t_0}^*(\mathcal{L}_X\alpha)$. Thus $\mathcal{L}_X\alpha\equiv 0$ if and only if $\phi_t^*\alpha\equiv \alpha$, \textit{i.e.} $\phi_t$ preserves $\alpha$. This can be generalized to time-dependent vector fields as follows.
\begin{lemma} \label{l:preservedform}
Let $X_t$ be a t-dependent vector field on $M$. Let $\phi_t=\phi_t(x)$ be its flow, defined by
\begin{equation} \label{eq:tflow}
d/dt(\phi_t(x))=X_t(\phi_t(x)),\ \ \phi_0(x)=x.
\end{equation}
Let $\alpha$ be a k-form on $M$. Then $d/dt(\phi_t^*\alpha)_{|t_0}=\phi_{t_0}^*(\mathcal{L}_{X_{t_0}}\alpha)$. In particular,
$\phi_t^*\alpha\equiv\alpha$ iff $\mathcal{L}_{X_t}\alpha\equiv 0$.
\end{lemma}
\proof{} For any fixed $s$, let $\psi_t^s$ be the flow of $X_s$, \textit{i.e.}
$$d/dt(\psi_t^s(x))=X_s(\psi_t^s(x)),\ \ \psi_0^s(x)=x.$$
Then $\psi_t^{t_0}\circ\phi_{t_0}(x)$ satisfies
$$d/dt(\psi_t^{t_0}\circ\phi_{t_0}(x))_{|t=0}=X_{t_0}(\phi_{t_0}(x)),\ \ \psi_0^{t_0}\circ\phi_{t_0}(x)=\phi_{t_0}(x)$$
so $\psi_t^{t_0}\circ\phi_{t_0}(x)$ at $t=0$ and $\phi_t$ at $t=t_0$ coincide up to first order, showing that
$$d/dt(\phi_t^*\alpha)_{|t=t_0}=d/dt((\psi_t^{t_0}\circ\phi_{t_0})^*\alpha)_{|t=0}=\phi_{t_0}^*(d/dt((\psi_t^{t_0})^*\alpha)_{|t=0})=\phi_{t_0}^*(\mathcal{L}_{X_{t_0}}\alpha).$$
\endproof

Notice that if we define $\phi^*Y:=(\phi^{-1})_*Y$ and we define $\mathcal{L}_XY:=d/dt(\phi_t^*Y)_{|t=0}$, then Equation \ref{eq:liebracket} shows that $\mathcal{L}_XY=[X,Y]$.

\begin{remark} Various formulae relate the above operations, leading to quick proofs of useful facts. For example, the fact
\begin{equation} \label{eq:liederivativebracket}
\mathcal{L}_{[X,Y]}\alpha=\mathcal{L}_X(\mathcal{L}_Y\alpha)-\mathcal{L}_Y(\mathcal{L}_X\alpha)
\end{equation}
shows that if the flows of $X$ and $Y$ preserve $\alpha$ then so does the flow of $[X,Y]$. Also,
\begin{equation} \label{eq:liederivativepullback}
\phi^*\mathcal{L}_X\alpha=\mathcal{L}_{\phi^*X}\phi^*\alpha.
\end{equation}
\end{remark}

\begin{remark} \label{r:dirderivative}
Notice that $\mathcal{L}_XY$ is not a ``proper" directional derivative in the sense that it is of first order also in the vector field $X$. In general the same is true for the Lie derivative of any tensor. The case of 0-forms, \textit{i.e.} functions, is an exception. In this case $\mathcal{L}_Xf=df(X)$ is of order zero in $X$ and coincides with the usual notion of directional derivative. We will often simplify the notation by writing it as $Xf$.
\end{remark}

We now want to introduce the \textit{exterior differentiation} operator on smooth forms. Let $\alpha$ be a k-form on $M$. Fix any point $x\in M$ and tangent vectors $X_0,\dots,X_k\in T_xM$. Choose any extension of each $X_j$ to a global vector field which we will continue to denote $X_j$. Then, at $x$,
\begin{align}\label{eq:usuald}
d\alpha(X_0,\dots,X_k)&:=\sum_{j=0}^k (-1)^jX_j\alpha(X_0,\dots,\hat{X_j},\dots,X_k)\\
&\quad + \sum_{j<l}(-1)^{j+l}\alpha([X_j,X_l],X_0,\dots,\hat{X_j},\dots,\hat{X_l},\dots,X_k)\nonumber
\end{align} 
where on the right hand side the superscript $\hat{}$ denotes an omitted term and we adopt the notation for directional derivatives introduced in Remark \ref{r:dirderivative}.

\begin{lemma} \label{l:d}
$d\alpha$ is a well-defined (k+1)-form, \textit{i.e.} at any point $x\in M$ it is independent of the choice of the extension. Exterior differentiation defines a first-order linear operator
\begin{equation}
d:\Lambda^kM\rightarrow \Lambda^{k+1}M
\end{equation}
satisfying $d\circ d=0$.
\end{lemma}

\begin{remark}\label{r:dof1form} It is not clear from the above definition that $d\alpha$ is \textit{tensorial} in $X_0,\dots,X_k$, \textit{i.e.} that it is independent of the choice of extensions. The point is that cancelling occurs to eliminate the first derivatives of $X_j$ which appear in Equation \ref{eq:usuald}. This is the main content of Lemma \ref{l:d}, which is proved by showing that Equation \ref{eq:usuald} is equivalent to the usual, local-coordinate, definition of $d\alpha$. For example, let $\alpha=\sum_{i=1}^D \alpha_i(x)dx^i$ be a smooth 1-form on $\RD$. Then 
$d\alpha=\sum_{j<i}\left(\frac{\partial\alpha_i}{\partial x^j}-\frac{\partial\alpha_j}{\partial x^i}\right) dx^j\wedge dx^i$. If we identify $\alpha$ with the vector field $x \rightarrow (\alpha_1(x), \cdots, \alpha_D(x))^T$ then  $d\alpha(X,Y)=\langle (\nabla \alpha- \nabla \alpha^T)  X, Y \rangle.$ 
\end{remark}

Given a k-form $\alpha$ and a vector field $X$, let $i_X\alpha$ denote the (k-1)-form $\alpha(X,\cdot,\dots,\cdot)$ obtained by \textit{contraction}. Then the Lie derivative and exterior differentiation are related by \textit{Cartan's formula}:
\begin{equation}\label{eq:cartan}
\mathcal{L}_X\alpha=d\,i_X\alpha+i_Xd\alpha.
\end{equation}


\subsection{Lie groups and group actions}\label{ss:groupactions}
Recall that a group $G$ is a \textit{Lie group} if it has the structure of a smooth manifold and group multiplication (respectively, inversion) defines a smooth map $G\times G\rightarrow G$ (respectively, $G\rightarrow G$). We denote by $e$ the identity element of $G$.

\begin{definition}
We say that $G$ has a \textit{left action} or \textit{acts on the left} or, more simply, \textit{acts} on a smooth manifold $M$ if there is a smooth map
$$G\times M\rightarrow M,\ \ (g,x)\mapsto g\cdot x$$
such that $g\cdot(h\cdot x)=(gh)\cdot x$. To simplify the notation we will often write $gx$ rather than $g\cdot x$. It is simple to see that if $G$ has a left action on $M$ then every $g\in G$ defines a diffeomorphism of $M$. More specifically, the action defines a group homomorphism $G\rightarrow\diffm$.

We say that $G$ has a \textit{right action} or \textit{acts on the right} on $M$ if the opposite composition rule holds: $g\cdot(h\cdot x)=hg\cdot x$. In this case it is standard to change the notation, writing $x\cdot g$ rather than $g\cdot x$: this makes the composition rule seem more natural but does not affect the substance of the definition, \textit{i.e.} the fact that the induced map $G\rightarrow\diffm$ is now a group anti-homomorphism.
\end{definition}

\begin{remark} \label{r:rightactions}
Notice that any left action induces a natural right action as follows: $x\cdot g:=g^{-1}\cdot x$. Conversely, any right action induces a natural left action: $g\cdot x:=x\cdot g^{-1}$.
\end{remark}

For any group action we can repeat the constructions of Equations \ref{eq:nablaphi} and \ref{eq:dualnablaphi}. For example a left action of $G$ on $M$ induces a \textit{lifted left action} of $G$ on $TM$ as follows:
$$G\times TM\rightarrow TM,\ \ g(x,v):=(gx,\nabla g\cdot v).$$
However, we need to apply the trick introduced in Remark \ref{r:rightactions} to obtain a coherent lifted action on $T^*M$ or $\Lambda^kM$. For example we can define a \textit{lifted left action} by setting
$$G\times \Lambda^kM\rightarrow\Lambda^kM,\ \ g(x,\alpha):=(gx, (\nabla g^{-1})^*\alpha)$$
or a \textit{lifted right action} by setting
$$G\times \Lambda^kM\rightarrow\Lambda^kM,\ \ g(x,\alpha):=(g^{-1}x, (\nabla g)^*\alpha).$$
We can also repeat the constructions of Equations \ref{eq:pushforward} and \ref{eq:pullback}. We thus find an induced action of $G$ on vector fields, defined by
\begin{equation}\label{eq:Gonfields}
G\times\fields (M)\rightarrow\fields (M),\ \ g\cdot X:=g_*X.
\end{equation}
Likewise, there is an induced action of $G$ on $k$-forms. On the other hand, with respect to Section \ref{ss:calculus} there now exists a new operation, as follows. Choose $v=d/dt(g_t)_{|t=0}\in T_eG$. For any $x\in M$ we can define the tangent vector $\tilde{v}(x):=d/dt(g_t\cdot x)_{|t=0}$. This defines a global vector field $\tilde{v}$ on $M$, called the \textit{fundamental vector field} associated to $v$. We have thus built a map $T_eG\rightarrow\fields(M)$.

Let us now specialize to the case $M=G$. Any Lie group $G$ admits two natural left actions on itself. Studying these actions leads to a deeper understanding of the geometry of Lie groups, thus of group actions. The first action is given by \textit{left translations}, as follows:
$$L:G\times G\rightarrow G, \ \ (g,h)\mapsto L_g(h):=gh.$$
Let $e\in G$ denote the identity element. Fix $v=d/dt(g_t)_{|t=0}\in T_eG$. The differential $\nabla L_g$ maps $T_eG$ to $T_gG$. We may thus define a global vector field $X_v$ on $G$ by setting $X_v(g):=\nabla L_g\cdot v=d/dt(gg_t)_{|t=0}$. This vector field has the property of being \textit{left-invariant} with respect to the action of $G$, \textit{i.e.} $L_{g*} X_v=X_v$. Viceversa, any left-invariant vector field arises this way.

\begin{lemma}\label{l:leftinvariant}
The set of left-invariant vector fields is a finite dimensional vector space isomorphic to $T_eG$. The Lie bracket of left-invariant vector fields is a left-invariant vector field.
\end{lemma}

It follows from Lemma \ref{l:leftinvariant} that $T_eG$ admits a natural operation $[v,w]$ such that $X_{[v,w]}=[X_v,X_w]$. It follows from the Jacobi identity on vector fields that $T_eG$ equipped with this structure is a Lie algebra: we call it the \textit{Lie algebra of $G$} and denote it by $\mathfrak{g}$.

\begin{remark}\label{r:rightinvariant}
Given any $v\in T_eG$, we have now defined two constructions of a global vector field on $G$ associated to $v$: the fundamental vector field $\tilde{v}$ and the left-invariant vector field $X_v$. It is simple to check that $\tilde{v}$ is invariant with respect to the \textit{right translations}
$$R:G\times G\rightarrow G, \ \ (g,h)\mapsto R_g(h):=hg.$$
This implies that the space of fundamental vector fields coincides with the space of right-invariant vector fields. The analogue of Lemma \ref{l:leftinvariant} holds for right-invariant fields and can be used to define a second Lie bracket on $T_eG$. It can be checked that this new bracket is simply the negative of the old one, \textit{i.e.} the two brackets differ only by sign.
\end{remark}

The second action of $G$ on itself is the \textit{adjoint action} defined by the \textit{inner automorphisms} $I_g(h):=ghg^{-1}$. Each of these fixes the identity and thus defines a map
\begin{equation}\label{eq:Gong}
Ad_g:=\nabla I_g: T_eG\rightarrow T_eG,
\end{equation}
\textit{i.e.} an automorphism of $T_eG$. In other words the adjoint action of $G$ on $G$ induces a left action of $G$ on $T_eG$ called the \textit{adjoint representation} of $G$. 

The adjoint representation of $G$ provides a useful way to calculate Lie brackets on $\mathfrak{g}$, as follows.
\begin{lemma}\label{l:bracketviaad}
Fix $v,w\in\mathfrak{g}$. Assume $v=d/dt (g_t)_{|t=0}$ for some $g_t\in G$. Then $[v,w]=d/dt(Ad_{g_t}(w))_{|t=0}$. 
\end{lemma}
\proof{}
Assume $w=d/ds (h_s)_{|s=0}$. By definition,
\begin{equation}
d/dt (Ad_{g_t}(w))_{|t=0} = d/dt\,d/ds(g_th_sg^{-1}_t)_{|t,s=0}.
\end{equation}
Notice that
\begin{equation}
X_v(g)=\nabla L_g(v)=d/dt (gg_t)_{|t=0}=d/dt (R_{g_t}(g))_{|t=0}.
\end{equation}
In particular this shows that, for $t=0$, $R_{g_t}$ coincides with the flow of $X_v$ up to first order. Thus 
\begin{align*}
[v,w] &= (\mathcal{L}_{X_v}X_w)_{|e} = d/dt ((R_{g_t})^*{X_w})_{|e;\,t=0} = d/dt ((R_{g_t^{-1}})_*{X_w})_{|e;\,t=0}\\
&= d/dt ((\nabla R_{g_t^{-1}})_{|g_t} {X_w}_{|g_t})_{|t=0}=d/dt((\nabla R_{g_t^{-1}})_{|g_t} d/ds (g_th_s)_{|s=0})_{t=0}\\
&= d/dt\,d/ds(g_th_sg_t^{-1})_{|s,t=0}.
\end{align*}
\endproof
\begin{remark} \label{r:tangentspacevsalgebra}
It is sometimes useful to distinguish the vector space $T_eG$ from the Lie algebra $\mathfrak{g}$, so as to distinguish between maps or constructions which involve the Lie bracket and those which do not. Our notation will sometimes reflect this.

For example, assume $G$ has a left action on $M$. One can then show that the construction of fundamental vector fields defines a Lie algebra anti-homomorphism $\mathfrak{g}\rightarrow\fields (M)$. Analogously one can show that every $Ad_g$ is an automorphism of $\mathfrak{g}$, \textit{i.e.} it preserves the Lie algebra structure: $Ad_g([v,w])=[Ad_g(v),Ad_g(w)]$.
\end{remark}

Let us now return to the general case of a Lie group acting on a manifold $M$. We can apply the above information on the geometry of Lie groups to develop a better understanding of the geometric aspects of the group action.

\begin{definition}
Assume $G$ acts on $M$. Fix $x\in M$. The \textit{orbit of $x$} in $M$ is the subset
$$\mathcal{O}_x:=\{g\cdot x:g\in G\}\subseteq M.$$
Notice that $\mathcal{O}_{gx}=\mathcal{O}_x$. The \textit{stabilizer of $x$} in $G$ is the closed subgroup
$$G_x:=\{g\in G: g\cdot x=x\}\subseteq G.$$
This is again a Lie group. We denote its Lie algebra $\mathfrak{g}_x$: it is a subalgebra of $\mathfrak{g}$. It is simple to check that $G_{gx}=I_g(G_x)=g\cdot G_x\cdot g^{-1}$ and that $\mathfrak{g}_{gx}=Ad_g(\mathfrak{g}_x)$.

We say that a subset $\mathcal{O}\subseteq M$ is an \textit{orbit} of the action if $\mathcal{O}=\mathcal{O}_x$, for some $x\in M$.
\end{definition}

The differential geometry of an orbit can be studied via the theory of \textit{homogeneous manifolds}, \textit{i.e.} manifolds obtained as quotients of Lie groups, as follows.

\begin{lemma}\label{l:orbits}
Let $G$ be a Lie group and $H$ be a closed subgroup. Then:
\begin{enumerate}
\item The quotient space $G/H$ has a natural smooth structure such that the projection $\pi:G\rightarrow G/H$ is a smooth map. The differential $\nabla\pi:T_eG\rightarrow T_{[e]}(G/H)$ is surjective with kernel $T_eH$ so it yields an identification $T_{[e]}(G/H)\simeq T_eG/T_eH$.
\item Left multiplication defines a natural action of $G$ on the manifold $G/H$ such that $\pi$ is $G$-equivariant. Choose $v\in T_eG$. Then the corresponding fundamental vector field on $G/H$, evaluated at $[e]$, coincides with $\nabla\pi(v)$.
\item Now assume $G$ acts on a manifold $M$. Choose $x\in M$ and set $H:=G_x$. Then the group action defines a smooth 1:1 equivariant immersion (not necessarily an embedding)
\begin{equation*}
j: G/H\rightarrow M,\ \ j([g]):=g\cdot x
\end{equation*}
with image $\mathcal{O}_x$. Using this immersion we can thus identify $\mathcal{O}_x$ with $G/G_x$.
\end{enumerate}
\end{lemma}

Lemma \ref{l:orbits} identifies the geometry of $\mathcal{O}_x$ with the geometry of the homogeneous space $G/G_x$. The choice of point $x$ plays an important role in this identification. However, given an orbit $\mathcal{O}$ of $G$ in $M$, the choice of $x\in\mathcal{O}$ is not canonical. Furthermore, given any other point $y\in\mathcal{O}$, the choice of $g\in G$ such that $y=gx$ is also not unique. The following lemma shows how things change under different such choices.

\begin{lemma}\label{l:changepoint}
Let $G$ be a Lie group and $H$ be a closed subgroup. Choose $g\in G$ and let $I_g$ denote the corresponding inner automorphism of $G$. It is simple to check that setting $[I_g]([k]):=[I_g(k)]$ yields a well-defined commutative diagram 
$$\begin{CD}\label{cd:Igdescends}
G   @>{I_g}>>   G\\
@V{\pi}VV   @V{\pi}VV\\
G/H @>{[I_g]}>> G/gHg^{-1}
\end{CD}$$
Now assume $G$ acts on $M$. Choose an orbit $\mathcal{O}$ and points $x,gx\in \mathcal{O}$. Set $H:=G_x$ so that we can identify $G/H\simeq\mathcal{O}_x$, $G/gHg^{-1}\simeq\mathcal{O}_{gx}$. In terms of these identifications, the map $[I_g]$ corresponds to the map $g:\mathcal{O}_x\rightarrow \mathcal{O}_{gx}$. Taking the differential of the maps in the above diagram thus leads to the commutative diagram
$$\begin{CD}\label{cd:nablagisadg}
T_eG   @>{Ad_g}>>   T_eG\\
@V{}VV   @V{}VV\\
T_x\mathcal{O} @>{\nabla g}>> T_{gx}\mathcal{O}
\end{CD}$$
where the vertical maps are those defined by the construction of fundamental vector fields.
\end{lemma}
\proof{}
Choose any $k\in G$. The identification $j:G/G_x\simeq \mathcal{O}_x$ implies $[k]\simeq k\cdot x$. Using the analogous identifications for $gx$ we find
\begin{equation*}
[I_g]([k])=[gkg^{-1}]\simeq gkg^{-1}\cdot gx=g(k\cdot x).
\end{equation*}
This proves that under these identifications $[I_g]$ corresponds to $g$. Now choose $v=d/dt(g_t)_{|t=0}\in T_eG$. Then, using the identification $G/G_x\simeq\mathcal{O}_x$,
\begin{align*}
\nabla\pi(v)&=d/dt(\pi(g_t))_{|t=0}=d/dt([g_t])_{|t=0}=d/dt([g_t\cdot e])_{|t=0}=d/dt(g_t\cdot[e])_{|t=0}\\
&\simeq d/dt(g_t\cdot x)_{|t=0}.
\end{align*}
This proves that $\nabla\pi(v)$ corresponds to $\tilde{v}(x)$, where $\tilde{v}$ is the fundamental vector field on $\mathcal{O}_x$ defined by $v$.
\endproof


\subsection{Cohomology and invariant cohomology}\label{ss:cohomology}
Let $M$ be a manifold. Recall that the \textit{de Rham cohomology} groups of $M$ are defined as the quotient spaces
\begin{equation*}
H^k(M;\R):=\frac{\mbox{Ker}(d:\Lambda^k(M)\rightarrow \Lambda^{k+1}(M))}{\mbox{Im}(d:\Lambda^{k-1}(M)\rightarrow \Lambda^k(M))}.
\end{equation*}
Given an action of a group $G$ on $M$, one can restrict one's attention to the space of k-forms which are invariant under the action of $G$. An analogous construction then leads to the definition of the \textit{invariant de Rham cohomology} groups of the pair $(M,G)$, cf. \textit{e.g.} \cite{bredon} Section V.12. For our purposes it is sufficient to only consider pairs of the form $(G,H)$, where $H$ is a closed subgroup of $G$ acting via right multiplication, \textit{i.e.} $R_h(g):=gh$. The construction is then as follows.

Consider the space of $H$-invariant k-forms on $G$,
\begin{equation*}
\Lambda^k(G^H):=\{\alpha\in\Lambda^k(G):R_h^*\alpha=\alpha,\ \ \forall h\in H\}.
\end{equation*}
Notice that the standard exterior differentiation operator $d$ on $\Lambda^k(G)$ is $H$-equivariant, \textit{i.e.} $d(R_h^*\alpha)=R_h^*(d\alpha)$. It thus restricts to an operator
\begin{equation*}
d:\Lambda^k(G^H)\rightarrow \Lambda^{k+1}(G^H),
\end{equation*}
defining the invariant cohomology groups 
\begin{equation*}
H^k(G^H;\R):=\frac{\mbox{Ker}(d:\Lambda^k(G^H)\rightarrow \Lambda^{k+1}(G^H))}{\mbox{Im}(d:\Lambda^{k-1}(G^H)\rightarrow \Lambda^k(G^H))}.
\end{equation*}
There is a natural relationship between the invariant cohomology of the pair $(G,H)$ and the cohomology of the manifold $G/H$, as follows. The projection $\pi:G\rightarrow G/H$ satisfies $\pi\circ R_h=\pi$. This implies that the pull-back operation induces injections
\begin{equation}\label{eq:homogeneouspullback}
\pi^*:\Lambda^k(G/H)\rightarrow \Lambda^k(G^H).
\end{equation}
Since $\pi^*$ commutes with $d$ it defines homomorphisms between the corresponding cohomology groups 
\begin{equation*}
\pi^*:H^k(G/H;\R)\rightarrow H^k(G^H;\R),\ \ \pi^*[\alpha]:=[\pi^*\alpha].
\end{equation*}
In the special case $k=1$, this map is an injection. Indeed, given $[\alpha]\in H^1(G/H;\R)$, assume $\pi^*[\alpha]=0$. Then $\pi^*\alpha\in\Lambda^1(G^H)$ is exact, \textit{i.e.} $\pi^*\alpha=df$ for some $f\in \Lambda^0(G/H)$. However it is clear that Equation \ref{eq:homogeneouspullback} is an isomorphism for $k=0$, \textit{i.e.} $f=\pi^*f'$ for some $f'\in\Lambda^0(G/H)$. Thus $\pi^*(\alpha-df')=0$ so $\alpha=df'$, \textit{i.e.} $[\alpha]=0$.

Now assume given a left action of $G$ on a manifold $M$. Choose an orbit $\mathcal{O}$ of this action. According to Lemma \ref{l:orbits}, $\mathcal{O}$ is a smooth submanifold of $M$. Choosing $x\in\mathcal{O}$ allows us to define the invariant cohomology of the pair $(G,G_x)$. Using the point $y=gx$ leads us instead to the invariant cohomology of the pair $(G,G_{y})$. We can use Lemma \ref{l:changepoint} to build isomorphisms between these groups. In this sense, these cohomology groups depend only on $\mathcal{O}$. It thus makes sense to look for a construction of cohomology groups which is independent of the choice of point. This can be done as follows.

Consider the set of smooth maps from $\mathcal{O}$ into the vector space $\Lambda^k(\mathfrak{g})$,
\begin{equation*}
\Lambda^k(\mathcal{O};\mathfrak{g}):=C^\infty(\mathcal{O},\Lambda^k(\mathfrak{g})).
\end{equation*}
Notice that, given $\alpha\in\Lambda^k(\mathcal{O};\mathfrak{g})$ and $v\in \mathfrak{g}$, we obtain by contraction an element
\begin{equation*}
i_v\alpha:=\alpha(v,\cdot,\dots,\cdot)\in\Lambda^{k-1}(\mathcal{O};\mathfrak{g}).
\end{equation*}
Likewise, given $v_1,\dots,v_k\in\mathfrak{g}$, iterated contractions define an element
\begin{equation*}
\alpha(v_1,\dots,v_k)\in \Lambda^0(\mathcal{O};\mathfrak{g})=C^\infty(\mathcal{O},\R).
\end{equation*}
For any given $v_0,\dots,v_k\in\mathfrak{g}$ it thus makes sense to define 
\begin{align}\label{eq:orbitd}
d\alpha(v_0,\dots,v_k)&:=\sum_{j=0}^k (-1)^j\tilde{v}_j\alpha(v_0,\dots,\hat{v}_j,\dots,v_k)\\
&\quad + \sum_{j<l}(-1)^{j+l}\alpha(-[v_j,v_l],v_0,\dots,\hat{v}_j,\dots,\hat{v}_l,\dots,v_k)\nonumber
\end{align} 
where $[\cdot,\cdot]$ denotes the Lie bracket on $\mathfrak{g}$, the superscript $\hat{}$ denotes an omitted term and $\tilde{v}_i$ denotes the fundamental vector field associated to $v_i$. One can check (or it follows from Proposition \ref{prop:2cohoms}) that $d\alpha\in\Lambda^{k+1}(\mathcal{O};\mathfrak{g})$ and that $d(d\alpha)=0$. We thus obtain cohomology groups 
\begin{equation*}
H^k(\mathcal{O};\mathfrak{g}):=\frac{\mbox{Ker}(d:\Lambda^k(\mathcal{O};\mathfrak{g})\rightarrow\Lambda^{k+1}(\mathcal{O};\mathfrak{g}))}{\mbox{Im}(d:\Lambda^{k-1}(\mathcal{O};\mathfrak{g})\rightarrow \Lambda^k(\mathcal{O};\mathfrak{g}))}.
\end{equation*}
Now recall that, for any $x\in\mathcal{O}$, Lemma \ref{l:orbits} defines a projection $\mathfrak{g}\rightarrow T_x\mathcal{O}$. Dually, this implies that there exist natural injections $\Lambda^k(T_x\mathcal{O})\rightarrow \Lambda^k(\mathfrak{g})$. We can use these to define injections
\begin{equation}\label{eq:orbitpullback}
\Lambda^k(\mathcal{O})\rightarrow \Lambda^k(\mathcal{O};\mathfrak{g}),\ \ \Lambda\mapsto\bar{\Lambda}. 
\end{equation} 
\begin{prop}\label{prop:2cohoms}
Given any $x\in\mathcal{O}$ and setting $H:=G_x$, there exists a canonical isomorphism $\Lambda^k(\mathcal{O};\mathfrak{g})\rightarrow \Lambda^k(G^H)$ leading to a commutative diagram
$$\begin{CD}\label{cd:pullback}
\Lambda^k(\mathcal{O}) @>j^*>> \Lambda^k(G/H)\\
@VVV   @VV{\pi^*}V\\
\Lambda^k(\mathcal{O};\mathfrak{g}) @>>> \Lambda^k(G^H)
\end{CD}$$
where the vertical arrow on the left denotes the map of Equation \ref{eq:orbitpullback}. This isomorphism also leads to a canonical isomorphism between the corresponding cohomology groups, \textit{i.e.} an isomorphism $H^k(\mathcal{O};\mathfrak{g})\rightarrow H^k(G^H;\R)$. In particular, $H^1(\mathcal{O};\R)$ can be canonically viewed as a subgroup of $H^1(\mathcal{O};\mathfrak{g})$.
\end{prop}
\proof{}
Fundamental vector fields provide an identification $T_eG\rightarrow T_gG$ for any $g\in G$, \textit{i.e.} a \textit{parallelization} of $G$. Using this parallelization we can identify the space $\Lambda^k(G)$ of all k-forms on $G$ with the space of smooth maps $G\rightarrow \Lambda^k(T_eG)$. Restricting this identification gives identifications
\begin{align*}
\Lambda^k(G^H) &\simeq \{\alpha:G\rightarrow \Lambda^k(T_eG):\alpha(gh)=\alpha(g),\ \forall h\in H\}\\
&\simeq C^\infty(G/H,\Lambda^k(T_eG))\\
&= C^\infty(\mathcal{O},\Lambda^k(\mathfrak{g}))\\
&= \Lambda^k(\mathcal{O};\mathfrak{g}).
\end{align*}
It is simple to check that, up to these identifications, the above diagram commutes.

Now choose $\alpha\in\Lambda^k(\mathcal{O};\mathfrak{g})$. Let $\alpha'$ denote the corresponding element of $\Lambda^k(G^H)$. As usual let us denote by $\tilde{v}$ the fundamental vector field generated by $v$. We now want to prove that $d(\alpha')=(d\alpha)'$, \textit{i.e.} that, for all $v_0,\dots,v_k\in T_eG$, 
\begin{equation}\label{eq:dcoincides}
d(\alpha')(\tilde{v}_0,\dots,\tilde{v}_k)=d\alpha(v_0,\dots,v_k).
\end{equation}
According to Equation \ref{eq:usuald}, we can calculate the left hand side using the usual bracket on $\fields(G)$. However, recall from Remark \ref{r:rightinvariant} that $[\tilde{v}_i,\tilde{v}_j]=-\widetilde{[v_i,v_j]}$. The change of sign here is cancelled by the choice of signs in Equation \ref{eq:orbitd}. This proves the claim on $d$, thus on the cohomology groups.

Clearly there also exists an identification $j^*:H^1(\mathcal{O};\R)\simeq H^1(G/H;\R)$. We can now use the injection $\pi^*:H^1(G/H;\R)\rightarrow H^1(G^H;\R)$ to prove the last claim.
\endproof



\subsection{The group of diffeomorphisms}\label{ss:diffM}
Let $\cdiffrD$ denote the set of diffeomorphisms of $\R^D$ with compact support, \textit{i.e.} those which coincide with the identity map $Id$ outside of a compact subset of $\R^D$. Composition of maps clearly yields a group structure on $\cdiffrD$. It is possible to endow $\cdiffrD$ with the structure of an infinite-dimensional Lie group in the sense of \cite{milnor:liegroups}. A local model is provided by the space $\cxrD$, endowed as in Section \ref{ss:distributions} with the structure of a topological vector space. More specifically, we can apply the construction outlined in Remark \ref{r:charts} below to build a local chart $\mathcal{U}$ for $\cdiffrD$ near the identity element $Id$. This yields by definition an isomorphism $T_{Id}\cdiffrD\simeq\cxrD$. We can then use right multiplication to build charts $\mathcal{U}_\phi:=\{u\circ\phi:u\in\mathcal{U}\}$ around any $\phi\in\cdiffrD$, leading to $T_\phi\cdiffrD\simeq\{X\circ\phi:X\in\cxrD\}$. Thoughout this article we will generally restrict our attention to the connected component of $\cdiffrD$ containing $Id$.

\begin{remark}
It may be useful to emphasize that defining charts on $\cdiffrD$ as above leads to the following interpretation of Equation \ref{eq:tflow}: $\phi_t$ is a smooth path on $\cdiffrD$ and $X_t\circ\phi_t\in T_{\phi_t}\cdiffrD$ is its tangent vector field.

\end{remark}

As usual one can define the Lie algebra to be the tangent space at $Id$. The Lie bracket $[\cdot,\cdot]_\mathfrak{g}$ on this space can then be defined as in Section \ref{ss:groupactions}, cf. \cite{milnor:liegroups}.

\begin{lemma}\label{l:adjointispushforward}
The adjoint representation of $\cdiffrD$ on $\cxrD$ coincides with the push-forward operation: $Ad_\phi(X)=\phi_*(X)$. 
Furthermore, the Lie bracket on $\cxrD$ induced by the Lie group structure on $\cdiffrD$ is the negative of the standard Lie bracket on vector fields.
\end{lemma}
\proof{}
Assume that $X$ integrates to $\phi_t\in\cdiffrD$. Then
$$Ad_\phi(X)=d/dt (\phi\circ\phi_t\circ\phi^{-1})_{|t=0}=\nabla\phi_{|\phi^{-1}}\cdot X_{|\phi^{-1}}=\phi_*(X).$$
As in Lemma \ref{l:bracketviaad} we can calculate the Lie bracket by differentiating the adjoint representation. Thus: 
$$[X,Y]_\mathfrak{g}=d/dt(Ad_{\phi_t}Y)_{|t=0}=d/dt(\phi_{t*}Y)_{|t=0}=-[X,Y].$$
\endproof

\begin{remark}\label{r:negativeliebracket} 
Lemma \ref{l:adjointispushforward} explains why the map of Remark \ref{r:tangentspacevsalgebra} is an algebra anti-homomorphism. Indeed, any left group action defines a homomorphism $G\rightarrow\diffm$, thus a homomorphism between the corresponding Lie algebras. However, we now see that the bracket used for $\fields(M)$ in Remark \ref{r:tangentspacevsalgebra} is the negative of the bracket induced by the Lie group structure. Lemma \ref{l:adjointispushforward} is also related to Remark \ref{r:rightinvariant}.
\end{remark}

\begin{remark}\label{r:charts}
A similar construction proves that for any compact (respectively, noncompact) manifold $M$ the group of diffeomorphisms $\diffm$ (respectively, $\cdiffm$) is an infinite-dimensional Lie group in the sense of \cite{milnor:liegroups}. Some care has to be exercised however in all these constructions, specifically in the definition of the local chart near $Id$. The naive choice
$$\xm\rightarrow \diffm,\ \ X\mapsto \phi_1,$$
where $\phi_1$ is the time $t=1$ diffeomorphism obtained by integrating $X$ to the flow $\phi_t$, is not possible as it does not cover an open neighbourhood of $Id$, cf. \cite{milnor:liegroups} Warning 1.6. Instead, the standard trick is to notice that diffeomorphisms near $Id$ are in a 1:1 relationship (via their graphs) with smooth submanifolds close to the diagonal $\Delta\subset M\times M$. These submanifolds can then be parametrized as follows. Assume $E\rightarrow M$ is a vector bundle over $M$. Let $Z$ denote its zero section and $U$ denote an open neighbourhood of $Z$. Assume one can find a diffeomorphism $\zeta:U\rightarrow M\times M$ sending $Z$ to $\Delta$. Then diffeomorphisms of $M$ near $Id$ correspond to smooth submanifolds of $E$ near $Z$, \textit{i.e.} smooth sections. For example, to construct a chart for diffeomorphisms close to $Id$ we would use $E:=TM$ setting $\zeta$ to be the Riemannian exponential map (with respect to a fixed metric on $M$).

Good choices of $E$ and $\zeta$ for $\diffm$ can yield as a by-product the fact that specific subgroups $G$ of $\diffm$ also admit Lie group structures such that the natural immersion $G\rightarrow\diffm$ is smooth. For example, to prove this fact for the subgroups of symplectomorphisms or Hamiltonian diffeomorphisms of a symplectic manifold $(M,\omega)$ (see Section \ref{ss:hamdiffs}) one can choose $E:=T^*M$ and the $\zeta$ defined by Weinstein's ``Lagrangian neighbourhood theorem", cf. \cite{weinstein:symplecticlagrangian} Section 6 or \cite{mcduffsalamon:book} Proposition 3.34.
\end{remark}


{\bf Acknowledgments.} The authors wish to thank L. Ambrosio, Y. Brenier, A. Fathi, E. Ghys, M. Loss and  C. Villani for fruitful conversations. They also thank B. Khesin, P. Lee and J. Lott for preliminary versions of their work \cite{khesinlee} and \cite{lott:somecalcs}, A. Weinstein for suggesting the relevance of the work \cite{marsdenweinstein:maxwellvlasov}, and the referee for stimulating comments. TP would also like to thank D. Burghelea, J. Ebert, D. Fox, A. Ghigi and D. Joyce for useful discussions, and his wife Lynda for her encouragement and support.

WG gratefully acknowledges the support provided by NSF grants DMS-03-54729 and DMS-06-00791. HKK gratefully acknowledges RA support provided by NSF grants DMS-03-54729 and DMS-06-00791. TP is grateful to the Georgia Institute of Technology, Imperial College and the University of Oxford for their hospitality during various stages of this project, with support provided by a NSF VIGRE fellowship (2003-2006), an EPSRC fellowship (2006-2007) and a Marie Curie EIF fellowship (2007-2009).


\bibliographystyle{amsplain}
\bibliography{geomeasure}
\end{document}